\documentclass[11pt]{article}
\usepackage{amssymb,latexsym,amsmath,amsbsy,amsthm,amsxtra,amsgen,dsfont,pdfsync,graphicx,bm,mathtools,tikz,longtable}
\usepackage{extpfeil}
\usepackage{mathrsfs}
\usepackage{bbm}
\usepackage{indentfirst}
\usepackage{cite}
\usepackage{mathrsfs} 
\usepackage{xcolor}
\usepackage{url}
\usepackage{authblk}
\usepackage{blindtext}
\usepackage{booktabs}
\usepackage{titling}
\usepackage{float, placeins, afterpage}
\usepackage{enumitem}
\predate{}
\postdate{}

\usetikzlibrary{decorations.pathreplacing}

\newfont{\msbmsm}{msbm10 at 8pt}

\def\Var{\textup{Var}}

\def\eps{\varepsilon}

\newtheorem{Theo}{Theorem}
\newtheorem{Lemma}[Theo]{Lemma}
\newtheorem{Cor}[Theo]{Corollary}
\newtheorem{Prop}[Theo]{Proposition}

\newtheorem{Rmk}[Theo]{Remark}

\newcommand{\qedwhite}{\hfill \ensuremath{\Box}}
\newcommand*\samethanks[1][\value{footnote}]{\footnotemark[#1]}

\title{Particle configurations for branching Brownian motion with an inhomogeneous branching rate}
\author{Jiaqi Liu\thanks{ Supported in part by NSF Grant DMS-1707953} \;and Jason Schweinsberg\samethanks}
\affil{\fontsize{9}{10.8}\itshape{Department of Mathematics, University of California San Diego, La Jolla, CA 92093-0112, USA}}
\date{}

\begin{document}

\maketitle
\begin{abstract}
Aiming to understand the distribution of fitness levels of individuals in a large population undergoing selection, we study the particle configurations of branching Brownian motion where each particle independently moves as Brownian motion with negative drift, particles can die or undergo dyadic fission, and the difference between the birth rate and the death rate is proportional to the particle's location. Under some assumptions, we obtain the limit in probability of the number of particles in any given interval and an explicit formula for the asymptotic empirical density of the fitness distribution. We show that after a sufficiently long time, the fitness distribution from the lowest to the highest fitness levels approximately evolves as a traveling wave with a profile which is asymptotically related to the Airy function. Our work complements the results in Roberts and Schweinsberg (2021), giving a fuller picture of the fitness distribution.
\end{abstract}

{\small MSC: Primary 60J80; secondary 92D15.

Keywords: branching Brownian motion; particle configurations; natural selection}

\section{Introduction}
To understand the evolution of populations undergoing selection, we study the distribution of fitness levels of individuals in the population. There is a well-known observation in the biology and physics literature, going back to the work of Tsimring, Levine and Kessler \cite{Tsimring}, that in a large population where various beneficial mutations compete for fixation simultaneously, the fitness distribution over time can be well approximated by a traveling wave. Roberts and Schweinsberg \cite{RS2020} used branching Brownian motion with an inhomogeneous branching rate to model a population undergoing selection. They showed rigorously that, in their model, the empirical distribution of the fitness levels of individuals is approximately Gaussian. Our work complements the results in \cite{RS2020}, giving a fuller picture of the fitness distribution.  Using the same model as in \cite{RS2020}, we show that after a sufficiently long time, the fitness distribution from the lowest to the highest fitness levels approximately evolves as a traveling wave.  The profile of this traveling wave is asymptotically equivalent to the profile obtained using nonrigorous methods in \cite{Tsimring, Cohen, Neher, Melissa}, which is expressed in terms of the Airy function. 

The most intuitive model of fitness is the fitness landscape, which is a mapping from the multidimensional genotype space to a real-valued fitness space. This model is constructed in a high-dimensional space where the number of dimensions is equal to the number of nucleotides in the genome. Each point in the genotype space corresponds to a particular genome, and each genome is assigned a fitness level.  Although the fitness landscape captures fully the relationship between genotypes and fitness, only limited quantitative analysis can be done in this model due to its high dimensional construction, see e.g. \cite{Szendro}. To study evolution on a smooth fitness landscape, Tsimring, Levine and Kessler \cite{Tsimring} introduced a model with a one-dimensional fitness space and characterized the population density as a function of time and the fitness level.  They found that the fitness distribution evolves over time like a traveling wave. Since then, the traveling wave model of fitness has been studied for more than two decades with different model assumptions, see e.g. \cite{Beerenwinkel, Desai, Durrett, DFisher, Good, Hallatschek, Kessler, Melissa, Neher, Park, Parksexual, Rouzine}. For a summary of the dynamical behavior of traveling wave fitness models in different settings, see \cite{Geyrhofer}.

In the mathematics literature, most of the work related to the dynamical behavior of fitness has been done under the framework of the Moran model, where the number of individuals in the population is $N$ at all times and individuals acquire beneficial mutations at a certain rate. Durrett and Mayberry \cite{Durrett} first rigorously established the traveling wave behavior of fitness when the selection rate is constant and the mutation rate is $N^{-\alpha}$ for $0<\alpha<1$. Schweinsberg~\cite{Schweinsberg1} considered slightly faster mutation rates and showed that the distribution of fitness has a Gaussian-like tail behavior, though it does not actually converge to a Gaussian distribution because at a typical time, the vast majority of individuals have the same number of mutations. The work \cite{Schweinsberg1} made rigorous the earlier work of Desai and Fisher \cite{Desai}.  In the work \cite{Durrett, Schweinsberg1, Desai}, mutations are sufficiently rare that one studies the process by considering the effects of individual mutations.  Yu, Etheridge and Cuthbertson \cite{Yu}, followed by Kelly \cite{Kelly} considered the case of strong selection and a very fast mutation rate and established upper and lower bounds for the rate at which the mean fitness of the population increases. 

When the rate of beneficial mutations is large but the selective advantage resulting from each mutation is small, each individual acquires many mutations with a small selective advantage, and an individual's fitness level will evolve like a random walk. After proper scaling, the fitness of each individual will move according to Brownian motion. It is therefore natural to model such populations using branching Brownian motion, an idea which goes back to \cite{brunet2006noisy, brunet2007effect}.  Mathematically rigorous work was later carried out in \cite{BBS2013, BBS2015, Maillard}. In these works, the branching rate is homogeneous in space but particles move according to one-dimensional Brownian motion on the positive real line and are killed at the origin.  The absorption at the origin models the deaths of individuals with low fitness.  Recently, Roberts and Schweinsberg \cite{RS2020} studied the evolution of a large population undergoing selection using branching Brownian motion with an inhomogeneous branching rate, so that particles with higher fitness are more likely to have offspring. We will work under the same setup and assumptions as \cite{RS2020}. 
We are interested in understanding the distribution of individual fitness levels from the least fit individuals to the most fit individuals, or in other words, the configuration of particles from the left edge to the right edge. We note that a detailed nonrigorous analysis of this model was provided in \cite{Neher}.  We note also the work \cite{DFisher, Melissa}, which aims to fill in the gap between the analysis in \cite{Desai} in which it is assumed that mutations are relatively rare, and the work in \cite{Neher}, which handles the case of very fast mutation rates.

Branching Brownian motion with a space-dependent branching rate was introduced by Harris and Harris \cite{HH2009}. In their model, a particle at location $y\in \mathbb{R}$ will split into two particles at rate $\beta |y|^p$, where $\beta>0$ and $p\in[0,2]$. They did not include the case $p>2$ because the process will explode in finite time if $p>2$. They studied the location of the right-most particle using martingales and the related spine changes of measure. They proved that for $p\in [0,2)$, the maximal displacement grows polynomially while for $p=2$, the maximal displacement grows exponentially. Later, Berestycki et al. \cite{BBHHR2015} studied the particle configurations in this model for all $p\in [0,2)$. By considering the large deviations probabilities for particles following certain rescaled paths, they obtained the logarithmic order for both the expected number and the typical number of particles whose rescaled trajectories follow paths in some set.  Tourniaire \cite{Tourniaire} considered branching Brownian motion with a space-dependent branching rate in which particles are killed at the origin, particles in $[0,1]$ branch at rate $\rho/2$, where $\rho > 1$, and particles in $(1, \infty)$ branch at rate $1/2$.  This process models the so-called semi-pushed traveling waves that can appear, for example, when a population invades a new habitat.

The particle configurations for branching Brownian motion with absorption at the origin are well understood.
Berestycki et al. \cite{BBS2015} proved, roughly speaking, that if the process settles into an equilibrium configuration with $N$ particles in total, then the density of particles near $y$ is proportional to $e^{-\sqrt{2}y}\sin(\sqrt{2}\pi y/\log N)$. This is very different from the behavior that we observe in our model, in which the density of particles near the origin where the bulk of particles are located follows approximately a Gaussian distribution, as shown in \cite{RS2020}.

\subsection{The model}

Because we aim to prove limit theorems that will be valid when the size of the population tends to infinity, we need to define a sequence of branching Brownian motion processes indexed by $n$. In the $n$-th process, each particle moves independently as one-dimensional Brownian motion with drift $-\rho_n$. A particle at position $x$ will die at rate $d_n(x)$ and branch into two particles at rate $b_n(x)$, where
\begin{equation}\label{betandef}
b_n(x)-d_n(x)=\beta_n x.
\end{equation}
Here, each particle corresponds to an individual in the population. The positions of particles represent  fitness levels of individuals, and the movement of particles illustrates changes in fitness levels over generations. Branching events represent births.  Note that (\ref{betandef}) indicates that the difference between the birth and death rates will be larger for particles with higher fitness.  As indicated in section 1.5 of \cite{RS2020}, this model has the potential to describe the evolution of populations undergoing selection under a fairly wide range of conditions.  It is therefore of interest to obtain a detailed and mathematically rigorous understanding of this model.

We assume that
\begin{equation}\label{assump1}
\lim_{n\rightarrow\infty}\frac{\rho_n^3}{\beta_n}=\infty,
\end{equation}
\begin{equation}\label{assump2}
\lim_{n\rightarrow\infty}\rho_n=0,
\end{equation}
and there exists $\alpha\in (0,1)$ such that
\begin{equation}\label{assump3}
d_n(x)\geq \alpha \;\text{for all}\;x\in \mathbb{R},\; n\in\mathbb{N}\quad \text{and}\quad b_n(x)\leq 1/\alpha\;\text{for all}\;x\leq 1/\beta_n, n\in\mathbb{N}.
\end{equation}
We assume that the conditions (\ref{assump1}), (\ref{assump2}) and (\ref{assump3}) hold true throughout the rest of this paper, even when they are not explicitly stated.

We will also make some assumptions on the initial configuration of particles at time zero.  First, we introduce some notation.  In the $n$-th process, we denote by $N_{t,n}$ the total number of particles at time $t$.  We also let $\mathcal{N}_{t,n}$ be the set of particles alive at time $t$, and $N_{t,n}(\mathcal{I})$ will denote the number of particles in the interval $\mathcal{I}$ at time $t$. The set of positions of particles at time $t$ is written as $\{X_{i,n}(t), i\in \mathcal{N}_{t,n}\}$. For $i\in\mathcal{N}_{t,n}$, we denote by $\{X_{i,n}(r), 0\leq r\leq t\}$ the past trajectory of the particle $i$ which is alive at time $t$. Denote the Airy function by
\[
Ai(x)=\frac{1}{\pi}\int_{0}^{\infty}\cos \bigg(\frac{y^3}{3}+xy\bigg)dy.
\]
The Airy function has an infinite number of zeros, all of which are negative. We denote the zeros of the Airy function by  $(\gamma_k)_{k=1}^{\infty}$ such that $\dots<\gamma_2<\gamma_1<0$. It is known \cite{NIST} that to three decimal places,
\begin{equation}\label{gamma1}
\gamma_1\approx -2.338.
\end{equation}
We define
\begin{equation}\label{L*}
L_n^*=\frac{\rho_n^2}{2\beta_n}, \qquad L_n^\dagger=-\frac{5\rho_n^2}{8\beta_n}.
\end{equation}
Roughly speaking, most particles will stay within $[L_n^\dagger, L_n^*]$. We call $L_n^*$ the right edge of the process and $L_n^\dagger$ the left edge. We define $L_n$ by
\begin{equation}\label{L}
L_n=\frac{\rho_n^2}{2\beta_n}-(2\beta_n)^{-1/3}\gamma_1,
\end{equation}
which is slightly larger than $L^*_n$ because $\gamma_1 < 0$.  We refer to $L_n$ as the right boundary because only rarely will particles exceed $L_n$ and we will often use truncation arguments in which particles are killed at $L_n$. Let
\begin{equation}\label{defY}
Y_n(t)=\sum_{i\in\mathcal{N}_{t,n}}e^{\rho_n X_{i,n}(t)},
\end{equation}
and
\begin{equation}\label{defZ}
Z_n(t)=\sum_{i\in\mathcal{N}_{t,n}}e^{\rho_n X_{i,n}(t)}Ai\big((2\beta_n)^{1/3}(L_n-X_{i,n}(t))+\gamma_1)\big)1_{\{X_{i,n}(t)<L_n\}}.
\end{equation}
It is explained in \cite{RS2020} that $Z_n(t)$ provides a natural measure of the ``size" of the process at time $t$. 

We make the following assumptions regarding the initial configuration expressed in terms of $Y_n(0)$ and $Z_n(0)$. We assume that
\begin{equation}\label{assumpY}
\rho_n^2e^{-\rho_nL_n}Y_n(0)\rightarrow_{p}0,
\end{equation}
where here and throughout the paper we use the notation $\rightarrow_p$ to denote convergence in probability as $n \rightarrow \infty$.
We also assume that for all $\eps>0$, there exists a $\delta>0$ such that for $n$ sufficiently large,
\begin{equation}\label{assumpZ}
P\bigg(\delta\frac{\beta_n^{1/3}}{\rho_n^3}e^{\rho_nL_n}\leq Z_n(0)\leq \frac{1}{\delta}\frac{\beta_n^{1/3}}{\rho_n^3}e^{\rho_nL_n}\bigg)>1-\eps.
\end{equation}
Roughly speaking, the assumption (\ref{assumpY}) requires that the ``size" of the process at early times will not be dominated by the descendants of a single particle in the initial configuration or by particles that are far from $L_n$ at time $0$. The assumption (\ref{assumpZ}) is roughly saying that the ``size" of the initial configuration will be around $\beta_n^{1/3}e^{\rho_nL_n}/\rho_n^3$.

\subsection{Main results}
We will first introduce some notation that will be used throughout the paper. For two sequences of positive numbers $(x_n)_{n=1}^{\infty}$ and $(y_n)_{n=1}^{\infty}$, if $x_n/y_n$ is bounded above by a positive constant, we write $x_n\lesssim y_n$ and if $\lim_{n\rightarrow\infty}x_n/y_n=0$, we write $x_n\ll y_n$. We define $x_n\gtrsim y_n$ and $x_n\gg y_n$ correspondingly. Moreover, the notation $x_n\asymp y_n$ means that $x_n/y_n$ is bounded above and below by positive constants, and the notation $x_n\sim y_n$ means that $\lim_{n \rightarrow \infty} x_n/y_n = 1.$ We write $x_n=O(y_n)$ if the sequence $(x_n/y_n)_{n=1}^{\infty}$ is bounded and $x_n=o(y_n)$ if $\lim_{n\rightarrow\infty}x_n/y_n=0$. 

Before stating our new results, we briefly recall the main results that Roberts and Schweinsberg \cite{RS2020} established for this model.  Under assumptions (\ref{assump1}), (\ref{assump2}), (\ref{assump3}), (\ref{assumpY}) and (\ref{assumpZ}), they showed that if $\rho_n^{2/3}/\beta_n^{8/9}\ll t_n-\rho_n/\beta_n\lesssim \rho_n/\beta_n$, then most particles are near the origin at time $t_n$ and the scaled empirical distribution of particles at time $t_n$ is Gaussian. More precisely, define the random probability measure which represents the empirical distribution of the particle locations at time $t_n$, scaled in space, to be
\begin{equation}\label{Gaussian}
\zeta_n(t_n)=\frac{1}{N_{t_n,n}}\sum_{i\in\mathcal{N}_{t_n,n}}\delta_{X_{i,n}(t_n)\sqrt{\beta_n/\rho_n}}.
\end{equation}
They showed as $n\rightarrow\infty$, that the random measures $\zeta_n(t_n)$ converge weakly to the standard normal distribution in the Polish space of probability measures on $\mathbb{R}$ equipped with the weak topology. From the scaling in (\ref{Gaussian}), this result implies that the empirical distribution of particle locations at time $t$ is approximately normal with mean $0$ and variance $\rho_n/\beta_n$.  In particular, this result describes the configuration of particles whose distance to the origin is $O(\sqrt{\rho_n/\beta_n})$.

Roberts and Schweinsberg \cite{RS2020} also provided an explicit characterization of the empirical distribution of particles close to the right edge. They considered the empirical measure where a particle at $x$ is weighted by $e^{\rho_nx}$. Define the random probability measure
\begin{equation}\label{xindef}
\xi_n(t_n)= \frac{1}{Y_n(t_n)} \sum_{i\in\mathcal{N}_{t_n,n}}e^{\rho_nX_{i,n}(t_n)}\delta_{(2\beta_n)^{-1/3}(L_n-X_{i,n}(t_n))}
\end{equation}
Thus, particles with a higher fitness level will contribute more to $\xi_n(t_n)$. Let $\mu$ be the probability measure on $(0,\infty)$ with probability density function
\[
h(y)=\frac{Ai(y+\gamma_1)}{\int_{0}^{\infty}Ai(z+\gamma_1)dz}.
\]
Roberts and Schweinsberg proved that under assumptions (\ref{assump1}), (\ref{assump2}), (\ref{assump3}), (\ref{assumpY}) and (\ref{assumpZ}), if 
\[
\beta_n^{-2/3}\log^{1/3} \bigg(\frac{\rho_n}{\beta_n^{1/3}}\bigg)\ll t_n\lesssim \frac{\rho_n}{\beta_n},
\]
then as $n\rightarrow\infty$, we have
\begin{equation}\label{RSrightedge}
\xi_n(t_n)\Rightarrow \mu,
\end{equation}
where $\Rightarrow$ refers to weak convergence in the Polish space of probability measures on $\mathbb{R}$ equipped with the weak topology.  From the scaling in (\ref{xindef}), we see that this convergence result describes the configuration of particles whose distance from the right edge $L_n^*$ is $O(\beta_n^{-1/3})$.

Our goal in this paper is to obtain a fuller understanding of the particle configurations from the left edge $L_n^\dagger$ to the right edge $L^*_n$. In other words, for this model, we aim to characterize the long-run empirical distribution of the fitness levels of individuals in a large population.

Consider a sequence of intervals $\{[a_n,b_n]\}_{n=1}^{\infty}$, where $-\infty \leq a_n < b_n \leq \infty$, satisfying the following three conditions:
\begin{align}
b_n-a_n &\gg1 \label{b-a} \\
L^*_n-a_n &\gg \beta_n^{-1/3} \label{an} \\
b_n-L_n^\dagger &\gg \beta_n^{-1/3} \label{bn}.
\end{align}
We are interested in the number of particles in the intervals $[a_n, b_n]$.
We include the conditions (\ref{an}) and (\ref{bn}) because we do not expect our results to describe the configuration of particles which are within distance $O(\beta_n^{-1/3})$ of the right edge $L^*_n$ or the left edge $L_n^\dagger$. Also, the particles within $O(\beta_n^{-1/3})$ distance of the right edge were studied in Theorem 1.2 in \cite{RS2020}. Define
\begin{equation}\label{defzn}
z_n= \begin{cases} 
       a_n&\qquad \text{if}\;a_n\geq 0,\\
       0&\qquad\text{if}\; a_n<0,\; b_n>0,\\
       b_n&\qquad \text{if}\;b_n\leq 0.
   \end{cases}
\end{equation}
Note that $z_n\in (L_n^\dagger,L_n^*)$, and the restrictions (\ref{an}) and (\ref{bn}) are equivalent to 
\begin{equation}\label{zn}
L_n^*-z_n\gg \beta_n^{-1/3}, \qquad z_n-L_n^\dagger\gg \beta_n^{-1/3}.
\end{equation}
Later, we will see that the asymptotic density of the number of particles in $[a_n,b_n]$ reaches its maximum at $z_n$. As a result, the number of particles near $z_n$ dominates the total number of particles in $[a_n,b_n]$. 

For every $n$, we will define two important functions in the domain $(-\infty,L_n^*]$.  First, we let
\begin{equation}\label{t_z}
t_{n}(y)=\sqrt{\frac{2}{\beta_n}(L^*_n-y)}.
\end{equation}
We will later see that particles near $y$ are most likely descended from ancestors that were near the right edge approximately $t_n(y)$ time units in the past.  For every $n$, we observe that $t_n(y)$ is a decreasing function of $y$.  We have $t_n(0)=\rho_n/\beta_n$. If $L^*_n-z_n\gg \beta_n^{-1/3}$, then
 \begin{equation}\label{t(z)>beta-3/2}
 t_n(z_n)\gg \beta_n^{-2/3}.
 \end{equation}
Also, for $y \in (-\infty, L_n^*]$, we define
\begin{equation}\label{g}
g_n(y)=\rho_n(L^*_n-y)-\frac{2\sqrt{2\beta_n}}{3}(L^*_n-y)^{3/2}.
\end{equation}
We will see shortly that in the long-run, the number of particles located near $y$ is roughly proportional to $e^{g_n(y)}$.  Note that $g_n(y)$ is decreasing in $[0,L_n^*]$ and increasing in $(-\infty,0]$.  The functions $g_n(y)$ and $t_n(y)$ were previously obtained in \cite{RS2020}.

We now state our main result, which describes the configuration of particles from the left edge to the right edge. 

\begin{Theo}\label{ThmGlobal}
Suppose assumptions (\ref{assump1}), (\ref{assump2}), (\ref{assump3}), (\ref{assumpY}) and (\ref{assumpZ}) hold. For every sequence of intervals $\{[a_n,b_n]\}_{n=1}^{\infty}$ satisfying (\ref{b-a})-(\ref{bn}), define $z_n$ according to (\ref{defzn}). If
\begin{equation}\label{thmassump1}
\frac{\rho_n^{2/3}}{\beta_n^{8/9}}\ll t_n-t_{n}(z_n)\ll \frac{\rho_n}{\beta_n},
\end{equation}
then as $n\rightarrow\infty$,
\begin{equation}\label{thm1}
N_{t_n,n}\big([a_n,b_n]\big)\bigg/\bigg(\frac{1}{Ai'(\gamma_1)^2}Z_n(0)e^{-\rho_n L_n^*}\int_{[a_n,b_n]\cap (-\infty,L_n^*]}\frac{1}{\sqrt{2\pi t_{n}(y)}}e^{g_n(y)}dy\bigg)\rightarrow_p 1.
\end{equation}
If 
\begin{equation}\label{thmassump2}
 t_n-t_{n}(z_n)\asymp \frac{\rho_n}{\beta_n},
\end{equation}
then as $n\rightarrow\infty$,
\begin{equation}\label{thm2}
N_{t_n,n}\big([a_n,b_n]\big)\bigg/\bigg(\frac{1}{Ai'(\gamma_1)^2}Z_n\big(t_n-t_n(z_n)\big)e^{-\rho_n L_n^*}\int_{[a_n,b_n]\cap (-\infty,L_n^*]}\frac{1}{\sqrt{2\pi t_{n}(y)}}e^{g_n(y)}dy\bigg)\rightarrow_p 1.
\end{equation}
\end{Theo}

Theorem \ref{ThmGlobal} describes the number of particles in any given interval in the long run. The randomness is characterized by the stochastic process $\{Z_n(t),t\geq 0\}$, which measures how the overall ``size" of the process changes over time. The deterministic part has a density formula proportional to $e^{g_n(y)}/\sqrt{2\pi t_n(y)}$. To be more precise, shortly after time $t_n(z_n)$, the number of particles in the interval $[a_n, b_n]$ depends on the initial configuration of particles through the value of $Z_n(0)$. For much later times $t_n$, when $t_n - t_n(z_n)$ is of the order $\rho_n/\beta_n$, the number of particles in the interval $[a_n, b_n]$ depends on $Z(t_n - t_n(z_n))$, which is the ``size" of the process
$t_n(z_n)$ time units in the past. Here $z_n$ is the point where the density of the number of particles in $[a_n,b_n]$ is maximized. Later it will be shown in the proof that the number of particles in any interval $[a_n,b_n]$ is dominated by the number of particles that are close to $z_n$. The proof of Theorem~\ref{ThmGlobal} indeed shows that most of the particles in the interval $[a_n,b_n]$ at time $t_n$ are descendants of particles that are close to the right edge $t_n(z_n)$ time units in the past. This also explains why the number of particles in $[a_n,b_n]$ depends on $Z(t_n-t_n(z_n))$. 


\begin{Cor}\label{Tight}
Suppose assumptions (\ref{assump1}), (\ref{assump2}), (\ref{assump3}), (\ref{assumpY}) and (\ref{assumpZ}) hold. For every sequence of intervals $\{[a_n,b_n]\}_{n=1}^{\infty}$ satisfying (\ref{b-a})-(\ref{bn}), define $z_n$ according to (\ref{defzn}). Suppose
\begin{equation}\label{tightassumptn}
\frac{\rho_n^{2/3}}{\beta_n^{8/9}}\ll t_n-\max\big\{t_{n}(z_n),t_n(0)\big\}\lesssim \frac{\rho_n}{\beta_n}.
\end{equation}
For $y\in (-\infty, L_n^*]$, define 
\[
f_n(y)=\frac{1}{\sqrt{2\pi t_{n}(y)}}e^{g_n(y)-\rho_n^3/6\beta_n}.
\]
The sequence
\begin{equation*}
(D_n)_{n=1}^{\infty}:=\left\{\frac{N_{t_n,n}\big([a_n,b_n]\big)}{N_{t_n}}\bigg/\bigg(\int_{[a_n,b_n]\cap (-\infty,L^*]}f_n(y)dy\bigg)\right\}_{n=1}^{\infty}
\end{equation*}
is tight. If $0\in [a_n,b_n]$ for all $n$, then $D_n$ converges to $1$ in probability as $n\rightarrow\infty$.
\end{Cor}
Corollary \ref{Tight} shows that the ratio of the number of particles in any given interval to the total number of particles is comparable to the integral of $f_n(y)$ over the given interval. We can therefore regard $f_n(y)$ as the density of the limiting empirical distribution of the process, or in other words, the asymptotic empirical density of the fitness levels of individuals.

To understand the connection with results on traveling waves in the biology and physics literature, consider a translation of the model where each particle independently moves as standard Brownian motion without drift. A particle at  location $y$ can either die or split into two particles, and the difference between the birth rate and the death rate is $\beta_n(y-\rho_nt)$. Corollary \ref{Tight} shows that after a sufficiently long time, the empirical density of individual fitness levels is
\[
f^*_n(t, y)=f_n(y-\rho_nt)=\frac{1}{\sqrt{2\pi t_n(y-\rho_nt)}}e^{g_n(y-\rho_nt)-\rho_n^3/6\beta_n}, \quad\text{for}\; y\in (L_n^\dagger+\rho_nt,L_n^*+\rho_nt),
\]
which is a traveling wave with profile $f_n(y)$.

The asymptotic empirical density $f_n(y)$ is closely related to the Airy function. For $y<L_n^*$, define
\begin{equation}\label{airydensity}
f^A_n(y)=(2\beta_n)^{1/3}e^{-\rho_n y+\rho_n^3/3\beta_n}Ai\big((2\beta_n)^{1/3}(L_n^*-y)\big).
\end{equation}
According to (2.45) in \cite{Airy}, 
\begin{equation}\label{airyasymp}
\lim_{x\rightarrow\infty}2\sqrt{\pi}x^{1/4} e^{(2/3)x^{3/2}}Ai(x) =1.
\end{equation}
Therefore, if $L_n^*-y_n \gg\beta_n^{-1/3}$, then as $n\rightarrow\infty$,
\[
f_n(y_n)\sim f_n^A(y_n).
\]
Note that the restriction $L_n^*-y_n\gg\beta_n^{-1/3}$ is consistent with our requirement (\ref{an}) on the interval.
The idea that the traveling wave profile should have a shape given by the Airy function goes back to the early work of Tsimring, Levine, and Kessler \cite{Tsimring}, and this Airy shape also appears, for example, in \cite{Cohen, Neher, Melissa}.  Theorem \ref{ThmGlobal} and Corollary \ref{Tight} therefore provide rigorous justification for this result in the biology and physics literature.

We also observe that the shape of $f_n(y)$ near $0$ is very much like the Gaussian density function with standard deviation $\sqrt{\rho_n/\beta_n}$.  Let
\begin{equation}\label{fnG}
f^G_n(y)=\frac{1}{\sqrt{2\pi\rho_n/\beta_n}}\exp\bigg(-\frac{\beta_ny^2}{2\rho_n}\bigg).
\end{equation}
As noted in \cite{RS2020}, the Taylor expansions of $g_n(y)$ and $t_n(y)$ around $0$ give
\begin{equation*}
g_n(y)\approx \frac{\rho_n^3}{6\beta_n}-\frac{\beta_ny^2}{2\rho_n}, \qquad t_n(y)\approx\sqrt{\frac{\rho_n}{\beta_n}}.
\end{equation*}
Therefore, the asymptotic empirical density $f_n(y)$ can be approximated by the Gaussian density formula $f^G_n(y)$. This is consistent with Theorem 1.1 in \cite{RS2020}. 

\begin{figure}
\includegraphics[width=14cm]{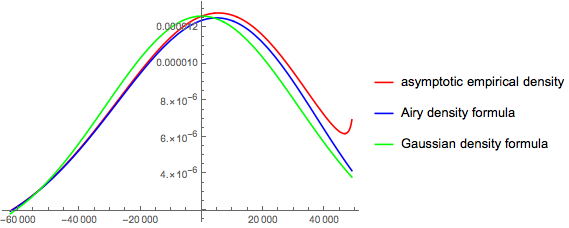}
\centering
\caption{Graph of the asymptotic empirical density, Airy density formula and Gaussian density formula when $\rho_n=10^{-4}$ and $\beta_n=10^{-13}$}
\label{comparison}
\end{figure}

Figure \ref{comparison} illustrates the graphs of the asymptotic empirical density $f_n(y)$, the Airy density formula $f^A_n(y)$ and the Gaussian density formula $f^G_n(y)$ from $L^\dagger_n$ to $L_n^*$ when $\rho_n=10^{-4}$ and $\beta_n=10^{-13}$. We see that all three functions have similar shapes. The asymptotic empirical density is very close to the Airy density formula in the bulk, especially in the negative real line where $y$ is far away from the right boundary $L^*$. 

Let $M_{t,n}=\max\{X_{i,n}(t),i\in\mathcal{N}_{t,n}\}$ be the position of the right-most particle at time $t$ and $m_{t,n}=\min\{X_{i,n}(t),i\in\mathcal{N}_{t,n}\}$ be the position of the left-most particle at time $t$. Propositions \ref{Corollaryr} and \ref{Corollaryl} show that, under certain assumptions, with high probability the right-most particle is close to $L_n^*$ and the left-most particle is close to $L_n^\dagger$. This explains why we are able to refer to $L_n^*$ as the right edge and $L_n^\dagger$ as the left edge of the process.
\begin{Prop}\label{Corollaryr}
Suppose assumptions (\ref{assump1}), (\ref{assump2}), (\ref{assump3}), (\ref{assumpY}) and (\ref{assumpZ}) hold and $(t_n)_{n=1}^{\infty}$ satisfies
\begin{equation}\label{corassumpr}
\beta_{n}^{-2/3}\log^{1/3}\bigg(\frac{\rho_n}{\beta_n^{1/3}}\bigg)\ll t_n\lesssim \frac{\rho_n}{\beta_n}.
\end{equation}
For any positive constant $C_1$, we have
\begin{equation}\label{rlower}
\lim_{n\rightarrow\infty}P\bigg(M_{t_n,n}\geq L_n-\frac{C_1}{\beta_n^{1/3}}\bigg)=1.
\end{equation}
If in addition, the birth rate $b_n(x)$ is non-decreasing and the death rate $d_n(x)$ is non-increasing, then for any constant $C_2\in\mathbb{R}$,
\begin{equation}\label{rupper}
\lim_{n\rightarrow\infty}P\bigg(M_{t_n,n}\leq L_n+\frac{C_2}{\rho_n}\bigg)=1.
\end{equation}
Therefore, we have as $n\rightarrow\infty$,
\begin{equation}\label{corright}
\frac{M_{t_n,n}}{L_n^*}\rightarrow_{p}1.
\end{equation}
\end{Prop}

Define
\begin{equation}
\bar{L}_n=-\frac{5}{8}\frac{\rho_n^2}{\beta_n}+2(2\beta_n)^{-1/3}\gamma_1,
\end{equation}
which is slightly smaller than $L_n^\dagger$. The following proposition shows that $\bar{L}_n$ is the approximate position of the left-most particle.
\begin{Prop}\label{Corollaryl}
Suppose assumptions (\ref{assump1}), (\ref{assump2}), (\ref{assump3}), (\ref{assumpY}) and (\ref{assumpZ}) hold and $(t_n)_{n=1}^{\infty}$ satisfies
\begin{equation}\label{corassumpl}
\frac{\rho_n^{2/3}}{\beta_n^{8/9}}\ll t_n-t_{n}(\bar{L}_n)\lesssim \frac{\rho_n}{\beta_n}.
\end{equation}
For any $\kappa>0$, there exists a positive constant $C_3$ such that for $n$ large enough,
\begin{equation}\label{lupper}
P\bigg(m_{t_n,n}\leq \bar{L}_n+\frac{C_3}{\beta_n^{1/3}}\bigg)>1-\kappa.
\end{equation}
If in addition, the birth rate $b_n(x)$ is non-decreasing and the death rate $d_n(x)$ is non-increasing, then for any $\kappa>0$, there exists a positive constant $C_4$ such that for $n$ large enough,
\begin{equation}\label{llower}
P\bigg(m_{t_n,n}\geq \bar{L}_n-\frac{C_4}{\rho_n}\bigg)>1-\kappa.
\end{equation}
Therefore, we have as $n\rightarrow\infty$,
\begin{equation}\label{corleft}
\frac{m_{t_n,n}}{L_n^\dagger}\rightarrow_{p}1.
\end{equation}
\end{Prop}

\subsection{Heuristics for understanding the density formula}


The large deviations heuristics proposed in \cite{BBHHR2015} inspired the derivation of the functions $g_n(z_n)$ and $t_n(z_n)$ in \cite{RS2020}, although the techniques used in \cite{BBHHR2015} are not sufficient to derive the exact asymptotic rate of the number of particles as we did in Theorem \ref{ThmGlobal}. Roberts and Schweinsberg \cite{RS2020} conjectured that the number of particles near $z_n$ in the long run is proportional to $e^{g_n(z_n)}$ and proved this conjecture when $|z_n|\lesssim\sqrt{\rho_n/\beta_n}$. Because these heuristics are essential for understanding the behavior of the process and the main strategy of the proof, we will briefly recall their calculations. 

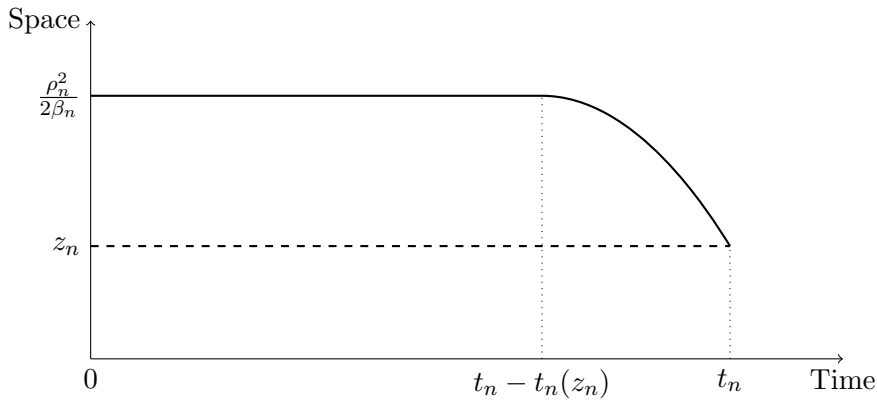
\begin{figure}[h]
\centering
\begin{tikzpicture}
\draw[->] (0,0) -- (10,0) node[anchor=north] {Time};
\draw	(0,0) node[anchor=north] {0}
		(6,0) node[anchor=north] {$t_n-t_n(z_n)$}
		(8.5,0) node[anchor=north] {$t_n$};

\draw[->] (0,0) -- (0,4.5) node[anchor=east] {Space};
\draw   (0,1.5) node[anchor=east]{$z_n$}
            (0,3.5) node[anchor=east]{$\frac{\rho_n^2}{2\beta_n}$};

\draw[dotted] (6,0) -- (6,3.5);
\draw[dotted] (8.5,0) -- (8.5,1.5);
\draw[thick] (0,3.5) -- (6,3.5) parabola[bend at start](8.5,1.5);
\draw[thick,dashed] (0,1.5) -- (8.5,1.5);
\end{tikzpicture}
\caption{Trajectory of $f_n^{z_n}$}
\label{figure}
\end{figure}

For every $n$, consider a large time $t_n$ and a path $f_n:[0,t_n]\rightarrow\mathbb{R}$. By Schilder's theorem and the many-to-one lemma, if the process starts with one particle at $f_n(0)$, the expected number of particles that stay ``close'' to $f_n$ during $[0,t_n]$ is approximately
\begin{equation}\label{Schilder}
\exp\bigg(\int_0^{t_n}\Big(\beta_n f_n(u)-\frac{1}{2}(f'_n(u)+\rho_n)^2\Big)du\bigg).
\end{equation}
Note that if $f_n(u)\equiv \rho_n^2/2\beta_n$, then the integrand is $0$. Thus, the right-most particle should stay close to $\rho_n^2/2\beta_n$, which is the right edge $L^*_n$. We next consider the optimal trajectory $f_n^{z_n}$ followed by particles that are near $z_n$ at time $t_n$. This path is optimal in the sense that particles which end up near $z_n$ must follow this trajectory to achieve the maximum almost sure growth rate. According to Theorem 7 in \cite{BBHHR2015},  there exists a cutoff time $t_n(z_n)$ such that the optimal path will follow the trajectory of the right-most particle up to some time $t_n-t_n(z_n)$ and then move towards $z_n$ by following a path that satisfies a certain differential equation. Therefore, $f_n^{z_n}$ satisfies 
\[
f_n^{z_n}(u)=\rho_n^2/2\beta_n \;\text{for}\; u\in [0,t_n-t_n(z_n)],
\]
\[
(f_n^{z_n})''(u)=-\beta_n\; \text{for} \;u\in [t_n-t_n(z_n),t_n],
\]
\[
 f_n^{z_n}(t_n)=z_n.
\]
Solving the above equations, we get the expression (\ref{t_z}) for $t_n(z_n)$, as well as the formula
$$f_{z_n}(u) = \frac{\rho_n^2}{2 \beta_n} - \frac{\beta_n}{2} \big(u - (t_n - t_n(z_n)) \big)^2 \qquad \mbox{for }u \in [t_n - t_n(z_n), t_n].$$
Combining these results with (\ref{Schilder}), we get that the number of particles near $z_n$ at time $t_n$ is approximately 
\[
\exp(g_n(z_n))=\exp\bigg(\int_0^{t_n}\Big(\beta f_n^{z_n}(u)-\frac{1}{2}\big((f_n^{z_n})'(u)+\rho_n\big)^2\Big)du\bigg),
\]
which gives (\ref{g}) for all $z_n$. Figure \ref{figure} is an illustration of the trajectory of $f_n^{z_n}$.

It is also worth mentioning that the expressions for the left edge $L_n^\dagger$ and the right edge $L_n^*$ emerge from $g_n(z_n)$. Solving $g_n(z_n)=0$, we get two solutions
\[
z_n=-\frac{5\rho_n^2}{8\beta_n}, \qquad z_n=\frac{\rho_n^2}{2\beta_n},
\]
which correspond to $L_n^\dagger$ and $L_n^*$ respectively.

These heuristics also explain why we need assumptions (\ref{assump1}) and (\ref{assump2}). Note that the Taylor expansion of $e^{g(z)}$ is proportional to the Gaussian density with mean $0$ and variance $\rho_n/\beta_n$. The standard deviation of the Gaussian distribution should be much smaller than the right-most position, which leads to (\ref{assump1}). Moreover, the branching rate should be small around the right edge $\rho_n^2/2\beta$, which leads to (\ref{assump2}).

We can also get some insight into the formula for $g_n(z_n)$ by considering the density for the branching Brownian motion process.  We denote by $p_{t}^n(x,y)$ the density at location $y$ and time $t$ for the process when it starts from a single particle at $x$ at time $0$. This means that if there is a single particle located at $x$ at time $0$, then the expected number of particles in the measurable set $U$ at time $t$ is
\[
\int_U p_{t}(x,y)dy.
\]According to formula (2.11) in \cite{RS2020}, by the many-to-one lemma,
\begin{equation}\label{density2.8}
p_{t}^n(x,y)=\frac{1}{\sqrt{2\pi t}}\exp\bigg(\rho_n x-\rho_n y-\frac{(x-y)^2}{2t}-\frac{\rho_n^2 t}{2}+\frac{\beta_n(x+y)t}{2}+\frac{\beta_n^2t^3}{24}\bigg).
\end{equation}
In general, the density formula $p_t^n(x,y)$ may not yield a good estimate of the actual number of particles near $y$ at time $t$ because the expectation may be dominated by rare events in which one particle drifts very far to the right and generates a large number of descendants around $y$.  For example, if $x = \rho_n^2/2 \beta_n$ and $t_n$ is much larger than $t_n(z_n)$, then $p_{t_n}(x,z_n)$ is dominated by particles whose trajectories start out by going above the trajectory $f_n^{z_n}$ that is pictured in Figure \ref{figure}.  Because it is rare that any particle follows this trajectory, $p_{t_n}(x,z_n)$ will overestimate the actual number of particles near $z_n$.  On the other hand, the density formula reflects well how the process evolves between times $t_n - t_n(z_n)$ and time $t_n$.  In particular, it is possible to calculate that for $z_n < L_n^*$, we have the exact equality
$$p^n_{t_n(z_n)}(L^*, z_n) = \frac{1}{\sqrt{2 \pi t_n(z_n)}} \: e^{g_n(z_n)}.$$

One way to obtain a density which approximates well the actual number of particles near $y$ at larger times is to use a truncation argument.  Suppose we modify our branching Brownian motion process by killing any particle that reaches $L_n$.  Lemma 2.5 in \cite{RS2020} shows that if $t$ is large enough, then as a function of $y$, the density for this modified process is approximately proportional to
\begin{equation}\label{RSAiry}
e^{-\rho_n y} Ai((2 \beta_n)^{1/3}(L_n - y) + \gamma_1),
\end{equation}
which closely resembles (\ref{airydensity}).

To prove Theorem \ref{ThmGlobal}, we use first and second moment estimates.  We note that $N_{t_n,n}([a_n,b_n])$ is dominated by the number of particles that are close to $z_n$.  We therefore look back to time $t_n - t_n(z_n)$.  According to the large deviations estimates, particles that are near $z_n$ at time $t_n$ will be descended from particles that are near the right edge at time $t_n - t_n(z_n)$, and we can use (\ref{RSrightedge}) to understand what the configuration of particles near the right edge looks like at time $t_n-t_n(z_n)$.  Conditional on the configuration at time $t_n - t_n(z_n)$, we can estimate the expected value of $N_{t_n,n}([a_n, b_n])$ by expressing $p_{t_n}^n(x,z_n)$ in terms of $g_n(z_n)$ and some other error terms that can be controlled.  To bound the variance of $N_{t_n,n}([a_n,b_n])$, we use a truncated second moment estimate in which we kill particles when they reach $L_n$.  We show that the dominant contribution to $N_{t_n,n}([a_n, b_n])$ comes from particles which are close to $L_n$ at time $t_n-t_n(z_n)$ and do not hit the right boundary $L_n$ during the time interval $[t_n-t_n(z_n),t_n]$. 

\subsection{Connections with the traveling wave literature}

According to the discussion after Corollary \ref{Tight}, the asymptotic empirical density $f_n(y)$ is closely related to the Airy density formula $f_n^A(y)$ in (\ref{airydensity}).  The Airy function for the shape of the traveling wave was previously derived nonrigorously in \cite{Tsimring, Cohen, Melissa, Neher}.  Here we review this nonrigorous derivation and explain in more detail the connections with our work.

Suppose there are $N$ individuals in a population. Each individual is subject to new mutations at rate $\mu$, and the selective advantage $s$ of each mutation is random and has a distribution with probability density function $\rho(s)$. Let $q(x,t)$ be the ``density" of particles with fitness $x$ at time $t$. Define $m(t)$ to be the average fitness at time $t$, with $m(0)=0$. Let $\nu(s)=\mu\rho(s)$. Then (see equation (4) in \cite{Melissa} or equation (2) in the supplementary information to \cite{Neher}), $q(x,t)$ can be approximated by the equation
\begin{equation}\label{SPDE}
\frac{\partial}{\partial t}q(x,t)=(x-m(t))q(x,t)+\int \big(q(x-s,t)-q(x,t)\big)\nu(s)ds+\sqrt{
\frac{q(x,t)}{N}}\eta(x,t),
\end{equation}
where $\eta$ is Gaussian white noise.  Note that the first term models selection, the second term models mutation, and the third term models the noise from the randomness in the births and deaths.  One can look for traveling wave solutions to (\ref{SPDE}) of the form
\[
q(x,t)=\omega(x-vt),
\]
where $v=m(t)/t$ denotes the velocity of the traveling wave, or the average rate at which the mean fitness of the population changes. Writing $y=x-vt$ for the relative fitness and neglecting the noise term, equation (\ref{SPDE}) becomes
\[
-v\omega'(y)=y\omega(y)+\int \big(\omega(y-s)-\omega(y)\big)\nu(s)ds.
\]
Assuming the selective advantage $s$ is small, we can use a Taylor expansion to approximate 
\[
\omega(y-s)-\omega(y)\approx -s\omega'(y)+\frac{1}{2} s^2\omega''(y).
\]
which leads to
\[
-v\omega'(y)=y\omega(y)-\mu E[s]\omega'(y)+\frac{1}{2} \mu E[s^2]\omega''(y).
\]
Now letting $D = \mu E[s^2]/2$ and $\sigma^2 = v - \mu E[s]$, we get
\begin{equation}\label{biodiff}
D\omega''(y)+\sigma^2\omega'(y)+y\omega(y)=0.
\end{equation}
Recalling that the Airy function satisfies the differential equation $Ai''(y) = yAi(y)$, we see that a solution to this differential equation is given by
\begin{equation}\label{biosol}
\omega(y)=Ce^{-\sigma^2 y/2D}Ai\Big(\frac{\sigma^4}{4D^{4/3}}-\frac{y}{D^{1/3}}\Big),
\end{equation}
which matches equation (6) in the supplementary information to \cite{Neher} and equations (25) and (28) in \cite{Melissa}.
We note that the above equation will lead to a solution which takes negative values for some large $y$, which is impossible for a fitness distribution. To avoid this, it is assumed in \cite{Neher, Melissa, Tsimring} that there is a cutoff value $y_{cut}$, which can be understood as the maximum fitness of the individuals in the population, such that $\omega(y)$ is given by (\ref{biosol}) for $y < y_{cut}$ and $\omega_y = 0$ for $y > y_{cut}$.

We now relate the parameters $\sigma^2$ and $D$ with $\rho_n$ and $\beta_n$ in our model.  The parameter $\sigma^2$ represents the variance of the fitness distribution.  Recall from (\ref{Gaussian}) and (\ref{fnG}) that the long-run empirical distribution of particles in our model is approximately normal with variance $\rho_n/\beta_n$. In view of (\ref{betandef}), one unit of space in our model corresponds to $\beta_n$ units of fitness.
Thus the variance of the fitness distribution in our model is $\beta_n^2 (\rho_n/\beta_n) = \rho_n \beta_n$, leading to the correspondence
\[
\sigma^2=\rho_n\beta_n.
\]
Also, as explained in equation (1.20) in \cite{RS2020}, we expect $\beta_n^2$ to correspond to $\mu E[s^2]$, which implies that
\[
D=\frac{\beta_n^2}{2}.
\]
Plugging the above two formulas into (\ref{biosol}) and taking the scaling into consideration, we get
\[
\omega (\beta_n y)=Ce^{-\rho_n y}Ai\bigg(\frac{\rho^2}{(2\beta)^{2/3}}-(2\beta)^{1/3}y\bigg),
\]
matching (\ref{airydensity}).  Indeed, the derivation above which was adapted from \cite{Neher, Melissa} is quite similar to the derivation of (\ref{RSAiry}) in \cite{RS2020}. The derivation of (\ref{RSAiry}) leaned on work of Salminen \cite{Salminen}, which involved solving a differential equation for the density.

\subsection{Table of notation}
We summarize some of the notation that is used throughout the rest of the paper in the following table.
\clearpage
\begin{longtable}[!htbp]{p{.12\textwidth} | p{.82\textwidth} } 
\caption{\textit{Index of notation}}\\
\toprule
    $n$ & Index of a sequence of processes.  \\
    $\rho_n$   & Particles move according to Brownian motion with drift $-\rho_n$. \\
    $\beta_n$ & Selection parameter. The difference between the birth rate and the death rate for a particle at $x$ is $\beta_n x$. \\
    $N_{t,n}$  & Total number of particles at time $t$.\\
    $\mathcal{N}_{t,n}$ & The set of particles alive at time $t$.\\
     $N_{t,n}(\mathcal{I})$  & Number of particles in the interval $\mathcal{I}$ at time $t$.\\
    $X_{i,n}(t)$ & Positions of the particle $i$ at time $t$ for $i\in\mathcal{N}_{t,n}$.\\
    $\gamma_1$ & The largest zero of the Airy function.\\
    $L_n$ & The approximate position of right-most particle, $L_n=\rho_n^2/\beta_n-(2\beta_n)^{-1/3}\gamma_1$. \\
    $L_n^A$ & Defined to equal $L_n-A/\rho_n$ for $A\in\mathbb{R}$.\\
     $\bar{L}_n$ & The approximate position of left-most particle, $\bar{L}_n=-5\rho_n^2/8\beta_n-2(2\beta_n)^{-1/3}\gamma_1$. \\
    $L_n^*$ & The position that is near the position of the right-most particle. We call it the right edge. Explicitly, $L_n^*=\rho_n^2/\beta_n$.  \\
    $L_n^\dagger$ & The position that is near the position of the left-most particle. We call it the left edge. Explicitly, $L_n^\dagger=-5\rho_n^2/8\beta_n$.\\
    $Y_n(t)$ & Sum of $e^{\rho_nX_{i,n}(t)}$ for all $i=1,..., N_{t,n}$. Defined in (\ref{defY}). \\
    $Z_n(t)$ & Weighted sum used to characterize the size of the configuration at time $t$. Defined in (\ref{defZ}).\\
    $\lesssim$ & Write $x_n\lesssim y_n$ if $x_n/y_n$ is bounded above by a positive constant. Define $\gtrsim$ similarly.\\
    $\ll$ & Write $x_n\ll y_n$ if $\lim_{n\rightarrow\infty}x_n/y_n=0$. Define $\gg$ similarly.\\
    $\asymp$ & Write $x_n\asymp y_n$ if $x_n/y_n$ is bounded above  and below by positive constants.\\
    $O$ & Write $x_n=O(y_n)$ if the sequence $(x_n/y_n)_{n=1}^{\infty}$ is bounded.\\
    $o$  &Write $x_n=o(y_n)$ if $\lim_{n\rightarrow\infty}x_n/y_n=0$.\\
    $\{[a_n,b_n]\}_{n=1}^{\infty}$ & A sequence of intervals satisfying (\ref{b-a})-(\ref{bn}).\\
    $z_n$ & Roughly speaking, the asymptotic density of the number of particles in $[a_n,b_n]$ is maximized at $z_n$. Defined in (\ref{defzn}).\\
    $l_n$ &Measures the length of the interval in which we are counting the number of particles.\\
    $t_n(y)$ & For particles near $y$ at time $t_n$, $t_n-t_n(y)$ is the time when their ancestors start to leave the right boundary and drift toward $y$. Defined in (\ref{t_z}).\\
    $g_n(y)$ &Function used to approximate the density of particles. Defined in (\ref{g}).\\
    $M_{t,n}$ & Position of the left-most particle at time $t$.\\
    $p_t^n(x,y)$ &Density at location $y$ and time $t$ for the process which starts from a single particle at $x$ at time $0$. Defined in (\ref{density2.8}).\\
    $p_t^{L_n}(x,y)$ &Density at location $y$ and time $t$ for the process where there is only one particle at $x$ at time $0$ and particles are killed upon hitting $L_n$.\\
    $c_{0,n}$ & The ratio between $z_n$ and $L_n^*$. Defined in (\ref{defc0}).\\
    $c_n$ & Measures the distance between $z_n$ and $L_n^*$. Defined in (\ref{defcc0}).\\
    $r_{x}^{L_n}(v)$ &Rate at which particles hit $L_n$ at time $v$.  Defined in \ref{RS5}. \\
    $N_{t}^{L_n}(\mathcal{I})$ &Number of particles in the interval $\mathcal{I}$ at time $t$ for the process in which particles are killed at $L_n$.\\
    $(\mathcal{F}_t,t\geq 0)$ & Natural filtration associated with the branching Brownian motion process.\\
    $d$ & Used to divide the length of the interval in the proof of Proposition \ref{Global}. Defined in (\ref{dl}).\\
    $s$ & Constants used to adjust time. In the proof of Propositions \ref{Local} and \ref{Global}, we define $s=C_1\beta_n^{-2/3}$. \\ 
    $s_y$ & Constants used to adjust time for each $y\in [z-l,z+l]$ based on the choice of $s$. Defined to be $t(y)-t(z)+s$.\\
    $u_1$ & The first cutoff time in the second moment calculation. See Lemma \ref{Small2ndmoment}.\\
    $u_2$  & The second cutoff time in the second moment calculation. Defined in (\ref{u2}).\\
 \bottomrule
\end{longtable}

\subsection{Organization of the paper}

The rest of this paper is organized as follows.  In Section 2, we show how to obtain Theorem~\ref{ThmGlobal} and Corollary \ref{Tight} from two other propositions, one of which controls the number of particles in narrow intervals and one of which controls the number of particles in longer intervals. In Section~3, we prove Propositions \ref{Corollaryr} and \ref{Corollaryl}, and give the most important arguments for the proofs of the two propositions that lead to Theorem \ref{ThmGlobal}.  Proofs of some technical lemmas are postponed until Section \ref{Lemmasec}, and the second moment calculations are presented in Section \ref{2momsec}.

\section{Structure of the proof of the main results}\label{sec2}

In this section, we show how Theorem \ref{ThmGlobal} and Corollary \ref{Tight} follow from Propositions \ref{Local} and \ref{Global} below.  We also introduce some notation that will be used throughout the paper.

\subsection{Division into larger and smaller intervals}

The proof of Theorem \ref{ThmGlobal} will be divided into two cases. First, we will deal with intervals with smaller length. In such intervals, we will control the number of particles using a second moment argument. Indeed, we will show that most particles that end up near $z_n$ at time $t_n$ stay close to $L_n$ up to time $t_n-t_n(z_n)$ and then drift towards $z_n$. Trajectories of such particles are illustrated in Figure 1. Second, we will consider longer intervals.  We will show that the number of particles in the interval $[a_n, b_n]$ that are far away from $z_n$ is negligible using a first moment argument, allowing us to estimate the number of particles in the entire interval by the number of particles in a smaller interval around $z_n$.  The first step will lead to Proposition \ref{Local} while the second step will lead to Proposition \ref{Global}. 

Consider a sequence $(z_n)_{n=1}^{\infty}$ satisfying (\ref{zn}) such that 
\begin{equation}\label{zn1}
|z_n|\gtrsim\sqrt{\frac{\rho_n}{\beta_n}}\qquad \text{or}\qquad |z_n|\ll\sqrt{\frac{\rho_n}{\beta_n}}.
\end{equation}
We further assume that
\begin{equation}\label{zn2}
z_n\geq 0\quad\text{for all}\; n\qquad \text{or}\qquad z_n\leq 0\quad\text{for all}\;n.
\end{equation}
Denote
\begin{equation}\label{defc0}
c_{0,n}=\frac{z_n}{L_n^*}.
\end{equation}
We consider intervals of the forms $[z_n,z_n+l_n]$ and $[z_n-l_n,z_n]$ where $l_n$ is the length of the interval. By convention, if $l_n=\infty$, then $[z_n,z_n+l_n]=[z_n,\infty)$ and $[z_n-l_n,z_n]=(-\infty,z_n]$.

\begin{Prop}\label{Local}
Suppose assumptions (\ref{assump1}), (\ref{assump2}), (\ref{assump3}), (\ref{assumpY}) and (\ref{assumpZ}) hold. For every sequence $(z_n)_{n=1}^{\infty}$ satisfying (\ref{zn}), (\ref{zn1}) and (\ref{zn2}), choose $(l_n)_{n=1}^{\infty}$ such that
\begin{equation}\label{ln}
 \begin{cases} 
      1\ll l_n \lesssim \frac{1}{|c_{0,n}|\rho_n} &\qquad \text{if}\;|z_n|\gtrsim \sqrt{\frac{\rho_n}{\beta_n}},\\
      1\ll l_n \lesssim \sqrt{\frac{\rho_n}{\beta_n}}  &\qquad \text{if}\;|z_n|\ll \sqrt{\frac{\rho_n}{\beta_n}}.
   \end{cases}
\end{equation}
Consider intervals of the form
\begin{equation}\label{I}
 \mathcal{I}_n=\begin{cases} 
     [z_n,z_n+l_n] &\qquad \text{if}\;z_n\geq 0,\\
     [z_n-l_n,z_n]  &\qquad \text{if}\;z_n\leq 0.
   \end{cases}
\end{equation}
If $t_n$ satisfies 
\begin{equation}\label{assumptn}
\frac{\rho_n^{2/3}}{\beta_n^{8/9}}\ll t_n-t_n(z_n)\ll \frac{\rho_n}{\beta_n},
\end{equation}
then for any $\kappa>0$,  we have
\begin{align}\label{mainpropl}
\lim_{n\rightarrow\infty}P\bigg(\frac{1-\kappa}{Ai'(\gamma_1)^2}e^{-\rho_n L_n^*}Z_n(0)&\int_{\mathcal{I}_n}\frac{1}{\sqrt{2\pi t_n(y)}}e^{g_n(y)}dy\leq N_{t_n,n}(\mathcal{I}_n)
\nonumber\\
&\leq\frac{1+\kappa}{Ai'(\gamma_1)^2}e^{-\rho_n L_n^*}Z_n(0)\int_{\mathcal{I}_n}\frac{1}{\sqrt{2\pi t_n(y)}}e^{g_n(y)}dy\bigg)=1.
\end{align}
\end{Prop}

\begin{Prop}\label{Global}
Suppose assumptions (\ref{assump1}), (\ref{assump2}), (\ref{assump3}), (\ref{assumpY}) and (\ref{assumpZ}) hold. For every sequence of ~$(z_n)_{n=1}^{\infty}$ satisfying (\ref{zn}), (\ref{zn1}) and (\ref{zn2}), choose $(l_n)_{n=1}^{\infty}$ such that
\begin{equation}\label{lng}
 \begin{cases} 
      l_n \gg \frac{1}{|c_{0,n}|\rho_n} &\qquad \text{if}\;\;|z_n|\gtrsim \sqrt{\frac{\rho_n}{\beta_n}},\\
      l_n \gg \sqrt{\frac{\rho_n}{\beta_n}}  &\qquad \text{if}\;\;|z_n|\ll \sqrt{\frac{\rho_n}{\beta_n}}.
   \end{cases}
\end{equation}
Consider intervals of the form
\begin{equation}\label{J}
 \mathcal{J}_n=\begin{cases} 
     [z_n,z_n+l_n] &\qquad \text{if}\;\;z_n\geq 0,\\
     [z_n-l_n,z_n]  &\qquad \text{if}\;\;z_n\leq 0.
   \end{cases}
\end{equation}
If $t_n$ satisfies (\ref{assumptn}), then for any $\kappa>0$,  we have
\begin{align}\label{mainpropg}
&\lim_{n\rightarrow\infty}P\bigg(\frac{1-\kappa}{Ai'(\gamma_1)^2}e^{-\rho_n L_n^*}Z_n(0)\int_{\mathcal{J}_n\cap (-\infty,L_n^*]}\frac{1}{\sqrt{2\pi t_n(y)}}e^{g_n(y)}dy\leq N_{t_n,n}(\mathcal{J}_n)
\nonumber\\
&\hspace{1.5in}\leq\frac{1+\kappa}{Ai'(\gamma_1)^2}e^{-\rho_n L_n^*}Z_n(0)\int_{\mathcal{J}_n\cap (-\infty,L_n^*]}\frac{1}{\sqrt{2\pi t_n(y)}}e^{g_n(y)}dy\bigg)=1.
\end{align}
\end{Prop}

Next, we will explain heuristically why the interval length $l_n$ is divided into the above two cases (\ref{ln}) and (\ref{lng}). Let us take the case $z_n\geq 0$ as an example. The case  when $z_n\leq 0$ is essentially the same. Our hope is to find a cutoff length $l_n$ depending on $z_n$ such that the number of particles in $[z_n,\infty)$ is dominated by the number of particles in $[z_n,z_n+l_n]$.  Since the number of particles near $z_n$ is approximately proportional to $e^{g_n(z_n)}/\sqrt{2\pi t_n(z_n)}$, this boils down to finding a cutoff length $l_n$ such that for any $\eta>0$, if $n$ is sufficiently large, then
\[
\int_{z_n}^{\infty}\frac{1}{\sqrt{2\pi t_n(y)}}e^{g_n(y)}dy<(1+\eta)\int_{z_n}^{z_n+l_n}\frac{1}{\sqrt{2\pi t_n(y)}}e^{g_n(y)}dy.
\] 
It turns out that if $z_n\gtrsim\sqrt{\rho_n/\beta_n}$, then we can take $l_n\asymp1/c_{0,n}\rho_n$, as shown in Lemma \ref{Integralofg}, while if $z_n\ll \sqrt{\rho_n/\beta_n}$, then we take $l_n\asymp \sqrt{\rho_n/\beta_n}$.

\subsection{Proof of Theorem \ref{ThmGlobal}}

In this subsection, we deduce Theorem \ref{ThmGlobal} from Propositions \ref{Local} and \ref{Global}.  We first review an important result from \cite{RS2020} which will be needed in the proof.

\begin{Rmk}\label{RS1}
Proposition 2.3 in \cite{RS2020} states that if $\beta_nt_n/\rho_n$ converges to a positive real number as $n$ goes to infinity, then with probability tending to $1$ as $n\rightarrow\infty$, conditions (\ref{assumpY}) and (\ref{assumpZ}) hold with $Y_n(t_n)$ and $Z_n(t_n)$ in place of $Z_n(0)$ and $Y_n(0)$ respectively. Furthermore, if $t_n \asymp \rho_n/\beta_n$, then for every subsequence $(n_j)_{j=1}^{\infty}$, there exists a sub-subsequence $(n_{j_k})_{k=1}^{\infty}$ such that 
\[
\lim_{k\rightarrow\infty}\frac{\beta_{n_{j_k}}t_{n_{j_k}}}{\rho_{n_{j_k}}}=\tau\in (0,\infty).
\]
Consequently, by Proposition 2.3 in \cite{RS2020}, with probability tending to $1$ as $k\rightarrow\infty$, conditions (\ref{assumpY}) and (\ref{assumpZ}) hold with $Y_{n_{j_k}}(t_{n_{j_k}})$ and $Z_{n_{j_k}}(t_{n_{j_k}})$ in place of $Z_n(0)$ and $Y_n(0)$ respectively. 
\end{Rmk}

\bigskip
\noindent\textit{Proof of Theorem \ref{ThmGlobal}.} First, we consider the case when $t_n$ satisfies (\ref{thmassump1}). To prove (\ref{thm1}), it suffices to show that for every subsequence $(n_j)_{j=1}^{\infty}$, there exists a sub-subsequence $(n_{j_k})_{k=1}^{\infty}$, such that for any $0<\kappa<1$,
\begin{align}\label{goal1}
&\lim_{k\rightarrow\infty}P\bigg(\frac{1-\kappa}{Ai'(\gamma_1)^2}e^{-\rho_{n_{j_k}} L^*_{n_{j_k}}}Z_{n_{j_k}}(0)\int_{[a_{n_{j_k}},b_{n_{j_k}}]\cap (-\infty,L^*_{n_{j_k}}]}\frac{1}{\sqrt{2\pi t_{n_{j_k}}(y)}}e^{g_{n_{j_k}}(y)}dy\nonumber\\
&\hspace{0.4in}\leq N_{t_{n_{j_k}},n_{j_k}}\big([a_{n_{j_k}},b_{n_{j_k}}]\big)\nonumber\\
&\hspace{0.6in}\leq\frac{1+\kappa}{Ai'(\gamma_1)^2}e^{-\rho_{n_{j_k}} L_{n_{j_k}}^*}Z_{n_{j_k}}(0)\int_{[a_{n_{j_k}},b_{n_{j_k}}]\cap (-\infty,L_{n_{j_k}}^*]}\frac{1}{\sqrt{2\pi t_{n_{j_k}}(y)}}e^{g_{n_{j_k}}(y)}dy\bigg)=1.
\end{align}
Given a subsequence $(n_j)_{j=1}^{\infty}$, there exists a further subsequence $(n_{j_k})_{k=1}^{\infty}$ such that one of the following holds:
\begin{enumerate}
\item We have $a_{n_{j_k}}\geq 0$ for all $k$.  Let $z_{n_{j_k}}= a_{n_{j_k}}$ and $l_{n_{j_k}}=b_{n_{j_k}}-a_{n_{j_k}}$. The subsequence $(z_{n_{j_k}})_{k=1}^{\infty}$ satisfies (\ref{zn}), (\ref{zn1}) and (\ref{zn2}), and the subsequence $(l_{n_{j_k}})_{k=1}^{\infty}$ satisfies (\ref{ln}) or (\ref{lng}).
\item We have $b_{n_{j_k}}\leq 0$ for all $k$. Let $z_{n_{j_k}}= b_{n_{j_k}}$ and $l_{n_{j_k}}=b_{n_{j_k}}-a_{n_{j_k}}$. The subsequence $(z_{n_{j_k}})_{k=1}^{\infty}$ satisfies (\ref{zn}), (\ref{zn1}) and (\ref{zn2}), and the subsequence $(l_{n_{j_k}})_{k=1}^{\infty}$ satisfies (\ref{ln}) or (\ref{lng}).
\item We have $a_{n_{j_k}}< 0$ and $b_{n_{j_k}}>0$ for all $k$. Let $z_{n_{j_k}}= 0$, $l_{1,n_{j_k}}=-a_{n_{j_k}}$ and $l_{2,n_{j_k}}=b_{n_{j_k}}$. Both the subsequences $(l_{1,n_{j_k}})_{k=1}^{\infty}$ and $(l_{2,n_{j_k}})_{k=1}^{\infty}$ satisfy (\ref{ln}) or (\ref{lng}).
\end{enumerate}
In cases 1 and 2, since $[a_{n_{j_k}},b_{n_{j_k}}]$ satisfies the hypotheses of either Proposition \ref{Local} or Proposition~\ref{Global}, equation (\ref{goal1}) follows from (\ref{mainpropl}) or (\ref{mainpropg}). As for case 3,  we see that both $[a_{n_{j_k}},0]$ and $[0,b_{n_{j_k}}]$ satisfy the hypotheses of Proposition \ref{Local} or Proposition \ref{Global}. Thus, both $[a_{n_{j_k}},0]$ and $[0,b_{n_{j_k}}]$ satisfy (\ref{goal1}) with $[a_{n_{j_k}},0]$ and $[0,b_{n_{j_k}}]$ in place of $[a_{n_{j_k}},b_{n_{j_k}}]$ respectively. Consequently, equation (\ref{goal1}) also holds in this case. Therefore, equation (\ref{thm1}) follows.

Next, consider the case when $t_n$ satisfies (\ref{thmassump2}). Choose a sequence $(h_{n})_{n=1}^{\infty}$ for which 
\begin{equation}\label{hnk}
\frac{\rho_{n}^{2/3}}{\beta_{n}^{8/9}}\ll h_{n}\ll \frac{\rho_{n}}{\beta_{n}}.
\end{equation}
Let 
\begin{equation}\label{rn}
r_{n}=t_{n}-t_{n}(z_{n})-h_{n}.
\end{equation} Note that $r_{n}\asymp\rho_n/\beta_n$. By Remark \ref{RS1}, for every subsequence $(n_j)_{j=1}^{\infty}$, we can choose a sub-subsequence $(n_{j_k})_{k=1}^{\infty}$ such that assumptions (\ref{assumpY}) and (\ref{assumpZ}) hold when $Y_{n}(0)$ and $Z_{n}(0)$ are replaced by $Y_{n_{j_k}}(r_{n_{j_k}})$ and $Z_{n_{j_k}}(r_{n_{j_k}})$. By using the Markov property at time $r_{n_{j_k}}$ and applying the previous argument, there exists a further sub-subsequence $(n_{j_{k_m}})_{m=1}^{\infty}$ such that equation (\ref{goal1}) holds with $Z_{n_{j_{k_m}}}(r_{n_{j_{k_m}}})$ in place of $Z(0)$. As a result, we have for any $0<\kappa<1$,
\begin{align}\label{goal2}
\lim_{n\rightarrow\infty}P\bigg(\frac{1-\kappa}{Ai'(\gamma_1)^2}&e^{-\rho_{n} L^*_{n}}Z_{n}(r_{n})\int_{[a_{n},b_{n}]\cap (-\infty,L^*_{n}]}\frac{1}{\sqrt{2\pi t_{n}(y)}}e^{g_{n}(y)}dy\leq N_{t_{n},n}\big([a_{n},b_{n}]\big)\nonumber\\
&\leq\frac{1+\kappa}{Ai'(\gamma_1)^2}e^{-\rho_{n} L_{n}^*}Z_{n}(r_{n})\int_{[a_{n},b_{n}]\cap (-\infty,L_{n}^*]}\frac{1}{\sqrt{2\pi t_{n}(y)}}e^{g_{n}(y)}dy\bigg)=1.
\end{align}
Note that equation (\ref{goal2}) holds for all choices of $(t_n)_{n=1}^{\infty}$ satisfying (\ref{thmassump2}) and $(h_{n})_{n=1}^{\infty}$ satisfying (\ref{hnk}). Thus for every $(z_{n})_{n=1}^{\infty}$ satisfying (\ref{zn}), $(t_n)_{n=1}^{\infty}$ satisfying (\ref{thmassump2}) and any two sequences $(h_{1,n})_{n=1}^{\infty}$, $(h_{2,n})_{n=1}^{\infty}$ satisfying (\ref{hnk}), we have
\begin{equation}\label{Z}
\lim_{n\rightarrow\infty}P\bigg(\frac{1-\kappa}{1+\kappa}Z_{n}(r_{1,n})\leq Z_{n}(r_{2,n})\leq \frac{1+\kappa}{1-\kappa}Z_{n}(r_{1,n})\bigg)=1,
\end{equation}
where $r_{i,n}=t_{n}-t_{n}(z_{n})-h_{i,n}$ for $i=1,2$. Choose $(z^*_{n})_{n=1}^{\infty}$ satisfying (\ref{zn}), $(t^*_n)_{n=1}^{\infty}$ satisfying (\ref{thmassump2}) and $(h_{1,n})_{n=1}^{\infty}$, $(h_{2,n})_{n=1}^{\infty}$ satisfying (\ref{hnk}) such that
\[
t_n-t_n(z_n)=t_n^*-t_n(z_n^*)-h_{1,n},\qquad t_n-t_n(z_n)-h_n=t_n^*-t_n(z_n^*)-h_{2,n}.
\]
For example, for any sequence of $(h_{1,n})_{n=1}^{\infty}$ satisfying (\ref{hnk}), we can take $z_n^*=z_n$, $t_n^*=t_n+h_{1,n}$ and $h_{2,n}=h_{1,n}+h_n$. By (\ref{rn}) and (\ref{Z}), we have
\begin{equation}\label{Z(r)Z(t-t(z))}
\lim_{n\rightarrow\infty}P\bigg(\frac{1-\kappa}{1+\kappa}Z_{n}\big(t_{n}-t_n(z_n)\big)\leq Z_{n}(r_{n})\leq \frac{1+\kappa}{1-\kappa}Z_{n}\big(t_{n}-t_n(z_n)\big)\bigg)=1.
\end{equation}
Finally, equation (\ref{thm2}) follows from (\ref{goal2}) and (\ref{Z(r)Z(t-t(z))}).
 \qedwhite
\\
\begin{Rmk}
The argument leading to (\ref{Z(r)Z(t-t(z))}) can be modified to show that for $t_n\asymp \rho_n/\beta_n$ and $h_n\ll\rho_n/\beta_n$, as $n\rightarrow\infty$,
\begin{equation}\label{Zratio}
\frac{Z_n(t_n)}{Z_n(t_n+h_n)}\rightarrow_{p}1.
\end{equation}
To see this, note that we can choose $(z^*_{n})_{n=1}^{\infty}$ satisfying (\ref{zn}), $(t^*_n)_{n=1}^{\infty}$ satisfying (\ref{thmassump2}), and $(h_{1,n})_{n=1}^{\infty}$ and $(h_{2,n})_{n=1}^{\infty}$ satisfying (\ref{hnk}) such that
\[
t_n=t_n^*-t_n(z_n^*)-h_{1,n},\qquad t_n+h_n=t_n^*-t_n(z_n^*)-h_{2,n}.
\]
For example, we can take $(h_{2,n})_{n=1}^{\infty}$ to be any sequence satisfying (\ref{hnk}) and $(z_n^*)_{n=1}^{\infty}$ to be any sequence satisfying (\ref{zn}) such that $t_n(z_n^*) \ll \rho_n/\beta_n$.  Then we let $h_{1,n} = h_{2,n} + h_n$ and $t_n^* = t_n + t_n(z_n^*) + h_{1,n}$.
Letting $r_{i,n}=t_{n}^*-t_{n}(z_{n}^*)-h_{i,n}$ for $i=1,2$, equation (\ref{Zratio}) follows from (\ref{Z}).
\end{Rmk}

As a byproduct of the proof of Theorem \ref{ThmGlobal}, the following lemma shows that the number of particles in any given interval will not change much on a time scale shorter than $\rho_n/\beta_n$.

\begin{Lemma}\label{Difference}
Suppose (\ref{assump1}), (\ref{assump2}), (\ref{assump3}), (\ref{assumpY}) and (\ref{assumpZ}) hold. For every sequence $\{[a_n,b_n]\}_{n=1}^{\infty}$ satisfying (\ref{b-a})-(\ref{bn}), define $z_n$ according to (\ref{defzn}). Suppose
\begin{equation}\label{diffassump1}
\frac{\rho_n^{2/3}}{\beta_n^{8/9}}\ll t_n-t_{n}(z_n)\lesssim \frac{\rho_n}{\beta_n}, \qquad\frac{\rho_n^{2/3}}{\beta_n^{8/9}}\ll t'_n-t_{n}(z_n)\lesssim \frac{\rho_n}{\beta_n}.
\end{equation}
If 
\begin{equation}\label{diffasump2}
|t_n-t_n'|\ll\frac{\rho_n}{\beta_n},
\end{equation}
then as $n\rightarrow\infty$,
\begin{equation}\label{diffconv}
\frac{N_{t_n,n}([a_n,b_n])}{N_{t'_n,n}([a_n,b_n])}\rightarrow_{p}1.
\end{equation}
\end{Lemma}
\noindent\textit{Proof of Lemma \ref{Difference}.}
First, we consider the case $\rho_n^{2/3}/\beta_n^{8/9}\ll t_n-t_n(z_n)\ll\rho_n/\beta_n$. From (\ref{diffasump2}), we also have $\rho_n^{2/3}/\beta_n^{8/9}\ll t'_n-t_n(z_n)\ll\rho_n/\beta_n$. Therefore, both $N_{t_n,n}([a_n,b_n])$ and $N_{t'_n,n}([a_n,b_n])$ satisfy (\ref{thm1}) and equation (\ref{diffconv}) follows.

It remains to consider the case $t_n-t_n(z_n)\asymp\rho_n/\beta_n$. By (\ref{diffasump2}), we also have $t'_n-t_n(z_n)\asymp\rho_n/\beta_n$. Then both $N_{t_n,n}([a_n,b_n])$ and $N_{t'_n,n}([a_n,b_n])$ satisfy (\ref{thm2}). As a result, equation (\ref{diffconv}) follows from (\ref{thm2}) and (\ref{Zratio}).
\qedwhite

\subsection{Proof of Corollary \ref{Tight}}

In this subsection, we show how to obtain Corollary \ref{Tight} from Theorem \ref{ThmGlobal}.  We first recall that Proposition 2.2 in \cite{RS2020} proves that if 
\begin{equation}\label{prop2.2assump}
\frac{\rho_n^{2/3}}{\beta_n^{8/9}}\ll t_n-t_n(0)\ll \frac{\rho_n}{\beta_n},
\end{equation}
then for all $\kappa>0$, 
\begin{equation}\label{prop2.2}
\lim_{n\rightarrow\infty}P\bigg(\frac{1-\kappa}{Ai'(\gamma_1)^2}e^{-\rho_n^3/3\beta_n}Z_n(0)\leq N_{t_n,n}\leq \frac{1+\kappa}{Ai'(\gamma_1)^2}e^{-\rho_n^3/3\beta_n}Z_n(0)\bigg)=1.
\end{equation}
Equation (\ref{prop2.2}) also follows from (\ref{thm1}) with $a_n = -\infty$ and $b_n = \infty$, 
once it is established that for $n$ sufficiently large, we have
\begin{equation}\label{prop2.2claim}
1-\eta<e^{-\rho_n^3/6\beta_n}\int_{-\infty}^{L_n^*}\frac{1}{\sqrt{2\pi t_n(y)}}e^{g_n(y)}dy<1+\eta.
\end{equation}
One can obtain (\ref{prop2.2claim}) by comparing the density $f_n$ defined in the statement of Corollary \ref{Tight} to the Gaussian density $f_n^G$ defined in (\ref{fnG}).  We omit the details.

\bigskip
\noindent\textit{Proof of Corollary \ref{Tight}.}
By applying the same argument leading to (\ref{thm2}) in the proof of Theorem \ref{ThmGlobal}, equation (\ref{prop2.2}) implies that if $t_n-t_n(0)\asymp\rho_n/\beta_n$, then
\begin{equation}\label{prop2.2+}
\lim_{n\rightarrow\infty}P\bigg(\frac{1-\kappa}{Ai'(\gamma_1)^2}e^{-\rho_n^3/3\beta_n}Z_n\big(t_n-t_n(0)\big)\leq N_{t_n,n}\leq \frac{1+\kappa}{Ai'(\gamma_1)^2}e^{-\rho_n^3/3\beta_n}Z_n\big(t_n-t_n(0)\big)\bigg)=1.
\end{equation}
If 
\[
\frac{\rho_n^{2/3}}{\beta_n^{8/9}}\ll t_n-t_n(z_n)\ll\frac{\rho_n}{\beta_n}, \qquad \frac{\rho_n^{2/3}}{\beta_n^{8/9}}\ll t_n-t_n(0)\ll\frac{\rho_n}{\beta_n},
\]
then by (\ref{thm1}) and (\ref{prop2.2}), we have for any $\kappa>0$
\begin{equation}\label{tightcase1}
\lim_{n\rightarrow\infty}P\bigg(\frac{1-\kappa}{1+\kappa}<D_n<\frac{1+\kappa}{1-\kappa}\bigg)=1.
\end{equation}
Thus, the sequence $(D_n)_{n=1}^{\infty}$ is tight. If 
\[
t_n-t_n(z_n)\asymp\frac{\rho_n}{\beta_n}, \qquad \frac{\rho_n^{2/3}}{\beta_n^{8/9}}\ll t_n-t_n(0)\ll\frac{\rho_n}{\beta_n},
\]
then by (\ref{thm2}) and (\ref{prop2.2}), we have for any $\kappa>0$, 
\begin{equation}\label{tightcase2}
\lim_{n\rightarrow\infty}P\bigg(\frac{1-\kappa}{1+\kappa}\frac{Z_n(t_n-t_n(z_n))}{Z_n(0)}<D_n<\frac{1+\kappa}{1-\kappa}\frac{Z_n(t_n-t_n(z_n))}{Z_n(0)}\bigg)=1.
\end{equation}
Note that $Z_n(0)$ satisfies (\ref{assumpZ}) and $Z_n(t_n-t_n(z_n))$ also satisfies (\ref{assumpZ}) by Remark \ref{RS1}. Equation (\ref{tightcase2}) thus implies that $(D_n)_{n=1}^{\infty}$ is tight. The remaining two cases 
\[
\frac{\rho_n^{2/3}}{\beta_n^{8/9}}\ll  t_n-t_n(z_n)\ll\frac{\rho_n}{\beta_n}, \qquad t_n-t_n(0)\asymp\frac{\rho_n}{\beta_n},
\]
and
\[
t_n-t_n(z_n)\asymp\frac{\rho_n}{\beta_n}, \qquad t_n-t_n(0)\asymp\frac{\rho_n}{\beta_n},
\] 
follow by essentially the same argument as the second case, using (\ref{prop2.2+}) in place of (\ref{prop2.2}). If $0\in [a_n,b_n]$, then (\ref{tightcase1}) holds true and thus $D_n\rightarrow1$ in probability as $n\rightarrow\infty$.
\qedwhite

\section{Proof of Propositions \ref{Corollaryr}, \ref{Corollaryl}, \ref{Local}, and \ref{Global}}\label{sec3}

In this section, we give the main arguments in the proofs of Propositions \ref{Corollaryr}, \ref{Corollaryl}, \ref{Local}, and \ref{Global}.  We defer the proofs of several technical lemmas until Section \ref{Lemmasec} and the proof of the second moment estimates until Section \ref{2momsec}.

\subsection{A review of results from \cite{RS2020}}

In this subsection, we will collect some of the results in \cite{RS2020} that will be used in the proofs. Suppose (\ref{assumpY}) and (\ref{assumpZ}) hold.

\begin{enumerate}[label=\textbf{3.1.\arabic*},ref=3.1.\arabic*]
\item\label{RS2} Lemma 6.1 in \cite{RS2020} shows that 
\begin{equation}\label{lem6.1}
\lim_{r\rightarrow\infty}e^{-r^3/3}\int_{0}^{\infty}e^{r(\gamma_1+z)}Ai(\gamma_1+z)dz=1.
\end{equation}
Based on equations (6.5), (6.6) and Lemma 6.1 in \cite{RS2020}, for any $\eta>0$, there exists a constant $C_5\geq 2$  sufficiently large such that
\begin{equation}\label{C1cond1}
\bigg(1-\frac{\eta}{2}\bigg)e^{C_5^3/6}\leq \int_{0}^{\infty}e^{2^{-1/3}C_5(\gamma_1+y)}Ai(\gamma_1+y)dy
\leq (1+\eta)e^{C_5^3/6},
\end{equation}
and
\begin{equation}\label{C1cond2}
e^{-C_5^3/6}\frac{(1+\eta)\sqrt{2\pi}}{Ai'(\gamma_1)^2}\int_{0}^{\infty}Ai(\gamma_1+y)dy<\eta
\end{equation}
hold. Furthermore, there exists a constant $C_6\geq -2^{-1/3}\gamma_1$ sufficiently large such that
\begin{equation}\label{C2}
\int_{2^{1/3}C_6}^{\infty}e^{2^{-1/3}C_5(\gamma_1+y)}Ai(\gamma_1+y)dy<\frac{\eta}{2}e^{C_5^3/48}.
\end{equation}

\item\label{RS3} Fix $A\in\mathbb{R}$. Define
\begin{equation}\label{LnA}
L_n^A = L_n-\frac{A}{\rho_n}.
\end{equation}
Lemma 5.1 in \cite{RS2020} proves that there exists a constant $C_7$ such that the probability that some particle that is to the right of $L_n^A$ at time $0$ has a descendant alive in the population at time $C_7\rho_n^{-2}$ tends to $0$ as $n\rightarrow\infty$. Moreover, according to the argument leading to (5.9) in \cite{RS2020}, for $t_n \ll \rho_n/\beta_n$, it follows that with probability tending to $1$ as $n\rightarrow\infty$, no particle that hits $L_n^A$ before time $t_n-C_7\rho_n^{-2}$ has descendants alive at time $t_n$.

\item \label{RS4} Consider the process in which particles are killed upon hitting $L_n$. If this process starts from a single particle at $x$, we denote the density of this process at time $t$ by $p_{t}^{L_n}(x,y)$. Lemma~2.5 in \cite{RS2020} implies that if $x,y<L_n$ and 
\begin{equation}\label{2.24}
(2\beta_n)^{1/6}\big((L_n-x)^{1/2}+(L_n-y)^{1/2}\big)-2^{-1/3}\beta_n^{2/3}t_n \rightarrow-\infty,
\end{equation}
then there exits a constant $C_8$ such that
\begin{equation}\label{lem2.5}
p_{t_n}^{L_n}(x,y)\leq C_4\beta_n^{1/3}e^{\rho_n x}Ai((2\beta_n)^{1/3}(L_n-x)+\gamma_1)e^{-\rho_n y}Ai((2\beta_n)^{1/3}(L_n-y)+\gamma_1).
\end{equation}
Define 
\begin{equation}\label{H(u)}
H_n(t)=L_n-\frac{\beta_n t^2}{9}. 
\end{equation}
Equation (5.5) in \cite{RS2020} states that if $x\leq H_n(t_n)$, $0\leq \zeta_n \leq \beta_n t_n/2$ and $\beta_n^{-2/3}\ll t_n \ll \rho_n/\beta_n$, then 
\begin{equation}\label{(5.5)}
\int_{-\infty}^{L_n}p_{t_n}^{L_n}(x,y)e^{(\rho_n-\zeta_n)y}dy\ll e^{\rho_n x}e^{-\beta_n^2t_n^3/73}.
\end{equation}
Equation (5.6) in \cite{RS2020} states that if $x< L_n$, $0\leq \zeta_n \leq \beta_n t_n/2$ and $\beta_n^{-2/3}\ll t_n\ll \rho_n/\beta_n$, then 
\begin{equation}\label{(5.6)}
\int_{-\infty}^{H_n(t_n)}p_{t_n}^{L_n}(x,y)e^{(\rho_n-\zeta_n)y}dy\ll e^{\rho_n x}e^{-\beta_n^2 t_n^3/73}.
\end{equation}

\item \label{RS5} Consider the process in which particles are killed upon hitting $K$. If the process starts with a single particle at $x<K$, we denote by $r_{x,n}^{K}(v)$ the rate at which particles hit $K$ at time $v$ and $r_{x,n}^{K}(u,t)$ the expected number of particles that are killed at $K$ between times $u$ and $t$. Then
\[
r_{x,n}^{K}(u,t)=\int_{u}^{t}r_{x,n}^{K}(v)dv.
\] 
If $K=L_n$, then by (6.29) in \cite{RS2020}, there exists a constant $C_9$ such that
\begin{equation}\label{rxLv}
r_{x,n}^{L_n}(v)\leq \frac{C_9(L_n-x)}{v^{3/2}}\exp\bigg(\rho_n x-\rho_n L_n-\frac{(L_n-x)^2}{2v}-2^{-1/3}\beta_n^{2/3}\gamma_1v\bigg).
\end{equation}
Furthermore, for $K=L_n^A$, define $A_{-}=\max\{-A,0\}$. By Lemma 2.13 in \cite{RS2020}, we have for all $x<L_n$ and $0\leq u<t$,
\begin{equation}\label{lem2.13}
r^{L_n}_{x,n}(u,t)\lesssim e^{\rho_nx}e^{-\rho_nL_n^A}e^{-\beta_n^2u^3/9}+\beta_n^{2/3}(t-u)e^{-\rho_nL_n^A}e^{\beta_nA_{-}t/\rho_n}e^{\rho_nx}Ai\big((2\beta_n)^{1/3}(L_n-x)+\gamma_1\big).
\end{equation}

\item \label{RS6} Suppose 
\begin{equation}\label{assump5.8}
\beta^{-2/3}\log^{1/3}\bigg(\frac{\rho}{\beta^{1/3}}\bigg)\ll t_n\ll\frac{\rho_n}{\beta_n}.
\end{equation}
Let $f:\mathbb{R}\rightarrow[0,\infty)$ be a bounded measurable function. Define
\begin{equation}\label{Phi}
\Phi_n(f)=\sum_{i\in\mathcal{N}_{t_n,n}}e^{\rho X_{i,n}(t_n)}f\big((2\beta_n)^{1/3}(L_n-X_{i,n}(t_n))\big).
\end{equation}
According to (5.8) in \cite{RS2020}, we have for any $\kappa>0$,
\begin{align}\label{5.8}
\lim_{n\rightarrow\infty}P\bigg(\frac{1-\kappa}{Ai'(\gamma_1)^2}\bigg(\int_{0}^{\infty}&f(z)Ai(\gamma_1+z)\bigg)Z_n(0)<\Phi_n(f)\nonumber\\
&<\frac{1+\kappa}{Ai'(\gamma_1)^2}\bigg(\int_{0}^{\infty}f(z)Ai(\gamma_1+z)\bigg)Z_n(0)\bigg)=1.
\end{align}

\item \label{RS7} For $\beta_n^{-1/3}\ll x_n\ll\beta_n^{-1}$, consider the process started from a single particle at $x_n$. According to Lemma 2.14 in \cite{RS2020}, there is a positive constant $C_{10}$ such that for large enough $n$, the probability that the process survives until time $C_{10}/(\beta_nx_n)$ is bounded above by $2\beta_nx_n/\alpha$. Here, $\alpha$ is the constant appearing in assumption (\ref{assump3}).
\end{enumerate}

\subsection{Notation}
Here we introduce one more piece of notation which will be used throughout the rest of this paper. Recall that $c_{0,n}=z_n/L_n^*$. We denote 
\begin{equation}\label{defcc0}
c_n=\sqrt{1-c_{0,n}}. 
\end{equation}
We can now write 
\begin{equation}\label{tL-z}
t_n(z_n)=\frac{c_n\rho_n}{\beta_n}, \qquad L_n^*-z_n=\frac{c_n^2\rho_n^2}{2\beta_n}.
\end{equation}
The notation $c_{0,n}$ and $c_n$ will be useful in simplifying expressions involving $z_n$, $L_n^*-z_n$ and $t_n(z_n)$. Therefore, we list some of the most useful formulas involving $c_{0,n}$ and $c_n$ below. We see that for $z_n\in (L^\dagger_n, L^*_n)$,
\begin{equation}\label{range}
-\frac{5}{4}<c_{0,n}<1,\qquad 0<c_n<\frac{3}{2}.
\end{equation}
We also have the following equivalent asymptotic expressions: 
\begin{equation}\label{1-c=c_0}
|1-c_n|=\frac{|c_{0,n}|}{1+c_n}\asymp |c_{0,n}|,
\end{equation}
\begin{equation}\label{c0lb}
|z_n|\gtrsim \sqrt{\frac{\rho_n}{\beta_n}} \Longleftrightarrow |c_{0,n}|\gtrsim \frac{\beta_n^{1/2}}{\rho_n^{3/2}},
\end{equation}
\begin{equation}\label{clb}
L_n^*-z_n\gg \frac{1}{\beta_n^{1/3}}  \Longleftrightarrow c_n\gg \frac{\beta_n^{1/3}}{\rho_n},
\end{equation}
\begin{equation*}
z_n-L_n^\dagger\gg \frac{1}{\beta_n^{1/3}} \Longleftrightarrow c_{0,n}+\frac{5}{4}\gg \frac{\beta_n^{2/3}}{\rho_n^2} \Longleftrightarrow \frac{9}{4}-c_n^2\gg \frac{\beta_n^{2/3}}{\rho_n^2} \Longleftrightarrow \frac{3}{2}-c_n\gg \frac{\beta_n^{2/3}}{\rho_n^2}.
\end{equation*}
Moreover, if $|z_n|\gtrsim \sqrt{\rho_n/\beta_n}$ and $z_n$ satisfies assumption (\ref{zn}), then 
\begin{equation}\label{cc0}
|c_nc_{0,n}|\gtrsim \frac{\beta_n^{1/2}}{\rho_n^{3/2}}.
\end{equation}
Also if $|z_n|\ll\sqrt{\rho_n/\beta_n}$, then
\begin{equation}\label{cc02}
c_n\asymp1,\qquad |c_{0,n}|\ll\frac{\beta_n^{1/2}}{\rho_n^{3/2}}.
\end{equation}
Table \ref{tablec} might be helpful in keeping track of the asymptotic behavior of $c_{0,n}$ and $c_n$.

\begin{table}[!htbp]
\centering
\caption{\textit{Asymptotic behavior of $c_{0,n}$ and $c_n$ for different values of $z_n$}}
\vspace*{5mm}
\begin{tabular}{c|c|c|c}  
 \toprule
 $z_n$ & $c_{0,n}$ & $c_n$ & Other \\
 \hline
 $|z_n|\ll \sqrt{\rho_n/\beta_n}$ & $|c_{0,n}|\ll\beta_n^{1/2}/\rho_n^{3/2}$ & $c_n\asymp1$  &\\ 
 $z_n\gtrsim \sqrt{\rho_n/\beta_n}$, $L^*_n-z_n\gg\beta_n^{-1/3}$ & $c_{0,n}\gtrsim \beta_n^{1/2}/\rho_n^{3/2}$ & $\beta_n^{1/3}/\rho_n\ll c_{n}\leq 1$ & $c_nc_{0,n}\gtrsim \beta_n^{1/2}/\rho_n^{3/2}$\\ 
 $-z_n\gtrsim \sqrt{\rho_n/\beta_n}$, $z_n-L_n^\dagger\gg\beta_n^{-1/3}$ & $|c_{0,n}|\gtrsim \beta_n^{1/2}/\rho_n^{3/2}$& $3/2-c_n\gg\beta_n^{2/3}/\rho_n^2$  & $|c_nc_{0,n}|\gtrsim \beta_n^{1/2}/\rho_n^{3/2}$  \\ 
 \bottomrule
 \end{tabular}
 \label{tablec}
\end{table}

\bigskip
In the rest of this paper, to lighten the burden of notation, we will usually omit the subscript $n$ in the notation. For example, we will write $\rho$ in place of $\rho_n$, $z$ in place of $z_n$ and $g(z)$ in place of $g_n(z_n)$. However, it is important to keep in mind that these quantities do depend on $n$.

\subsection{Proof of Proposition \ref{Local}}
In this subsection, we will prove Proposition \ref{Local} using first and second moment estimates. First, we have the following lemma which shows that $y\in [z-l,z+l]$ satisfies the restriction (\ref{zn}) with $y$ in place of $z$.
\begin{Lemma}\label{YrestrictionZ}
For every $z$ satisfying (\ref{zn}), (\ref{zn1}) and (\ref{zn2}), choose $l$ according to (\ref{ln}). For all $y\in [z-l,z+l]$, we have
\begin{equation}\label{L^*-y}
L^*-y\gg\beta^{-1/3}
\end{equation}
and
\begin{equation}\label{y-Ldagger}
y-L^\dagger\gg\beta^{-1/3}.
\end{equation}
\end{Lemma}
Next, we have the following lemmas which control the difference between $t(z)$ and $t(y)$, and $g(z)$ and $g(y)$ for $y\in [z-l,z+l]$. 

\begin{Lemma}\label{Sy}
For every $z$ satisfying (\ref{zn}), (\ref{zn1}) and (\ref{zn2}), choose $l$ according to (\ref{ln}). For all $y\in [z-l,z+l]$, we have
\begin{equation}\label{sy}
|t(y)-t(z)|=o(\beta^{-2/3}).
\end{equation}
Moreover, uniformly for all $y\in [z-l,z+l]$, 
\begin{equation}\label{tyasymtz}
\lim_{n\rightarrow\infty}\frac{t(y)}{t(z)}=1.
\end{equation}
\end{Lemma}

\begin{Lemma}\label{Gy}
For every $z$ satisfying (\ref{zn}), (\ref{zn1}) and (\ref{zn2}), choose $l$ according to (\ref{ln}). Then for all $y\in [z-l,z+l]$, we have
\begin{equation}\label{gdifferconst}
|g(y)-g(z)|\lesssim 1.
\end{equation}
\end{Lemma}

The following lemma controls the first moment. 

\begin{Lemma}\label{Lemma2.4equivg}
For every $z$ satisfying (\ref{zn}), (\ref{zn1}) and (\ref{zn2}), choose $l$ according to (\ref{ln}). Let $s\geq0$ and
\[
t=t(z)-s, \quad x=L^*-w, \quad s_y=t(y)-t(z)+s.
\]
For all $w\in\mathbb{R}$, $s<t(z)$ and $y\in \mathbb{R}$, 
\begin{equation}\label{lem2.4egallwsmally}
p_t(x,y)\leq\frac{1}{\sqrt{2\pi t}}\exp \bigg(g(y)-\rho w+\beta ws_y-\frac{\beta^2}{6}s_y^3\bigg).
\end{equation}
Furthermore, if $s\asymp \beta^{-2/3}$, then for all $|w|\lesssim \beta^{-1/3}$ and $y\in [z-l,z+l]$,
\begin{equation}\label{lem2.4egsmallw}
p_t(x,y)=\frac{1}{\sqrt{2\pi t}}\exp \bigg(g(y)-\rho w+\beta ws_y-\frac{\beta^2}{6}s_y^3+o(1)\bigg).
\end{equation}
\end{Lemma}

A key step in the proof of Proposition \ref{Local} is the following second moment estimate. Note that it is rare for a particle to drift to the right of $L$ but once it does so, it will generate a large number of descendants in the interval $\mathcal{I}$ at time $t$, which ruins  the second moment argument. Therefore, we need to consider a truncated second moment estimate where particles are killed at $L$. For this process, we denote by $N_t^{L}(\mathcal{I})$ the number of particles in the interval $\mathcal{I}$ at time $t$.

\begin{Lemma}\label{Secondmomentg}
Consider the process which starts from a single particle at $x$ such that $0\leq L-x\lesssim \beta^{-1/3}$. For every $z$ satisfying (\ref{zn}), (\ref{zn1}) and (\ref{zn2}), choose $l$ according to (\ref{ln}). Consider intervals $\mathcal{I}$ defined in (\ref{I}). Suppose 
\[
s\asymp \beta^{-2/3},\qquad t=t(z)-s.
\]
Then for the process in which particles are killed upon hitting $L$, we have
\begin{align}\label{secondmoment}
E[N_t^L(\mathcal{I})^2]\lesssim \frac{\beta^{2/3}}{\rho^4}e^{\rho x+\rho L-2\rho L^*}\bigg(\int_{\mathcal{I}}\frac{1}{\sqrt{2\pi t(y)}}e^{g(y)}dy\bigg)^2.
\end{align}
\end{Lemma} 

To prove Proposition \ref{Local}, we need one more technical lemma. Let $\eta>0$.  Choose constants $C_5\geq 2$ and $C_6\geq -2^{-1/3}\gamma_1$ such that (\ref{C1cond1})-(\ref{C2}) hold. 
Then, $C_5$ and $C_6$ satisfy 
\begin{equation}\label{C1C2}
(1-\eta)e^{C_5^3/6}\leq \int_{0}^{2^{1/3}C_6}e^{2^{-1/3}C_5(\gamma_1+y)}Ai(\gamma_1+y)dy\leq (1+\eta)e^{C_5^3/6}.
\end{equation}
The next lemma is a slight generalization of Lemma 6.2 in \cite{RS2020}.

\begin{Lemma}\label{Unifgamma}
Suppose (\ref{assumpY}) and (\ref{assumpZ}) hold. Let $\eta>0$, and choose positive constants $C_5$ and $C_6$ such that (\ref{C1cond1})-(\ref{C2}) and (\ref{C1C2}) hold.
Let $s=C_5 \beta^{-2/3}$ and $u=t-t(z)+s$ for all $z$. If
\begin{equation}\label{unifgammaassump}
\beta^{-2/3}\log^{1/3}\Big(\frac{\rho}{\beta^{1/3}}\Big)\ll u\ll \frac{\rho}{\beta},
\end{equation}
then
\begin{align}\label{unifgammaz}
\lim_{n\rightarrow\infty}P\bigg( \frac{1-2\eta}{Ai'(\gamma_1)^2}Z(0)\leq \exp\bigg(\frac{\rho^2s}{2}-\frac{\beta^2s^3}{6}\bigg)\sum_{j\in\mathcal{N}_u}e^{(\rho-\beta s)X_j(u)}&1_{\{L-C_6\beta^{-1/3}<X_j(u)<L\}}\nonumber\\
&\leq  \frac{1+2\eta}{Ai'(\gamma_1)^2}Z(0)\bigg)=1.
\end{align}
Moreover, for every $z$ satisfying (\ref{zn}), (\ref{zn1}) and (\ref{zn2}), choose $l$ according to (\ref{ln}). For every $y\in[z-l,z+l]$, let $s_y=t(y)-t(z)+s$ and 
\begin{equation}\label{Gammaydef}
\Gamma_y=\exp\bigg(\frac{\rho^2s_y}{2}-\frac{\beta^2s_y^3}{6}\bigg)\sum_{j\in\mathcal{N}_u}e^{(\rho-\beta s_y)X_j(u)}1_{\{L-C_6\beta^{-1/3}<X_j(u)<L\}}.
\end{equation}
If (\ref{unifgammaassump}) holds, then uniformly for all $y\in [z-l,z+l]$,
\begin{equation}\label{unifgamma}
\lim_{n\rightarrow\infty}P\bigg(\frac{1-3\eta}{Ai'(\gamma_1)^2}Z(0)\leq \Gamma_y\leq \frac{1+3\eta}{Ai'(\gamma_1)^2}Z(0)\bigg)=1.
\end{equation}
\end{Lemma}

Lemmas \ref{Sy}, \ref{Lemma2.4equivg} and \ref{Unifgamma} will be proved in Section 4. Since the proof of Lemma \ref{Secondmomentg} is rather technical and tedious, we defer it until Section 5. With the help of the above lemmas, we will follow the same strategy as the proof of Proposition 2.2 of \cite{RS2020} to prove Proposition \ref{Local}.
\\

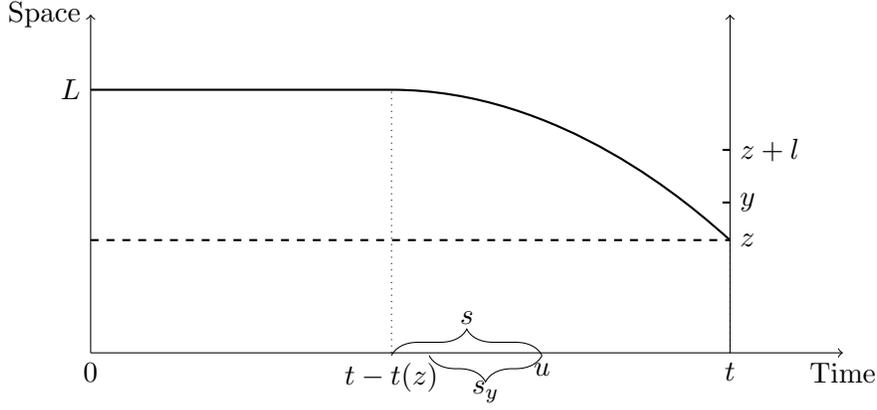
\begin{figure}[h]
\centering
\begin{tikzpicture}
\draw[->] (0,0) -- (10,0) node[anchor=north] {Time};
\draw	(0,0) node[anchor=north] {0}
		(4,0) node[anchor=north] {$t-t(z)$}
		(6,0) node[anchor=north] {$u$}
		(8.5,0) node[anchor=north] {$t$};

\draw[->] (0,0) -- (0,4.5) node[anchor=east] {Space};
\draw   (8.5,1.5) node[anchor=west]{$z$}
            (8.5,2) node[anchor=west]{$y$}
            (8.5,2.7) node[anchor=west]{$z+l$}
            (0,3.5) node[anchor=east]{$L$};

\draw[dotted] (4,0) -- (4,3.5);
\draw[dotted] (8.5,0) -- (8.5,1.5);
\draw[thick] (0,3.5) -- (4,3.5) parabola[bend at start](8.5,1.5);
\draw[thick] (8.4,2)-- (8.5,2);
\draw[thick] (8.4,2.7)-- (8.5,2.7);
\draw[->] (8.5,0) -- (8.5,4.5); 
\draw[thick,dashed] (0,1.5) -- (8.5,1.5);

\draw [decorate,decoration={brace,amplitude=10pt}, yshift=-1pt]
(4,0) -- (6,0)node [black,midway, above, yshift=8pt] {$s$};
\draw [decorate,decoration={brace,amplitude=10pt}, yshift=-1pt]
(6,0) -- (4.5,0)node [black,midway, below, yshift=-6pt] {
$s_y$};
\end{tikzpicture}
\caption{Notation when $z>0$}
\label{notation}
\end{figure}

\noindent\textit{Proof of Proposition \ref{Local}.}
Let 
\begin{equation}\label{defsu}
s=C_5\beta^{-2/3}, \qquad u=t-t(z)+s.
\end{equation}
By (\ref{t(z)>beta-3/2}), we have $t-t(z)<u<t$ for $n$ sufficiently large. For all $y\in \mathcal{I}$, denote
\[
s_y=t(y)-t(z)+s.
\]
By Lemma \ref{Sy}, for all $y\in \mathcal{I}$, we see that $s_y=(C_5\pm o(1))\beta^{-2/3}$ or equivalently, uniformly for all $y\in\mathcal{I}$,
\begin{equation}\label{syunif}
\lim_{n\rightarrow\infty} \beta^{2/3}s_y=C_5.
\end{equation} 
Figure \ref{notation} might be helpful for keeping track of notation.

Recall that $H(u)=L-\beta u^2/9$ by (\ref{H(u)}). For a particle $i\in\mathcal{N}_t$, recall that $\{X_i(v), 0\leq v\leq t\}$ denotes its past trajectory. Define
\begin{align*}
S_{1}&=\big\{i\in\mathcal{N}_t:X_{i}(u)\leq H(u), X_i(v)<L\;\text{for all}\;v\in [0,u]\big\},\\
S_{2}&=\Big\{i\in\mathcal{N}_t: H(u)<X_{i}(u)\leq L-C_6\beta^{-1/3}, X_i(v)<L\;\text{for all}\;v\in [0,u]\Big\},\\
S_{3}&=\Big\{i\in\mathcal{N}_t:L-C_6\beta^{-1/3}<X_{i}(u)< L, X_{i}(v)>L\;\text{for some}\;v\in (u,t)\Big\},\\
S_{4}&=\Big\{i\in\mathcal{N}_t:L-C_6\beta^{-1/3}<X_{i}(u)< L, X_{i}(v)\leq L\;\text{for all}\;v\in (u,t)\Big\},\\
S_{5}&=\mathcal{N}_t\setminus\big(S_1\cup S_2\cup S_3\cup S_4\big).
\end{align*}
For $j=1,..., 5$, write
\[
\Theta_{j}=\sum_{i\in S_{j}}1_{\{X_{i}(t)\in \mathcal{I}\}}.
\]
Then
\[
N_{t}(\mathcal{I})=\sum_{j=1}^{5}\Theta_{j}.
\]
We are going to show that the major contribution comes from $\Theta_4$, and $\Theta_4$ is concentrated around its mean. Define $(\mathcal{F}_t,t\geq 0)$ to be the natural filtration associated with the branching Brownian motion process.

Let us first consider $\Theta_1$. By inequality (\ref{lem2.4egallwsmally}) and Tonelli's theorem, we have
\begin{align}\label{Theta2first}
E[\Theta_{1}|\mathcal{F}_u]&=\sum_{j\in\mathcal{N}_u}\int_{\mathcal{I}}p_{t-u}(X_j(u),y)1_{\{\forall v\in [0,u], X_j(v)<L\}}1_{\{X_j(u)\leq H(u)\}}dy\nonumber\\
&\leq \int_{\mathcal{I}}\frac{1}{\sqrt{2\pi (t-u)}}\exp\bigg(g(y)-\rho L^*+\frac{\rho^2s_y}{2}-\frac{\beta^2s_y^3}{6}\bigg)\nonumber\\
&\qquad \times\sum_{j\in\mathcal{N}_u}e^{(\rho-\beta s_y)X_j(u)}1_{\{\forall v\in [0,u], X_j(v)<L\}}1_{\{X_j(u)\leq H(u)\}}dy.
\end{align}
We denote by $H_1$ the summation on the last line of (\ref{Theta2first}). By (\ref{assumptn}), we have $\beta^{-2/3}\ll u\ll \rho/\beta$. Furthermore, by equation (\ref{sy}), we see that $0\leq \beta s_y\ll \beta u$ for all $y\in \mathcal{I}$. Therefore, by (\ref{(5.6)}), we have for all $y\in \mathcal{I}$, 
\begin{equation}\label{unif5.6}
E\left[H_1|\mathcal{F}_0\right]\ll e^{-\beta^2u^3/73}Y(0).
\end{equation}
Since $t-u=t(z)-s\gg\beta^{-2/3}$, equation (\ref{tyasymtz}) implies that uniformly for all $y\in \mathcal{I}$,
\begin{equation}\label{ty/t-u}
\lim_{n\rightarrow\infty}\frac{t(y)}{t-u}=1.
\end{equation}
Also, for all $y\in \mathcal{I}$, we see from (\ref{assumptn}) that 
\begin{equation}\label{syu}
\rho^2s_y/2\ll \beta^2u^3.
\end{equation}
Thus, from (\ref{Theta2first})--(\ref{syu}), we have 
\[
E[\Theta_1|\mathcal{F}_0]\ll e^{-\rho L^*-\beta^2u^3/74}Y(0)\int_{\mathcal{I}}\frac{1}{\sqrt{2\pi t(y)}}e^{g(y)}dy.
\] 
Then by the conditional Markov's inequality, we can deduce that for any $\eta>0$, if $n$ is sufficiently large, 
\[
P\bigg(\Theta_1>\eta e^{-\rho L^*}Z(0)\int_{\mathcal{I}}\frac{1}{\sqrt{2\pi t(y)}}e^{g(y)}dy\bigg|\mathcal{F}_0\bigg)\leq \frac{Y(0)}{\eta Z(0)}e^{-\beta^2u^3/74}.
\]
Based on assumptions (\ref{assumpY}) and (\ref{assumpZ}), we have that for any $\eta>0$,
\begin{equation}\label{Theta2}
\lim_{n\rightarrow\infty}P\bigg(\Theta_1>\eta e^{-\rho L^*}Z(0)\int_{\mathcal{I}}\frac{1}{\sqrt{2\pi t(y)}}e^{g(y)}dy\bigg)=0.
\end{equation}

We next consider $\Theta_2$. By  (\ref{lem2.4egallwsmally}) and Tonelli's theorem again, we get
\begin{align}\label{Theta3first}
E[\Theta_{2}|\mathcal{F}_0]
&\leq \int_{\mathcal{I}}\frac{1}{\sqrt{2\pi (t-u)}}\exp\bigg(g(y)-\rho L^*+\frac{\rho^2s_y}{2}-\frac{\beta^2s_y^3}{6}\bigg)\nonumber\\
&\hspace{0.4in}\times E\left[\sum_{j\in\mathcal{N}_u}e^{(\rho-\beta s_y)X_j(u)}1_{\{\forall v\in [0,u], X_j(v)<L\}}1_{\{H(u)<X_j(u)\leq L-C_6\beta^{-1/3}\}}\bigg|\mathcal{F}_0\right]dy.
\end{align}
We can separate the expectation in the integrand into two parts by writing
\begin{align}\label{Theta3sumdef}
&E\left[\sum_{j\in\mathcal{N}_u}e^{(\rho-\beta s_y)X_j(u)}1_{\{\forall v\in [0,u], X_j(v)<L\}}1_{\{H(u)<X_j(u)\leq L-C_6\beta^{-1/3}\}}\bigg|\mathcal{F}_0\right]\nonumber\\
&\hspace{0.2in}= E\left[\sum_{j\in\mathcal{N}_u}e^{(\rho-\beta s_y)X_j(u)}1_{\{\forall v\in [0,u], X_j(v)<L\}}1_{\{ X_j(0)\leq H(u)\}}1_{\{H(u)<X_j(u)\leq L-C_6\beta^{-1/3}\}}\bigg|\mathcal{F}_0\right]\nonumber\\
&\hspace{0.4in}+E\left[\sum_{j\in\mathcal{N}_u}e^{(\rho-\beta s_y)X_j(u)}1_{\{\forall v\in [0,u], X_j(v)<L\}}1_{\{H(u)<X_j(0)<L\}}1_{\{H(u)<X_j(u)\leq L-C_6\beta^{-1/3}\}}\bigg|\mathcal{F}_0\right]\nonumber\\
&\hspace{0.2in}=:E[H_2|\mathcal{F}_0]+E[H_3|\mathcal{F}_0].
\end{align}
Note that when $H(u)<x<L$ and $H(u)<y<L-C_6\beta^{-1/3}$, 
\[
(2\beta)^{1/6}\big((L-x)^{1/2}+(L-y)^{1/2}\big)\leq (2\beta)^{1/6}\cdot 2\cdot (L-H(u))^{1/2}=2^{7/6}3^{-1}\beta^{2/3}u.
\]
Since $2^{7/6}3^{-1}<2^{-1/3}$ and $\beta^{2/3}u\gg 1$, equation (\ref{2.24}) is satisfied. Thus, in (\ref{Theta3sumdef}), we can upper bound the first expectation by (\ref{(5.5)}) and upper bound the second expectation by (\ref{lem2.5}). We have that for all $y\in\mathcal{I}$,
\begin{align}\label{Theta3sum}
E[H_2+H_3|\mathcal{F}_0]
&\leq e^{-\beta^2u^3/73}Y(0)+C_8Z(0)\int_{H(u)}^{L-C_6\beta^{-1/3}}e^{-\beta s_y v}\beta^{1/3}Ai\Big((2\beta)^{1/3}(L-v)+\gamma_1\Big)dv.
\end{align}
Substituting $v$ with $L-(2\beta)^{-1/3}r$, by (\ref{C2}) and (\ref{syunif}), we have for $n$ sufficiently large, for all ~$y\in \mathcal{I}$,
\begin{align}\label{Theta3dct}
&\int_{H(u)}^{L-C_6\beta^{-1/3}}e^{-\beta s_y v}\beta^{1/3}Ai\Big((2\beta)^{1/3}(L-v)+\gamma_1\Big)dv\nonumber\\
&\hspace{0.2in}\leq 2^{-1/3}e^{-\rho^2s_y/2}\int_{2^{1/3}C_6}^{\infty}e^{2^{-1/3}\beta^{2/3}s_y(\gamma_1+r)}Ai(\gamma_1+r)dr\nonumber\\
&\hspace{0.2in} \leq(1+\eta)2^{-1/3}e^{-\rho^2s_y/2}\int_{2^{1/3}C_6}^{\infty}e^{2^{-1/3}C_5(\gamma_1+r)}Ai(\gamma_1+r)dr\nonumber\\
&\hspace{0.2in}\leq (1+\eta)2^{-1/3}e^{-\rho^2s_y/2}\frac{\eta}{2}e^{C_5^3/48}.
\end{align}
Combining the above formula with (\ref{Theta3first}) and (\ref{Theta3sum}), we have
\begin{align*}
E[\Theta_2|\mathcal{F}_0]
&\leq e^{-\beta^2u^3/73}Y(0)e^{-\rho L^*}\int_{\mathcal{I}}\frac{1}{\sqrt{2\pi (t-u)}}\exp\bigg(g(y)+\frac{\rho^2s_y}{2}-\frac{\beta^2s_y^3}{6}\bigg)dy\\
&\hspace{0.2in}+\frac{C_8\eta(1+\eta)}{2^{4/3}}Z(0)e^{-\rho L^*}\int_{\mathcal{I}}\frac{1}{\sqrt{2\pi (t-u)}}\exp\bigg(g(y)-\frac{\beta^2s_y^3}{6}+\frac{C_5^3}{48}\bigg)dy.
\end{align*}
By (\ref{syu}), if $n$ is sufficiently large, then for all $y\in\mathcal{I}$,
\begin{equation}\label{rhosbetau}
\frac{\rho^2 s_y}{2}-\frac{\beta^2 u^3}{73}\leq -\frac{\beta^2 u^3}{74}.
\end{equation}
Also note that for $n$ large enough, $\beta^2s_y^3/6\geq C_5^3/48$ for all $y\in\mathcal{I}$. Then by (\ref{ty/t-u}), we obtain that for $n$ sufficiently large
\[
E[\Theta_2|\mathcal{F}_0] \leq \bigg(e^{-\beta^2u^3/74}Y(0)+\frac{C_8\eta(1+\eta)}{2^{4/3}}Z(0)\bigg)e^{-\rho L^*}\int_{\mathcal{I}}\frac{1}{\sqrt{2\pi t(y)}}e^{g(y)}dy.
\]
By (\ref{assumptn}), we have $u \gg \rho^{2/3}/\beta^{8/9}$, and therefore (\ref{assumpY}) and (\ref{assumpZ}) imply that $e^{-\beta^2 u^3/74} Y(0)/Z(0) \rightarrow_p 0$.  Thus, by the conditional Markov's inequality,
\begin{equation}\label{Theta3}
\limsup_{n\rightarrow\infty}P\bigg(\Theta_2>\eta^{1/2}e^{-\rho L^*}Z(0)\int_{\mathcal{I}}\frac{1}{\sqrt{2\pi t(y)}}e^{g(y)}dy\bigg)\leq \frac{C_8\eta^{1/2}(1+\eta)}{2^{4/3}}.
\end{equation}

We now consider $\Theta_3$ and $\Theta_4$. According to (\ref{lem2.4egsmallw}), (\ref{ty/t-u}) and Tonelli's theorem, there is a sequence of $\mathcal{F}_u$-measurable random variables $\{\theta_n\}_{n=1}^{\infty}$ which converges uniformly to $0$ as $n$ goes to infinity such that
\begin{align*}
E[\Theta_3+\Theta_4|\mathcal{F}_u]&=(1+\theta_n)\int_{\mathcal{I}}\frac{1}{\sqrt{2\pi t(y)}}\exp\bigg(g(y)-\rho L^*+\frac{\rho^2s_y}{2}-\frac{\beta^2s_y^3}{6}\bigg)\\
&\hspace{0.2in}\times\sum_{j\in\mathcal{N}_u}e^{(\rho-\beta s_y)X_j(u)}1_{\{L-C_6\beta^{-1/3}<X_j(u)<L\}}dy.
\end{align*}
We can deduce from (\ref{unifgamma}) that 
\begin{align}\label{Theta45}
\lim_{n\rightarrow\infty}P\bigg(\frac{1-4\eta}{Ai'(\gamma_1)^2}e^{-\rho L^*}Z(0)&\int_{\mathcal{I}}\frac{1}{\sqrt{2\pi t(y)}}e^{g(y)}dy\leq E[\Theta_3+\Theta_4|\mathcal{F}_u]\nonumber\\ 
&\leq \frac{1+4\eta}{Ai'(\gamma_1)^2}e^{-\rho L^*}Z(0)\int_{\mathcal{I}}\frac{1}{\sqrt{2\pi t(y)}}e^{g(y)}dy \bigg)=1.
\end{align}
Next, we will estimate $\Theta_3$ and $\Theta_4$ individually. Note that $\Theta_{3}$ accounts for particles that reach $L$ between times $u$ and $t$ and then drift back to $\mathcal{I}$. Consider a process which starts from a single particle at $L-C_6\beta^{-1/3}<x<L$. Suppose we kill particles upon hitting $L$. For $v\in [0,t-u]$, recall $r_x^{L}(v)$ is the rate at which particles hit $L$ at time $v$. We further denote by $m^z(v)$ the expected number of descendants in $\mathcal{I}$ at time $t$ of a particle that reaches $L$ at time $u+v$. Consider the process in which there is one particle at $x$ at time $u$ without killing. Then the expected number of particles in $\mathcal{I}$ at time $t$ whose trajectories cross $L$ between times $u$ and $t$ is 
\[
\int_{0}^{t-u}r_x^{L}(v)m^z(v)dv.
\]
From the definition of $m^z(v)$, we have
\[
m^z(v)=\int_{\mathcal{I}}p_{t-u-v}(L,y)dy.
\]
Setting $w=(2\beta)^{-1/3}\gamma_1$ and substituting $s+v$ in place of $s$ in equation (\ref{lem2.4egallwsmally}), we have
\begin{align}\label{mzv}
m^z(v)
&\leq \int_{\mathcal{I}}\frac{1}{\sqrt{2\pi (t-u-v)}}\exp\bigg(g(y)-\rho(2\beta)^{-1/3}\gamma_1+\beta(2\beta)^{-1/3}\gamma_1\big(t(y)-(t-u-v)\big)\nonumber\\
&\hspace{0.2in}-\frac{\beta^2}{6}\big(t(y)-(t-u-v)\big)^3\bigg)dy.
\end{align}
Note that $t(y)-(t-u-v)=s_y+v$. Combining (\ref{mzv}) with (\ref{rxLv}), we have
\begin{align*}
r_x^L(v)m^z(v)
&\leq \frac{C_9(L-x)}{v^{3/2}}\frac{1}{\sqrt{2\pi (t-u-v)}}\int_{\mathcal{I}}\exp\bigg(\rho x-\rho L^*-\frac{(L-x)^2}{2v}+g(y)\\
&\hspace{0.2in}+2^{-1/3}\beta^{2/3}\gamma_1s_y-\frac{\beta^2}{6}(s_y+v)^3\bigg)dy.
\end{align*}
We are going to deal with the terms involving $s_y$ by using an argument similar to the one leading to (\ref{Theta3dct}). Notice that $(s_y+v)^3\geq s_y^3$. By (\ref{syunif}), the dominated convergence theorem and Tonelli's theorem, for $n$ large enough, we obtain
\begin{align}\label{intrm}
&\int_0^{t-u}r_x^L(v)m^z(v)dv\nonumber\\
&\leq \frac{C_9(1+\eta)}{\sqrt{2\pi}}\exp\bigg(\rho x-\rho L^*+2^{-1/3}\beta^{2/3}\gamma_1s-\frac{\beta^2s^3}{6}\bigg)\int_{\mathcal{I}}e^{g(y)}\nonumber\\
&\hspace{0.2in}\times\bigg(\int_{0}^{(t-u)/2}\frac{L-x}{v^{3/2}}\frac{1}{\sqrt{(t-u-v)}}\exp\bigg(-\frac{(L-x)^2}{2v}\bigg)dv+\int_{(t-u)/2}^{t-u}\frac{L-x}{v^{3/2}}\frac{1}{\sqrt{t-u-v}}dv\bigg).
\end{align}
Lemma 4.1 in \cite{RS2020} states that for $a>0$ and $b>0$,
\begin{equation*}
\int_{0}^{\infty}\frac{1}{v^{3/2}}e^{-b^2/av}dv=\frac{\sqrt{\pi a}}{b}.
\end{equation*}
Thus we have
\begin{align}\label{1stintrm}
\int_{0}^{(t-u)/2}\frac{L-x}{v^{3/2}}\frac{1}{\sqrt{(t-u-v)}}\exp\bigg(-\frac{(L-x)^2}{2v}\bigg)dv&\leq \sqrt{\frac{2}{t-u}}\int_{0}^{\infty}\frac{L-x}{v^{3/2}}\exp\bigg(-\frac{(L-x)^2}{2v}\bigg)dv\nonumber\\
&= \frac{2\sqrt{\pi}}{\sqrt{t-u}}.
\end{align}
Noticing that $L-x\leq C_6\beta^{-1/3}\ll  \sqrt{t-u}$, we see that for $n$ sufficiently large
\begin{equation}\label{2ndintrm}
\int_{(t-u)/2}^{t-u}\frac{L-x}{v^{3/2}\sqrt{t-u-v}}dv\leq \frac{2^{3/2}(L-x)}{(t-u)^{3/2}}\int_{(t-u)/2}^{t-u}\frac{1}{\sqrt{t-u-v}}dv=\frac{4(L-x)}{t-u}\leq \frac{1}{\sqrt{t-u}}.
\end{equation}
By equations (\ref{ty/t-u}), (\ref{intrm}), (\ref{1stintrm}) and (\ref{2ndintrm}), since $\gamma_1<0$, we have for $n$ large enough,
\begin{align*}
&\int_0^{t-u}r_x^L(v)m^z(v)dv\\
&\quad \leq \big(2\sqrt{\pi}+1\big)C_9(1+\eta)\exp\bigg(\rho x-\rho L^*+2^{-1/3}\beta^{2/3}\gamma_1s-\frac{\beta^2s^3}{6}\bigg)\int_{\mathcal{I}}\frac{1}{\sqrt{2\pi t(y)}}e^{g(y)}dy\\
&\quad \leq \big(2\sqrt{\pi}+1\big)C_9(1+\eta)\exp\bigg(\rho x-\rho L^*-\frac{\beta^2s^3}{6}\bigg)\int_{\mathcal{I}}\frac{1}{\sqrt{2\pi t(y)}}e^{g(y)}dy.
\end{align*}
Summing over all particles at time $u$, we have for $n$ large enough
\begin{equation}\label{Theta4condFu}
E[\Theta_3|\mathcal{F}_u]\leq \big(2\sqrt{\pi}+1\big)C_9(1+\eta)Y(u)\exp\bigg(-\rho L^*-\frac{\beta^2s^3}{6}\bigg)\int_{\mathcal{I}}\frac{1}{\sqrt{2\pi t(y)}}e^{g(y)}dy.
\end{equation}
Furthermore, equation (6.31) in \cite{RS2020} states that
\begin{equation}\label{6.31inRS2020}
\lim_{n\rightarrow\infty}P\bigg(Y_n(u)<\frac{1+\eta}{Ai'(\gamma_1)^2}Z_n(0)\int_{0}^{\infty}Ai(\gamma_1+z)dz\bigg)=1.
\end{equation}
Recall that $s=C_5\beta^{-2/3}$. Combining (\ref{Theta4condFu}) with (\ref{C1cond2}) and (\ref{6.31inRS2020}), we have
\begin{equation}\label{Theta4cond}
\lim_{n\rightarrow\infty}P\bigg(E[\Theta_3|\mathcal{F}_u]\geq \Big(\sqrt{2}+\frac{1}{\sqrt{2\pi}}\Big)C_9\eta(1+\eta) e^{-\rho L^*}Z(0)\int_{\mathcal{I}}\frac{1}{\sqrt{2\pi t(y)}}e^{g(y)}dy\bigg)=0.
\end{equation}
By the conditional Markov's inequality,
\begin{equation}\label{Theta4}
\limsup_{n\rightarrow\infty}P\bigg(\Theta_3>\eta^{1/2}e^{-\rho L^*}Z(0)\int_{\mathcal{I}}\frac{1}{\sqrt{2\pi t_y}}e^{g(y)}dy\bigg)\leq \bigg(\sqrt{2}+\frac{1}{\sqrt{2\pi}}\bigg)C_9\eta^{1/2}(1+\eta).
\end{equation}

Next, we are going to show that $\Theta_4$ is concentrated around its mean. By (\ref{Theta45}) and (\ref{Theta4cond}), letting
\[
C(\eta)=4\eta+ \Big(\sqrt{2}+\frac{1}{\sqrt{2\pi}}\Big)C_9\eta(1+\eta)Ai'(\gamma_1)^2,
\]
we see that
\begin{equation}\label{Theta5lb}
\lim_{n\rightarrow\infty}P\bigg(E[\Theta_4|\mathcal{F}_u]\geq \frac{1-C(\eta)}{Ai'(\gamma_1)^2} e^{-\rho L^*}Z(0)\int_{\mathcal{I}}\frac{1}{\sqrt{2\pi t(y)}}e^{g(y)}dy\bigg)=1.
\end{equation}
Considering the process in which particles are killed upon hitting $L$, we can bound the conditional variance of $\Theta_4$ by Lemma \ref{Secondmomentg}. We get
\begin{align*}
\Var(\Theta_4|\mathcal{F}_u)&\leq\sum_{j\in\mathcal{N}_u} E_{X_j(u)}\big[N_t^L(\mathcal{I})^2\big]1_{\{L-C_6\beta^{-1/3}<X_j(u)<L\}}\\
&\lesssim \sum_{j\in\mathcal{N}_u} \frac{\beta^{2/3}}{\rho^4}e^{\rho X_i(u)+\rho L-2\rho L^*}\bigg(\int_{\mathcal{I}}\frac{1}{\sqrt{2\pi t(y)}}e^{g(y)}dy\bigg)^21_{\{L-C_6\beta^{-1/3}<X_j(u)<L\}}\\
&\leq \frac{\beta^{2/3}}{\rho^4}e^{-2\rho L^*+\rho L}Y(u)\bigg(\int_{\mathcal{I}}\frac{1}{\sqrt{2\pi t(y)}}e^{g(y)}dy\bigg)^2.
\end{align*}
Then by Chebyshev's inequality,
\begin{equation}\label{Chebyshevg}
P\bigg(\Big|\Theta_4-E[\Theta_4|\mathcal{F}_u]\Big|>\eta e^{-\rho L^*}Z(0)\int_{\mathcal{I}}\frac{1}{\sqrt{2\pi t(y)}}e^{g(y)}dy\bigg|\mathcal{F}_u\bigg)\leq \frac{\beta^{2/3}e^{\rho L}}{\eta^2\rho^4}\frac{Y(u)}{Z(0)^2}.
\end{equation}
On account of (\ref{assumpZ}) and (\ref{6.31inRS2020}), as $n\rightarrow\infty$,
\begin{equation}\label{YZconvto0}
\frac{\beta^{2/3}e^{\rho L}}{\rho^4}\frac{Y(u)}{Z(0)^2}\rightarrow_{p}0.
\end{equation}
As a result, by equations (\ref{Theta45}), (\ref{Theta5lb}), (\ref{Chebyshevg}) and (\ref{YZconvto0}), we obtain
\begin{align}\label{Theta5}
\lim_{n\rightarrow\infty}P\bigg(\bigg|\Theta_4-\frac{1}{Ai'(\gamma_1)^2}e^{-\rho L^*}&Z(0)\int_{\mathcal{I}}\frac{1}{\sqrt{2\pi t(y)}}e^{g(y)}dy\bigg|\nonumber\\
&\leq\bigg(\eta+\frac{C(\eta)}{Ai'(\gamma_1)^2}\bigg)e^{-\rho L^*}Z(0)\int_{\mathcal{I}}\frac{1}{\sqrt{2\pi t(y)}}e^{g(y)}dy\bigg)=1.
\end{align}

It remains to consider $\Theta_5$. Define $S_5^*$ to be the set consists of particles whose trajectories cross $L$ before time $u$, so
\[
S^*_5=\big\{i\in\mathcal{N}_t:X_i(v)\geq L\;\text{for some}\;v\in[0,u]\big\}.
\] 
We observe that $S_5\subseteq S_5^*$. Note that $u+2C_7\rho^{-2}\leq t$ since $t-u=t(z)-s\gg\beta^{-2/3}$ by (\ref{t(z)>beta-3/2}). According to \ref{RS3}, the probability that particles that either are to the right of $L$ at time $0$ or hit $L$ before time $u$ have descendants alive at time $t$ goes to $0$ as $n$ goes to infinity. Therefore,
\begin{equation}\label{Theta1}
\lim_{n\rightarrow \infty}P(\Theta_5=0)=1.
\end{equation}

Consequently, for any $\kappa>0$, by choosing $\eta$ appropriately, equation (\ref{mainpropl}) follows from (\ref{Theta2}), (\ref{Theta3}), (\ref{Theta4}), (\ref{Theta5}) and (\ref{Theta1}). 
\qedwhite

\subsection{Proof of Proposition \ref{Global}}
In this subsection, we will prove Proposition \ref{Global} with the help of Proposition \ref{Local}. Before starting the proof, we need two more lemmas to control the number of particles that are far away from $z$.

\begin{Lemma}\label{Integralofg}
Consider $z$ such that (\ref{zn}) holds and $|z|\gtrsim\sqrt{\rho/\beta}$. For any $\eta>0$, there exists a constant $C_{11}$ large enough such that for 
\begin{equation}\label{dl}
 d= \frac{C_{11}}{|c_0|\rho},
\end{equation}
the following hold:
\begin{enumerate}
\item The constant satisfies
\begin{equation}\label{C6eta}
C_{11}>4,\qquad e^{-C_{11}/2}<\eta.
\end{equation}
\item If $z> 0$ for all $n$, then for $n$ sufficiently large
\begin{equation}\label{lemintegralofg1}
\int_{z+d}^{L^*}\frac{1}{\sqrt{2\pi t(y)}}e^{g(y)}dy<\eta \int_{z}^{z+d}\frac{1}{\sqrt{2\pi t(y)}}e^{g(y)}dy.
\end{equation}
\item If $z< 0$ for all $n$, then for $n$ sufficiently large
\begin{equation}\label{lemintegralofg2}
\int_{-\infty}^{z-d}\frac{1}{\sqrt{2\pi t(y)}}e^{g(y)}dy<\eta \int_{z-d}^{z}\frac{1}{\sqrt{2\pi t(y)}}e^{g(y)}dy.
\end{equation}
\end{enumerate}
\end{Lemma}
Since we expect that the density of the number of particles near $y$ is roughly proportional to $e^{g(y)}/\sqrt{2\pi t(y)}$, Lemma \ref{Integralofg} indicates that most particles in $[z, L^*]$ are in $[z, z+d]$, while most particles in $(-\infty, z]$ are in $[z-d, z]$.

Suppose (\ref{zn}) holds and $|z|\gtrsim\sqrt{\rho/\beta}$. For any $\eta>0$,  choose $d$ according to Lemma \ref{Integralofg}. Note that (\ref{ln}) holds with $d$ in place of $l$, and for $l$ satisfying (\ref{lng}), we have $2d<l$ for $n$ sufficiently large. Denote 
\begin{equation*}
\zeta=\begin{cases} 
     z+2d&\qquad \text{if}\;z> 0,\\
     z-2d  &\qquad \text{if}\;z< 0,
   \end{cases}
\end{equation*}

\begin{Lemma}\label{Lemma2.4equivd}
Consider $z$ such that (\ref{zn2}) holds and $|z|\gtrsim\sqrt{\rho/\beta}$.  Let 
\[
s\asymp \beta^{-2/3},\quad t=t(z)-s, \quad x\leq L.
\]
\begin{enumerate}
\item Suppose $z>0$ for all $n$, and $z$ satisfies (\ref{zn}).  Choose $d$ according to (\ref{dl}). For all $y\in [\zeta,\infty)$, we have for $n$ sufficiently large,
\begin{equation}\label{lem2.4d1}
p_{t}(x,y)\leq p_{t}(x,\zeta)\exp\bigg(-\frac{\rho}{2}(1-c)(y-\zeta)\bigg).
\end{equation}
\item Suppose $z< 0$ for all $n$. Write $x=L^*-w$. For all $y$, we have for $n$ sufficiently large,
\begin{equation}\label{lem2.4d2}
p_{t}(x,y)\leq \frac{1}{\sqrt{2\pi t}}\exp\bigg(g(z)-(c-1)\rho (z-y)-\rho w +\beta sw-\frac{\beta^2s^3}{6}\bigg).
\end{equation}
Furthermore, let $s_\zeta=t(\zeta)-t(z)+s$. If $y\leq \zeta$, then for $n$ sufficiently large,
\begin{equation}\label{lem2.4d2var}
p_{t}(x,y)\leq \frac{1}{\sqrt{2\pi t}}\exp\bigg(g(\zeta)-(c-1)\rho (\zeta-y)-\rho w +\beta s_{\zeta}w-\frac{\beta^2s_{\zeta}^3}{6}\bigg).
\end{equation}
\end{enumerate}
\end{Lemma}

\noindent\textit{Proof of Proposition \ref{Global}.}
Let us first consider the case $|z|\gtrsim \sqrt{\rho/\beta}$. Define $s$, $u$ and $H(u)$ as in (\ref{defsu}) and (\ref{H(u)}). Denote
\begin{equation*}
 \mathcal{K}=\begin{cases} 
     [z,\zeta] &\qquad \text{if}\;z> 0,\\
     [\zeta,z]  &\qquad \text{if}\;z< 0.
   \end{cases}
\end{equation*}
 Define
\begin{align*}
S_{1}&=\big\{i\in\mathcal{N}_t:X_{i}(t)\in \mathcal{K}\big\},\\
S_{2}&=\big\{i\in \mathcal{N}_t\setminus S_1:X_{i}(v)\geq L\;\text{for some}\; v\in [0,u]\big\},\\
S_{3}&=\big\{i\in\mathcal{N}_t\setminus (S_1\cup S_2): X_i(u)\leq L-C_6\beta^{-1/3}\big\},\\
S_{4}&=\big\{i\in\mathcal{N}_t\setminus (S_1\cup S_2): L-C_6\beta^{-1/3}<X_i(u)<L\big\}.
\end{align*}
For $j=1, ..., 4$, write
\[
\Xi_{j}=\sum_{i\in S_{j}}1_{\{X_{i}(t)\in \mathcal{J}\}}.
\]
Then
\[
N_{t}(\mathcal{J})=\sum_{j=1}^{4}\Xi_{j}.
\]
We will show that compared with $\Xi_1$, the terms $\Xi_2$, $\Xi_3$ and $\Xi_4$ are negligible.

We first consider $\Xi_1$. Since $2d$ satisfies the restriction (\ref{ln}) in Proposition \ref{Local}, according to (\ref{mainpropl}), (\ref{lemintegralofg1}) and (\ref{lemintegralofg2})
\begin{align}\label{Xi1}
\lim_{n\rightarrow\infty}P\bigg(\frac{1-\eta}{(1+\eta)Ai'(\gamma_1)^2}e^{-\rho L^*}Z(0)&\int_{\mathcal{J}}\frac{1}{\sqrt{2\pi t(y)}}e^{g(y)}dy\leq \Xi_1\nonumber\\
&\leq\frac{1+\eta}{Ai'(\gamma_1)^2}e^{-\rho L^*}Z(0)\int_{\mathcal{J}}\frac{1}{\sqrt{2\pi t(y)}}e^{g(y)}dy\bigg)=1.
\end{align}

For $\Xi_2$, since $u+C_7\rho^{-2}\leq t$, according to \ref{RS3}, the probability that particles that either are to the right of $L$ at time $0$ or hit $L$ before time $u$ have descendants alive at time $t$ goes to $0$ as $n$ goes to infinity.  Therefore,
\begin{equation}\label{Xi2}
\lim_{n\rightarrow \infty}P(\Xi_2=0)=1.
\end{equation}

It remains to consider $\Xi_3$ and $\Xi_4$. Let us first consider the case when $z>0$. Recall the definition of $H_1$ in (\ref{Theta2first}) and $H_2$, $H_3$ in (\ref{Theta3sumdef}). Since (\ref{ln}) holds with $2d$ in place of $l$, by inequalities (\ref{lem2.4egallwsmally}), (\ref{lem2.4d1}) and Tonelli's theorem, for $n$ large enough, we have
\begin{align}\label{Xi3firstmoment1}
E[\Xi_{3}|\mathcal{F}_u]&=\sum_{j\in\mathcal{N}_u}\int_{\zeta}^{z+l}p_{t-u}(X_j(u),y)1_{\{\forall v\in [0,u],X_{j}(v)<L\}}1_{\{X_j(u)\leq L-C_6\beta^{-1/3}\}}dy\nonumber\\
&\leq \sum_{j\in\mathcal{N}_u}p_{t-u}(X_j(u),\zeta)1_{\{\forall v\in [0,u],X_{j}(v)<L\}}1_{\{X_j(u)\leq L-C_6\beta^{-1/3}\}}\int_{\zeta}^{z+l}e^{-\rho(1-c)(y-\zeta)/2}dy\nonumber\\
&\leq \frac{1}{\sqrt{2\pi (t-u)}}\exp\bigg(g(\zeta)-\rho L^*+\frac{\rho^2s_{\zeta}}{2}-\frac{\beta^2s_{\zeta}^3}{6}\bigg)\int_{\zeta}^{z+l}e^{-\rho(1-c)(y-\zeta)/2}dy\nonumber\\
&\hspace{0.2in}\times(H_1+H_2+H_3)\end{align}
According to the choice of $d$, since $c\in (0,1)$ when $z > 0$, we have
\[
 \frac{2}{\rho(1-c)}=\frac{2(1+c)}{\rho c_0}\leq \frac{4}{\rho c_0}<d.
\]
Also, by (\ref{tyasymtz}), uniformly for all $y\in [z,\zeta]$, we have $t-u>t(y)/(1+\eta)$ for sufficiently large $n$. Combining this observation with the fact that $g(y)$ is decreasing on $(0,L^*)$ and (\ref{lemintegralofg1}), we have
\begin{align}\label{remainint}
\frac{e^{g(\zeta)}}{\sqrt{2\pi(t-u)}}\int_{\zeta}^{z+l}e^{-\rho(1-c)(y-\zeta)/2}dy
&\leq \frac{1}{\sqrt{2\pi(t-u)}}\frac{2}{\rho(1-c)}e^{g(\zeta)}\nonumber\\
&\leq \int_{\zeta-d}^{\zeta}\frac{\sqrt{1+\eta}}{\sqrt{2\pi t(y)}}e^{g(y)}dy\nonumber\\
&\leq \eta\sqrt{1+\eta}\int_{\mathcal{J}\cap(-\infty,L^*]}\frac{1}{\sqrt{2\pi t(y)}}e^{g(y)}dy.
\end{align}
Since (\ref{ln}) in Proposition \ref{Local} holds with $2d$ in place of $l$, equations (\ref{unif5.6}), (\ref{Theta3sum}) and (\ref{Theta3dct}) hold with $s_{\zeta}$ in place of $s_{y}$. By (\ref{unif5.6}), (\ref{Theta3sum}), (\ref{Theta3dct}), (\ref{Xi3firstmoment1}) and (\ref{remainint}), along with (\ref{rhosbetau}) with $\zeta$ in place of $y$, we get
\begin{align}\label{Xi31pre}
E[\Xi_{3}|\mathcal{F}_0]&< \eta\sqrt{1+\eta}e^{-\rho L^*} \bigg(2e^{-\beta^2u^3/74}Y(0)+\frac{C_8\eta(1+\eta)}{2^{4/3}}Z(0)\bigg)\int_{\mathcal{J}\cap(-\infty,L^*]}\frac{1}{\sqrt{2\pi t(y)}}e^{g(y)}dy.
\end{align}
Therefore, equations (\ref{assumpY}), (\ref{assumpZ}), (\ref{Xi31pre}) and the conditional Markov's inequality imply
\begin{equation}\label{Xi31}
\limsup_{n\rightarrow\infty}P\bigg(\Xi_3>\sqrt{\eta(1+\eta)}e^{-\rho L^*}Z(0)\int_{\mathcal{J}\cap(-\infty,L^*]}\frac{1}{\sqrt{2\pi t(y)}}e^{g(y)}dy\bigg)\leq \frac{C_8\eta^{3/2}(1+\eta)}{2^{4/3}}.
\end{equation}
As for $\Xi_4$, according to the argument leading to (\ref{Xi3firstmoment1}), we have
\[
E[\Xi_4|\mathcal{F}_u]\leq  \frac{1}{\sqrt{2\pi (t-u)}}e^{g(\zeta)-\rho L^*}\int_{\zeta}^{z+l}e^{-\rho(1-c)(y-\zeta)/2}dy\times\Gamma_\zeta,
\]
where $\Gamma_\zeta$ was defined in (\ref{Gammaydef}).
By equations (\ref{unifgamma}) and  (\ref{remainint}), and the conditional Markov's inequality, we get
\begin{equation}\label{Xi41}
\limsup_{n\rightarrow\infty}P\bigg(\Xi_4>\sqrt{\eta(1+\eta)}e^{-\rho L^*}Z(0)\int_{\mathcal{J}\cap(-\infty,L^*]}\frac{1}{\sqrt{2\pi t(y)}}e^{g(y)}dy\bigg)
\leq \frac{(1+3\eta)\sqrt{\eta}}{Ai'(\gamma_1)^2}.
\end{equation}
As a result, when $z>0$, for any $\kappa>0$, by choosing $\eta$ appropriately, equation (\ref{mainpropg}) follows from (\ref{Xi1}), (\ref{Xi2}), (\ref{Xi31}) and (\ref{Xi41}). 

When $z<0$, by (\ref{lem2.4d2var}) and Tonelli's theorem, for $n$ large enough, we have
\begin{align}\label{Xi3firstmoment2}
E[\Xi_{3}|\mathcal{F}_u]&=\sum_{j\in\mathcal{N}_u}\int_{z-l}^{\zeta}p_{t-u}(X_j(u),y)1_{\{\forall v\in [0,u], X_j(v)<L\}}1_{\{X_j(u)\leq L-C_6\beta^{-1/3}\}}dy\nonumber\\
&\leq \frac{1}{\sqrt{2\pi (t-u)}}\exp\bigg(g(\zeta)-\rho L^*+\frac{\rho^2s_{\zeta}}{2}-\frac{\beta^2s_{\zeta}^3}{6}\bigg)\int_{z-l}^{\zeta}e^{-\rho(c-1)(\zeta-y)}dy\nonumber\\
&\hspace{0.2in}\times(H_1+H_2+H_3),
\end{align}
where $H_1$, $H_2$, and $H_3$ are defined as in \eqref{Theta2first} and \eqref{Theta3sumdef} but with $\zeta$ in place of $y$.
By (\ref{tyasymtz}), uniformly for all $y\in [\zeta,z]$, we have $t-u>t(y)/(1+\eta)$ for sufficiently large $n$. Combining this observation with the fact that $g(y)$ is increasing on $(-\infty,0)$, we have
\begin{align}\label{remainint2}
\frac{e^{g(\zeta)}}{\sqrt{2\pi(t-u)}}\int_{z-l}^{\zeta}e^{-\rho(c-1)(\zeta-y)}dy
&\leq \frac{1}{\sqrt{2\pi(t-u)}}\frac{e^{g(\zeta)}}{\rho(c-1)}\leq \frac{1}{d\rho(c-1)}\int_{\zeta}^{\zeta+d}\frac{\sqrt{1+\eta}}{\sqrt{2\pi t(y)}}e^{g(y)}dy.
\end{align}
By (\ref{range}), (\ref{1-c=c_0}) and (\ref{C6eta}), we see that
\begin{equation}\label{(c-1)rhod}
(c-1)\rho d=\frac{(c-1)C_{11}}{|c_0|}=\frac{C_{11}}{1+c}>\frac{4}{1+3/2}=\frac{8}{5}.
\end{equation}
Also, by (\ref{lemintegralofg2}), we have
\begin{equation}\label{3integrals}
\int_{\zeta}^{\zeta+d} \frac{1}{\sqrt{2 \pi t(y)}} e^{g(y)} \: dy \leq \int_{-\infty}^{z-d} \frac{1}{\sqrt{2 \pi t(y)}} e^{g(y)} \: dy \leq \eta \int_{z-d}^z \frac{1}{\sqrt{2 \pi t(y)}} e^{g(y)} \: dy.
\end{equation}
Therefore, by (\ref{remainint2}), (\ref{(c-1)rhod}) and (\ref{3integrals}), we obtain that for $n$ sufficiently large,
\begin{equation}\label{remaint3}
\frac{e^{g(\zeta)}}{\sqrt{2\pi (t-u)}}\int_{z-l}^{\zeta}e^{-\rho(c-1)(\zeta-y)dy}\leq \frac{5\eta\sqrt{1+\eta}}{8}\int_{\mathcal{J}\cap(-\infty, L^*]}\frac{1}{\sqrt{2\pi t(y)}}e^{g(y)}dy.
\end{equation}
Since $2d$ satisfies the restriction (\ref{ln}) in Proposition \ref{Local}, equations (\ref{unif5.6}), (\ref{Theta3sum}) and (\ref{Theta3dct}) hold with $s_{\zeta}$ in place of $s_y$.
By (\ref{unif5.6}), (\ref{Theta3sum}), (\ref{Theta3dct}), (\ref{Xi3firstmoment2}) and (\ref{remaint3}), we get
\begin{align}\label{Xi32pre}
E[\Xi_{3}|\mathcal{F}_0]&<\frac{5\eta\sqrt{1+\eta}}{8}e^{-\rho L^*} \bigg(2e^{-\beta^2u^3/74}Y(0)+\frac{C_8\eta(1+\eta)}{2^{4/3}}Z(0)\bigg)\int_{\mathcal{J}\cap(-\infty,L^*]}\frac{1}{\sqrt{2\pi t(y)}}e^{g(y)}dy.
\end{align}
Therefore, equations (\ref{assumpY}), (\ref{assumpZ}), (\ref{Xi32pre}) and the conditional Markov's inequality imply
\begin{equation}\label{Xi32}
\limsup_{n\rightarrow\infty}P\bigg(\Xi_3>\frac{5\sqrt{\eta(1+\eta)}}{8}e^{-\rho L^*}Z(0)\int_{\mathcal{J}\cap(-\infty,L^*]}\frac{1}{\sqrt{2\pi t(y)}}e^{g(y)}dy\bigg)\leq \frac{C_8\eta^{3/2}(1+\eta)}{2^{4/3}}.
\end{equation}
As for $\Xi_4$, according to the argument leading to (\ref{Xi3firstmoment2}), we have
\[
E[\Xi_4|\mathcal{F}_u]\leq  \frac{1}{\sqrt{2\pi (t-u)}}e^{g(\zeta)-\rho L^*}\int_{z-l}^{\zeta}e^{-\rho(c-1)(z-y)}dy\times\Gamma_\zeta.
\]
By equations (\ref{unifgamma}) and (\ref{remaint3}), and the conditional Markov's inequality, we get
\begin{equation}\label{Xi42}
\limsup_{n\rightarrow\infty}P\bigg(\Xi_4>\frac{5\sqrt{\eta(1+\eta)}}{8}e^{-\rho L^*}Z(0)\int_{\mathcal{J}\cap(-\infty,L^*]}\frac{1}{\sqrt{2\pi t(y)}}e^{g(y)}dy\bigg)
\leq \frac{(1+3\eta)\sqrt{\eta}}{Ai'(\gamma_1)^2}.
\end{equation}
As a result, when $z<0$, for any $\kappa>0$, by choosing $\eta$ appropriately, equation (\ref{mainpropg}) follows from (\ref{Xi1}), (\ref{Xi2}), (\ref{Xi32}) and (\ref{Xi42}). 

It remains to consider the case $|z|\ll \sqrt{\rho/\beta}$. Below we will only prove the result under the scenario $z>0$. The scenario when $z < 0$ can be proved using the same argument. The interval $[z,z+l]$ can be divided into two intervals
\[
[z,z+l]=\Big[z,z+\sqrt{\frac{\rho}{\beta}}\Big]\cup \Big[z+\sqrt{\frac{\rho}{\beta}},z+l\Big].
\]
It is obvious that the first interval fits in the setting of Proposition \ref{Local}. We further claim that the second interval fits in the setting of the previous case. Indeed, according to Lemma \ref{Sy}, we know that
\[
t(z)-t\Big(z+\sqrt{\frac{\rho}{\beta}}\Big)=o(\beta^{-2/3}).
\]
Thus (\ref{assumptn}) holds with $z+\sqrt{\rho/\beta}$ in place of $z$. Also, letting $c_0^*=(z+\sqrt{\rho/\beta})/L^*$, which is the same as $c_0$ but with $z+\sqrt{\rho/\beta}$ in place of $z$,
the length of the second interval satisfies
\[
l-\sqrt{\frac{\rho}{\beta}}\gg \sqrt{\frac{\rho}{\beta}}\asymp \frac{1}{c_0^*\rho}.
\]
According to Proposition \ref{Local} and the previous case, equation (\ref{mainpropl}) holds with $[z,z+\sqrt{\rho/\beta}]$ in place of $\mathcal{I}$ and equation (\ref{mainpropg}) holds with $[z+\sqrt{\rho/\beta},z+l]$ in place of $\mathcal{J}$. Combining these two equations, (\ref{mainpropg}) follows.
\qedwhite

\subsection{Proof of Proposition \ref{Corollaryr}}
In this subection, we will prove Proposition \ref{Corollaryr}, which gives the maximal displacement of the process. For any constant $C_2\in\mathbb{R}$, define
\begin{equation}\label{AC'}
A = \begin{cases}
              0 & C_2>0 \\
              1 & C_2=0\\
              -2C_2 & C_2<0
       \end{cases} \quad \text{and} \quad 
C'_2 = \begin{cases}
             C_2  & C_2>0 \\
             1  & C_2=0\\
             -C_2 & C_2<0.
       \end{cases}
\end{equation}
The proof of Proposition \ref{Corollaryr} requires the following lemma which concerns the maximal displacement of a slightly supercritical branching Brownian motion with constant branching rate. \begin{Lemma}\label{Dominant}
Consider a branching Brownian motion started from a single particle at $L^A$. Each particle moves as standard Brownian motion. Each particle independently dies at rate $d(2L)$, and splits into two particles at rate $b(2L)$. Let $M^*_t$ be the maximal position that is ever reached by a particle before time $t$. For any constant $C_2\in\mathbb{R}$, define $C'_2>0$ as in (\ref{AC'}). There exists a constant $C_{12}$ such that if $n$ is sufficiently large, then for all $t$,
\begin{equation}\label{domprob}
P\bigg(M_t^*>L^A+\frac{C'_2}{\rho}\bigg)\leq C_{12}\rho^2.
\end{equation}
\end{Lemma}

\noindent\textit{Proof of Lemma \ref{Dominant}.}
In this process, each individual lives for an exponentially distributed time with parameter $b(2L)+d(2L)$, and then gives birth to $0$ offspring with probability $d(2L)/(b(2L)+d(2L))$ and $2$ offspring with probability $b(2L)/(b(2L)+d(2L))$. Therefore, the generating function for the offspring distribution is
\[
f(s)=\frac{d(2L)}{b(2L)+d(2L)}+\frac{b(2L)}{b(2L)+d(2L)}s^2.
\] Let $B$ be the event of survival. By (\ref{assump3}) and the formula for the survival probability of the Galton-Watson process, there exists a constant $C_{13}$ such that for all $n$,
\[
P(B)=\frac{b(2L)-d(2L)}{b(2L)}\leq C_{13}\rho^2.
\]
For any time $t$, we get
\begin{equation}\label{domcond}
P\bigg(M_t^*>L^A+\frac{C'_2}{\rho}\bigg)\leq P\bigg(M_t^*>L^A+\frac{C'_2}{\rho}\bigg|B^c\bigg)P(B^c)+P(B)\leq P\bigg(M_t^*>L^A+\frac{C'_2}{\rho}\bigg|B^c\bigg)+C_{15}\rho^2.
\end{equation}
We are interested in the behavior of the process conditioned on the event $B^c$ of extinction. According to equation (4) of Gadag and Rajarshi \cite{Gadag}, the conditioned process is equivalent to a subcritical branching process with generating function
\[
\hat f(s)=\frac{b(2L)f(sd(2L)/b(2L))}{d(2L)}=\frac{b(2L)}{d(2L)+b(2L)}+\frac{d(2L)}{d(2L)+b(2L)}s^2.
\]
Thus, in the conditioned process, there is a single particle at $L^A$ at the beginning. Each individual moves as standard Brownian motion. It lives for an exponentially distributed time with parameter $b(2L)+d(2L)$, and then gives birth to $0$ offspring with probability $b(2L)/(b(2L)+d(2L))$ and $2$ offspring with probability $d(2L)/(b(2L)+d(2L))$. Consider a critical branching process started from a single particle at $L^A$. Each individual moves as standard Brownian motion. Each particle lives for an exponentially distributed time with parameter $b(2L)+d(2L)$, and then gives birth to $0$ offspring with probability $1/2$ and $2$ offspring with probability $1/2$. We observe that the right-most position that is ever reached by particles up to time $t$ in the conditioned process is stochastically dominated by the right-most position that is ever reached by particles up to time ~$t$ in the critical process. Letting $M$ be the all-time maximal displacement of the critical process, we have for all $C_2>0$ and all time $t$,
\begin{equation}\label{domcrit}
P\bigg(M_t^*>L^A+\frac{C'_2}{\rho}\bigg|B^c\bigg)\leq P\bigg(M>L^A+\frac{C'_2}{\rho}\bigg)
\end{equation}
According to equation (1.7) of Sawyer and Fleischman \cite{FS79}, we have for $n$ large enough,
\begin{equation}\label{fsmax}
P\bigg(M>L^A+\frac{C'_2}{\rho}\bigg)\leq \frac{6}{(C'_2)^2}\rho^2.
\end{equation}
Letting $C_{12}=C_{13}+6/(C'_2)^2$, equations (\ref{domcond})-(\ref{fsmax}) imply (\ref{domprob}).
\qedwhite
\\

\noindent\textit{Proof of Proposition \ref{Corollaryr}.}
Let us first consider the case when
\[
\beta^{-2/3}\log^{1/3}\bigg(\frac{\rho}{\beta^{1/3}}\bigg)\ll t\ll\frac{\rho}{\beta}.
\]
We start with the proof of equation (\ref{rlower}), which follows directly from results in \cite{RS2020}. For any constant $C_1>0$, define
\begin{equation*}
f(x) = \begin{cases}
             1  & x<2^{1/3}C_1 \\
             0  & \text{otherwise}.
       \end{cases}
\end{equation*}
Define $\Phi(f)$ as in (\ref{Phi}). By \ref{RS6}, we see that $\Phi(f)$ satisfies (\ref{5.8}). For all $C_1>0$ and $0<\kappa<1$, if $n$ is sufficiently large, then
\begin{align}\label{rlower1}
P\bigg(M_t\geq L-\frac{C_1}{\beta^{1/3}}\bigg)
&=P\bigg(\Phi(f)\geq e^{\rho L -C_1\rho/\beta^{1/3}}\bigg)\nonumber\\
&\geq P\bigg(\frac{1-\kappa}{Ai'(\gamma_1)^2}\bigg(\int_{0}^{\infty}f(z)Ai(\gamma_1+z)dz\bigg)Z(0)\geq e^{\rho L-C_1\rho/\beta^{1/3}}\bigg)-\frac{\kappa}{2}.
\end{align}
By (\ref{assumpZ}), we have for $n$ sufficiently large
\begin{equation}\label{rlower2}
P\bigg(\frac{1-\kappa}{Ai'(\gamma_1)^2}\bigg(\int_{0}^{\infty}f(z)Ai(\gamma_1+z)dz\bigg)Z(0)\geq e^{\rho L-C_1\rho/\beta^{1/3}}\bigg)
\geq 1-\frac{\kappa}{2}.
\end{equation}
Equation (\ref{rlower}) follows from (\ref{rlower1}) and (\ref{rlower2}).

We next prove equation (\ref{rupper}) under the additional assumption that the birth rate function $b(x)$ is non-decreasing and the death rate function $d(x)$ is non-increasing. For any constant $C_2$, define $A$ and $C'_2$ as in (\ref{AC'}). We divide particles at time $t$ into the following categories:
\begin{align*}
S_1&=\big\{i\in\mathcal{N}_t:X_i(v)< L^A\;\text{for all}\;v\in [0,t]\big\},\\
S_2&=\big\{i\in\mathcal{N}_t\setminus S_1:X_i(v)\geq L^A\;\text{for some}\;v\in[0,t-C_7\rho^{-2}]\big\},\\
S_3&=\mathcal{N}_t\setminus (S_1\cup S_2).
\end{align*}
For $j=1,2,3$, write
\[
M^{S_j}_t=\max\big\{X_{i}(t), i\in S_j\big\}.
\]
Note that 
\begin{equation}\label{maxS123}
M_t=\max\Big\{M_t^{S_1},M_t^{S_2},M_t^{S_3}\Big\}.
\end{equation}

For $S_1$, it is obvious that for all constants $C_2$,
\begin{equation}\label{maxS1}
P\bigg(M_t^{S_1}\leq L+\frac{C_2}{\rho}\bigg)=P\bigg(M_t^{S_1}\leq L^A+\frac{C'_2}{\rho}\bigg)=1.
\end{equation}

For $S_2$, according to \ref{RS3}, with probability tending to $1$, particles that either are to the right of $L^A$ at time $0$ or hit $L^A$ before time $t-C_7\rho^{-2}$ will not have descendants alive at time $t$. Thus for all constant $C_2$,
\begin{equation}\label{maxS2}
\lim_{n\rightarrow\infty}P\bigg(M_t^{S_2}\leq L+\frac{C_2}{\rho}\bigg)\geq\lim_{n\rightarrow\infty}P(S_2=0)=1.
\end{equation}

It remains to deal with $S_3$, which consists of particles whose trajectories cross $L^A$ in the last $C_7\rho^{-2}$ units of time. Consider the process in which particles are killed upon hitting ~$L^A$. Let $R$ be the number of particles that first hit $L^A$ between $t-C_7\rho^{-2}$ and $t$. We denote by $\{r_i\}_{i=1}^{R}$ the sequence of hitting times. For the process started from a single particle at $x$, recall that $r_x^{L^A}(u,t)$ is the expected number of particles hitting $L^A$ between time $u$ and $t$. By (\ref{lem2.13}), taking all the particles at time $0$ into consideration, we have
\[
E[R|\mathcal{F}_0]=\sum_{i\in\mathcal{N}_{0}}r_{X_i(0)}^{L^A}\bigg(t-\frac{C_7}{\rho^2},t\bigg)\lesssim \exp\bigg(-\rho L^A-\frac{\beta^2(t-C_7\rho^{-2})^3}{9}\bigg)Y(0)+\frac{C_7}{\rho^2}\beta^{2/3}e^{-\rho L^A}Z(0).
\]
From (\ref{assumpY}) and (\ref{assumpZ}), it now follows that for any $\kappa>0$, there exists a constant $C_{14}$ such that for $n$ sufficiently large,
\[
P\bigg(E[R|\mathcal{F}_0]<C_{14}\frac{\beta}{\rho^{3}}\frac{1}{\rho^2}\bigg)>1-\frac{\kappa}{4}.
\]
Thus by the conditional Markov's inequality, we have for $n$ sufficiently large,
\begin{equation}\label{rightnumhit}
P\bigg(R>C_{14}\frac{\beta^{2/3}}{\rho^{2}}\frac{1}{\rho^2}\bigg)\leq E\bigg[\frac{E[R|\mathcal{F}_0]}{C_{14}\beta^{2/3}\rho^{-4}}1_{\{E[R|\mathcal{F}_0]<C_{14}\beta/\rho^5\}}\bigg]+P\bigg(E[R|\mathcal{F}_0]\geq C_{14}\frac{\beta}{\rho^{3}}\frac{1}{\rho^2}\bigg)<\frac{\kappa}{2}.
\end{equation}
Therefore, with probability at least $1-\kappa/2$, the number of particles that hit $L^A$ in the last $C_7\rho^{-2}$ unit of time is at most $C_{14}\beta^{2/3}/\rho^4$, which is $o(\rho^{-2})$. 

For every $i=1,...,R$, we consider three branching Brownian motions. All three processes start from a single particle at $L^A$ at time $r_i$. The first process has inhomogeneous birth rate $b(x)$ and death rate $d(x)$. Each particle moves as Brownian motion with drift $-\rho$. The second process is constructed based on the first process with the extra restriction that particles are killed upon hitting $2L$. To be more precise, in the second process, particles give birth at rate $b(x)$ and die at rate $d(x)$. Particles move as Brownian motion with drift $-\rho$ and are absorbed at $2L$. In the third process, the birth rate is the constant $b(2L)$ and the death rate is the constant $d(2L)$. Each particle moves as standard Brownian motion. We denote by $\bar{M}_{t-r_i}$, $\bar{M}^{2L}_{t-r_i}$ and $M^*_{t-r_i}$ the maximal positions that are ever reached by particles before time $t$ in the three processes respectively. Because of the monotonicity of $b(x)$ and $d(x)$, we observe that $M^*_{t-r_i}$ stochastically dominates $\bar{M}^{2L}_{t-r_i}$. By Lemma \ref{Dominant}, we have for sufficiently large $n$,
\[
P\bigg(\bar{M}^{2L}_{t-r_i}>L^A+\frac{C'_2}{\rho}\bigg)\leq P\bigg(M^*_{t-r_i}>L^A+\frac{C'_2}{\rho}\bigg)\leq C_{12}\rho^2.
\]
Note that $L^A+C'_2/\rho<2L$ for $n$ sufficiently large. Thus, if $\bar{M}^{2L}_{t-r_i}\leq L^A+C'_2\rho^{-1}$, then the first process is identical to the second process up to time $t-r_i$. Therefore, for sufficiently large $n$,
\begin{equation}\label{rightcouple}
P\bigg(\bar{M}_{t-r_i}>L+\frac{C_2}{\rho}\bigg)=P\bigg(\bar{M}_{t-r_i}>L_A+\frac{C'_2}{\rho}\bigg)
=P\bigg(\bar{M}^{2L}_{t-r_i}>L^A+\frac{C'_2}{\rho}\bigg)\leq C_{12}\rho^2.
\end{equation}
Combining (\ref{rightnumhit}) with (\ref{rightcouple}), for any $\kappa>0$, we have for $n$ sufficiently large,
\begin{equation}\label{maxS3}
P\bigg(M_{t}^{S_3}>L+\frac{C_2}{\rho}\bigg)\leq \frac{\kappa}{2}+C_{14}\frac{\beta^{2/3}}{\rho^2}\frac{1}{\rho^2}\cdot C_{12}\rho^2<\frac{3\kappa}{4}.
\end{equation}
As a result, equation (\ref{rupper}) follows from (\ref{maxS123}), (\ref{maxS1}), (\ref{maxS2}) and (\ref{maxS3}).

Now let us consider the case when $t$ satisfies (\ref{corassumpr}). It suffices to show that for every subsequence $(n_j)_{j=1}^{\infty}$, there exists a sub-subsequence $(n_{j_k})_{k=1}^{\infty}$, such that
\begin{equation}\label{rlowersub}
\lim_{k\rightarrow\infty}P\bigg(M_{t_{n_{j_k}},n_{j_k}}\geq L_{n_{j_k}}-\frac{C_1}{\beta^{1/3}_{n_{j_k}}}\bigg)=1
\end{equation}
and under the additional assumption on the birth rate and the death rate,
\begin{equation}\label{ruppersub}
\lim_{k\rightarrow\infty}P\bigg(M_{t_{n_{j_k}},n_{j_k}}\leq L_{n_{j_k}}+\frac{C_2}{\rho_{n_{j_k}}}\bigg)=1.
\end{equation}
By (\ref{corassumpr}), given a subsequence $(n_j)_{j=1}^{\infty}$, we can choose a sub-subsequence  $(n_{j_k})_{k=1}^{\infty}$ for which
\[
 \lim_{k\rightarrow\infty}\frac{\beta_{n_{j_k}}t_{n_{j_k}}}{\rho_{n_{j_k}}}=\tau\in[0,\infty).
\]
If $\tau=0$, then according to the previous argument, (\ref{rlowersub}) and (\ref{ruppersub}) hold. If $\tau>0$, choose times $(v_{n_{j_k}})_{k=1}^{\infty}$ for which 
\[
\beta_{n_{j_k}}^{-2/3}\log^{1/3}\bigg(\frac{\rho_{n_{j_k}}}{\beta_{n_{j_k}}^{1/3}}\bigg)\ll v_{n_{j_k}}\ll \frac{\rho_{n_{j_k}}}{\beta_{n_{j_k}}}.
\] 
Let $r_{n_{j_k}}=t_{n_{j_k}}-v_{n_{j_k}}$. By Remark \ref{RS1}, assumptions (\ref{assumpY}) and (\ref{assumpZ}) hold with $Y_{n_{j_k}}(r_{n_{j_k}})$ and $Z_{n_{j_k}}(r_{n_{j_k}})$ in place of $Z(0)$ and $Y(0)$ respectively. Replacing $Y(0)$ and $Z(0)$ by $Y_{n_{j_k}}(r_{n_{j_k}})$ and $Z_{n_{j_k}}(r_{n_{j_k}})$, the previous argument also works. Therefore, equations (\ref{rlowersub}) and (\ref{ruppersub}) also hold in this case. Equation (\ref{corright}) follows from (\ref{rlower}) and (\ref{rupper}).
\qedwhite

\subsection{Proof of Proposition \ref{Corollaryl}}
In this subsection, we will prove Proposition \ref{Corollaryl}, which gives the position of the left-most particle of the process.
Denote 
\[
t_1=t-\frac{2C_{10}}{\rho^2},\qquad t_2=t+\frac{2C_{10}}{\rho^2},
\]
where $C_{10}$ is defined in \ref{RS7}. We have the following lemma which controls the number of particles below $\bar{L}$ at any time between $t_1$ and $t_2$.

\begin{Lemma}\label{Leftnumber}
Suppose
\begin{equation}\label{Leftnumberassumpt}
\frac{\rho^{2/3}}{\beta^{8/9}}\ll t-t(\bar{L})\ll \frac{\rho}{\beta}.
\end{equation}
For any $\eps>0$, if $n$ is sufficiently large, then
\begin{equation}\label{leftnumt1}
P\bigg(N_{t_1}\big((-\infty, \bar{L})\big)\leq \frac{1}{\rho^2}\frac{\beta^{3/4}}{\rho^{9/4}}\bigg)>1-\eps.
\end{equation}
Moreover, for any $\eps>0$, if $n$ is sufficiently large, then there exists an event $B\in \mathcal{F}_{t_1}$ satisfying
\begin{equation}\label{Aprob}
 P(B)>1-\eps
\end{equation}
such that
\begin{equation}\label{leftnumoccup}
E\bigg[\int_{t_1}^{t_2}N_r\big((-\infty, \bar{L})\big)dr\cdot 1_{B}\bigg]\leq \frac{1}{\rho^4}\frac{\beta^{3/4}}{\rho^{9/4}}.
\end{equation}
\end{Lemma}

\noindent\textit{Proof of Lemma \ref{Leftnumber}.}
 Define $s$ as in (\ref{defsu}) and 
\[
u=t_2-t(\bar{L})+s=t-t(\bar{L})+\frac{2C_{10}}{\rho^2}+\frac{C_5}{\beta^{2/3}}.
\]
Since $2C_{10}\rho^{-2}+C_5\beta^{-2/3}\ll\rho/\beta$, by (\ref{Leftnumberassumpt}), we have $\rho^{2/3}/\beta^{8/9}\ll u\ll\rho/\beta$. For any $r\in [t_1,t_2]$, define
\[
s_r=s-r+t_2.
\]
Note that 
\begin{equation}\label{lem2.4d2res}
s_r\asymp s,\qquad r-u=t(\bar{L})-s_r.
\end{equation}
For every $r\in [t_1,t_2]$, denote
\begin{align*}
S_1(r)&=\{i\in\mathcal{N}_r: \exists v\in [0,u], X_i(v)\geq L\},\\
S_2(r)&=\{i\in\mathcal{N}_r:X_i(u)\leq L-C_6\beta^{-1/3}, X_i(v)< L\;\text{for all}\; v\in [0,u]\},\\
S_3(r)&=\{i\in\mathcal{N}_r: L-C_6\beta^{-1/3}<X_i(u)< L, X_i(v)< L\;\text{for all}\; v\in[0,u]\}.
\end{align*}
For $j=1,2,3$, write
\[
\Sigma_j(r)=\sum_{i\in S_j(r)}1_{\{X_i(r)\leq \bar{L}\}}.
\]
Then
\begin{equation}\label{Sigma123}
N_{r}\big((-\infty,\bar{L})\big)=\Sigma_1(r)+\Sigma_2(r)+\Sigma_3(r).
\end{equation}

We first consider $\Sigma_1(r)$. Let $B_1$ be the event that particles that are to the right of $L$ before time $u$ have descendants alive at time $t_1$. Then $B_1\in\mathcal{F}_{t_1}$. Since $u+C_7\rho^{-2}\ll t_1$, by \ref{RS3}, we have that for any $\eta>0$, if $n$ is sufficiently large, then
\[
P(B_1)<\eta.
\]
Note that $\{\Sigma_1(r)\neq 0\}$ is a subset of $B_1$ for all $r\in [t_1,t_2]$. Therefore, we have for sufficiently large $n$, for all $r\in [t_1,t_2]$, 
\begin{equation}\label{Sigma1A1}
P(\Sigma_1(r)\neq 0)\leq P(B_1)<\eta.
\end{equation}

We now consider $\Sigma_2(r)$. Denote 
\[
c_{\bar{L}}=\sqrt{1-\frac{\bar{L}}{L^*}}.
\]
By (\ref{lem2.4d2res}), we see that for all $r\in [t_1,t_2]$, equation (\ref{lem2.4d2}) holds with $r-u$ in place of $t$, $\bar{L}$ in place of $z$, $c_{\bar{L}}$ in place of $c$ and $s_r$ in place of $s$. By (\ref{lem2.4d2}) and Tonelli's theorem, for $n$ large enough, we have for all $r\in [t_1,t_2]$,
\begin{align}\label{Sigma2firstmoment}
E[\Sigma_{2}(r)|\mathcal{F}_u]&\leq \frac{1}{\sqrt{2\pi (r-u)}}\exp\bigg(g(\bar{L})-\rho L^*+\frac{\rho^2s_r}{2}-\frac{\beta^2s_r^3}{6}\bigg)\int_{-\infty}^{\bar{L}}e^{-\rho(c_{\bar{L}}-1)(\bar{L}-y)}dy\nonumber\\
&\hspace{0.2in}\times(H_1+H_2+H_3),
\end{align}
where $H_1$, $H_2$ and $H_3$ are defined as in (\ref{Theta2first}) and (\ref{Theta3sumdef}) but with $s_r$ in place of $s_y$. Since $s_r\asymp s$, the upper bounds on $H_1$, $H_2$ and $H_3$ in (\ref{unif5.6}), (\ref{Theta3sum}), (\ref{Theta3dct}) also hold here.
For any $\eta>0$, since ~$r-u\sim t(\bar{L})\sim 3\rho/2\beta$ and $c_{\bar{L}}\sim 3/2$, we have for $n$ sufficiently large, 
\begin{equation}\label{sigmacontrol1}
\frac{1}{\sqrt{2\pi (r-u)}}\int_{-\infty}^{\bar{L}}e^{-\rho(c_{\bar{L}}-1)(\bar{L}-y)}dy
= \frac{1}{\sqrt{2\pi (r-u)}}\frac{1}{\rho(c_{\bar{L}}-1)}\leq \frac{2(1+\eta)}{\sqrt{3\pi}}\frac{\beta^{1/2}}{\rho^{3/2}}.
\end{equation}
Also, because $(1-x)^{3/2}\geq 1-3x/2$, we have
\begin{align}\label{sigmacontrol2}
g(\bar{L})&=\rho \left(L^*-\bar{L}\right)-\frac{2\sqrt{2\beta}}{3}\left(L^*-\bar{L}\right)^{3/2}\nonumber\\
&=\frac{9\rho^3}{8\beta}-\frac{2^{2/3}\gamma_1\rho}{\beta^{1/3}}-\frac{9\rho^3}{8\beta}\bigg(1-\frac{8\beta}{9\rho^2}\frac{2^{2/3}\gamma_1}{\beta^{1/3}}\bigg)^{3/2}\nonumber\\
&\leq \rho(2\beta)^{-1/3}\gamma_1.
\end{align}
We get from equations (\ref{unif5.6}), (\ref{Theta3sum}), (\ref{Theta3dct}), (\ref{Sigma2firstmoment})--(\ref{sigmacontrol2}) that for all $r\in [t_1,t_2]$, if $n$ is sufficiently large, then
\begin{equation}\label{Sigma2int}
E[\Sigma_{2}(r)|\mathcal{F}_0]<\frac{2(1+\eta)}{\sqrt{3\pi}}\frac{\beta^{1/2}}{\rho^{3/2}}e^{-\rho L} \bigg(2e^{-\beta^2u^3/74}Y(0)+\frac{C_8\eta(1+\eta)}{2^{4/3}}Z(0)\bigg).
\end{equation}
By (\ref{assumpY}) and (\ref{assumpZ}), for any $\eta>0$, there exists a $\delta>0$ such that for $n$ sufficiently large
\begin{equation}\label{Sigma2A2}
P\bigg(\Big\{Y(0)>\frac{1}{\rho^2}e^{\rho L}\Big\}\cup \Big\{Z(0)>\frac{1}{\delta}\frac{\beta^{1/3}}{\rho^3}e^{\rho L}\Big\}\bigg)<\eta.
\end{equation}
Define $B_2$ to be the union of the previous two events. We see that $B_2\in \mathcal{F}_0\subset \mathcal{F}_{t_1}$. Note that $e^{-\beta^2 u^3/74}\ll \beta^{1/3}/\rho$. From equation (\ref{Sigma2int}), we have for all $r\in [t_1,t_2]$, 
\begin{equation}\label{Sigma2}
E[\Sigma_2(r)1_{B_2^c}]\lesssim \frac{1}{\rho^2}\frac{\beta^{5/6}}{\rho^{5/2}}\ll\frac{1}{\rho^2}\frac{\beta^{3/4}}{\rho^{9/4}}.
\end{equation}
Specifically, for $r=t_1$, by (\ref{Sigma2A2}) and (\ref{Sigma2}), we have for $n$ large enough
\begin{equation}\label{Sigma2t1}
P\bigg(\Sigma_2(t_1)>\frac{1}{3\rho^2}\frac{\beta^{3/4}}{\rho^{9/4}}\bigg)<2\eta.
\end{equation}

Similarly for $\Sigma_3(r)$, by (\ref{lem2.4d2}), (\ref{sigmacontrol1}), (\ref{sigmacontrol2}) and Tonelli's theorem, we have for $n$ sufficiently large
\[
E[\Sigma_{3}(r)|\mathcal{F}_u]<\frac{2(1+\eta)}{\sqrt{3\pi}}\frac{\beta^{1/2}}{\rho^{3/2}}e^{-\rho L}  \exp\bigg(\frac{\rho^2s_r}{2}-\frac{\beta^2s_r^3}{6}\bigg)\sum_{j\in\mathcal{N}_u}e^{(\rho-\beta s_r)X_j(u)}1_{\{L-C_6\beta^{-1/3}<X_j(u)<L\}}.
\]
Because $s\leq s_r\leq s+4C_{10}/\rho^2$ for all $r\in [t_1,t_2]$, we have
\begin{equation}\label{Sigma3int}
E[\Sigma_{3}(r)|\mathcal{F}_u]\lesssim \frac{\beta^{1/2}}{\rho^{3/2}}e^{-\rho L}  \exp\bigg(\frac{\rho^2s}{2}-\frac{\beta^2s^3}{6}\bigg)\sum_{j\in\mathcal{N}_u}e^{(\rho-\beta s)X_j(u)}1_{\{L-C_6\beta^{-1/3}<X_j(u)<L\}}.
\end{equation}
By (\ref{assumpZ}) and (\ref{unifgammaz}), for any $\eps>0$, there exists a $\delta>0$ such that for $n$ sufficiently large,
\begin{equation}\label{Sigma3A3}
P\bigg(\exp\bigg(\frac{\rho^2s}{2}-\frac{\beta^2s^3}{6}\bigg)\sum_{j\in\mathcal{N}_u}e^{(\rho-\beta s)X_j(u)}1_{\{L-C_6\beta^{-1/3}<X_j(u)<L\}}>\frac{1}{\delta}\frac{1+2\eta}{Ai'(\gamma_1)^2}\frac{\beta^{1/3}}{\rho^3}e^{\rho L}\bigg)<\eta.
\end{equation}
Define $B_3$ to be the event in the previous equation. Then $B_3\in\mathcal{F}_u\subset\mathcal{F}_{t_1}$. From (\ref{Sigma3int}), we have for all $r\in [t_1,t_2]$,
\begin{equation}\label{Sigma3}
E[\Sigma_3(r)1_{B_3^c}]
\lesssim \frac{1}{\rho^2}\frac{\beta^{5/6}}{\rho^{5/2}}\ll\frac{1}{\rho^2}\frac{\beta^{3/4}}{\rho^{9/4}}.
\end{equation}
Specifically, for $r=t_1$, by (\ref{Sigma3A3}) and (\ref{Sigma3}), we have for $n$ sufficiently large,
\begin{equation}\label{Sigma3t1}
P\bigg(\Sigma_3(t_1)>\frac{1}{3\rho^2}\frac{\beta^{3/4}}{\rho^{9/4}}\bigg)<2\eta.
\end{equation}

As a result, for any $\eps>0$, by choosing $\eta$ appropriately, equation (\ref{leftnumt1}) follows from (\ref{Sigma123}), (\ref{Sigma1A1}), (\ref{Sigma2t1}) and (\ref{Sigma3t1}). Let 
\[
B=B_1^c\cap B_2^c\cap B_3^c.
\]
Then $B\in \mathcal{F}_{t_1}$, and by (\ref{Sigma1A1}), (\ref{Sigma2A2}) and (\ref{Sigma3A3}), for $n$ large enough,
\[
P(B)=1-P(B_1\cup B_2\cup B_3)>1-3\eta.
\]
From (\ref{Sigma123}), (\ref{Sigma1A1}), (\ref{Sigma2}) and (\ref{Sigma3}), we have for $n$ sufficiently large
\begin{align*}
&E\bigg[\int_{t_1}^{t_2}N_r\big((-\infty, \bar{L})\big)dr\cdot 1_B\bigg]\\
&\hspace{0.2in}\leq E\bigg[\int_{t_1}^{t_2}\Sigma_1(r)dr\cdot 1_{B_1^c}\bigg]+E\bigg[\int_{t_1}^{t_2}\Sigma_2(r)dr\cdot 1_{B_2^c}\bigg]+E\bigg[\int_{t_1}^{t_2}\Sigma_3(r)dr\cdot 1_{B_3^c}\bigg]\\
&\hspace{0.2in}\leq\frac{1}{\rho^4}\frac{\beta^{3/4}}{\rho^{9/4}}
\end{align*}
Letting $\eta=\eps/3$, equations (\ref{Aprob}) and (\ref{leftnumoccup}) follow.
\qedwhite
\\

\noindent\textit{Proof of Proposition \ref{Corollaryl}.}
Let us first consider the case when $\rho^{2/3}/\beta^{8/9}\ll t-t(\bar{L})\ll \rho/\beta$. We start with the proof of equation (\ref{lupper}). For any sequence $(d_n)_{n=1}^{\infty}$ satisfying $\beta^{-1/3}\ll d\ll \rho^2/\beta$,  we claim that 
\begin{equation}\label{leftclaim}
\lim_{n\rightarrow\infty}P\big(m_t\leq\bar{L}+d\big)=1.
\end{equation}
To prove the claim, we first note that $\bar{L}+d$ satisfies assumptions (\ref{zn}), (\ref{zn1}) and (\ref{zn2}) in Proposition ~\ref{Global}. Furthermore, according to the Taylor expansion $\sqrt{1-x}=1-x/2+O(x^2)$, we get
\begin{align}\label{cortdiff}
t(\bar{L})-t(\bar{L}+d)
&=\sqrt{\frac{2}{\beta}(L^*-\bar{L})}\bigg(1-\sqrt{1-\frac{d}{L^*-\bar{L}}}\bigg)\nonumber\\
&=\sqrt{\frac{2}{\beta}(L^*-\bar{L})}\bigg(1-1+\frac{d}{2(L^*-\bar{L})}+O\bigg(\frac{d^2}{(L^*-\bar{L})^2}\bigg)\bigg)\nonumber\\
&=\frac{d}{\sqrt{2\beta(L^*-\bar{L})}}+O\bigg(\frac{d^2}{\sqrt{\beta}(L^*-\bar{L})^{3/2}}\bigg).
\end{align}
Since $d\ll \rho^2/\beta$, we have 
\[
\frac{d^2}{\sqrt{\beta}(L^*-\bar{L})^{3/2}}\ll \frac{d}{\sqrt{2\beta(L^*-\bar{L})}}\ll\frac{\rho}{\beta}.
\] 
Therefore,
\[
\frac{\rho^{2/3}}{\beta^{8/9}}\ll t-t(\bar{L}+d)=t-t(\bar{L})+\Big(t(\bar{L})-t(\bar{L}+d)\Big)\ll\frac{\rho}{\beta},
\]
which is assumption (\ref{assumptn}). Since for $0<x<1$, $(1-x)^{3/2}=1-3x/2+O(x^2)$, we have for $n$ large enough,
\begin{align*}
g\bigg(L^\dagger+\frac{d}{2}\bigg)
&= \rho\bigg(L^*-L^\dagger-\frac{d}{2}\bigg)-\frac{2\sqrt{2\beta}}{3}(L^*-L^\dagger)^{3/2}\bigg(1-\frac{d}{2(L^*-L^\dagger)}\bigg)^{3/2}\\
&=\rho(L^*-L^\dagger)-\frac{\rho d}{2}-\frac{2\sqrt{2\beta}}{3}(L^*-L^\dagger)^{3/2}\bigg(1-\frac{3}{2}\frac{d}{2(L^*-L^\dagger)}+O\bigg(\frac{d^2}{(L^*-L^\dagger)^2}\bigg)\bigg)\\
&=-\frac{\rho d}{2}+\frac{3\rho d}{4}+O\bigg(\frac{\beta d^2}{\rho}\bigg)\\
&\geq \frac{\rho d}{8}.
\end{align*}
Thus for $n$ sufficiently large,
\begin{equation}\label{leftclaimint}
\int_{(-\infty,\bar{L}+d)}\frac{1}{\sqrt{2\pi t(y)}}e^{g(y)}dy\geq \int_{[L^\dagger+d/2,\bar{L}+d)}\frac{1}{\sqrt{2\pi t(y)}}e^{g(y)}dy\geq e^{\rho d/16}.
\end{equation}
By Proposition \ref{Global} and equations (\ref{assumpZ}) and (\ref{leftclaimint}), for any $\eps>0$, if $n$ is sufficiently large, then
\begin{equation}\label{leftclaimprop}
P\big(m_t> \bar{L}+d\big)= P\bigg(N_t\big((-\infty,\bar{L}+d]\big)=0\bigg)< \eps,
\end{equation}
which implies (\ref{leftclaim}). Now we use (\ref{leftclaim}) to prove (\ref{lupper}). Suppose (\ref{lupper}) is not true. Then there exists $\kappa>0$, such that for any constant $C_{15}>0$, we have for infinitely many $n$,
\[
P\bigg(m_t\leq \bar{L}+\frac{C_{15}}{\beta^{1/3}}\bigg)\leq1-\kappa.
\]
We can therefore choose a sequence of positive integers $(n_j)_{j=1}^{\infty}$ and another sequence of positive constants $(C_{15,j})_{j=1}^{\infty}$ satisfying 
\[
n_j\ll 1 \qquad\text{and}\qquad 1\ll C_{15,j}\ll\rho^2_{n_j}/\beta^{2/3}_{n_j}
\]
such that for all ~$j$,
\begin{equation}\label{contra}
P\bigg(m_{t_{n_j},n_j}\leq \bar{L}_{n_j}+\frac{C_{15,j}}{\beta^{1/3}_{n_j}}\bigg)\leq1-\kappa.
\end{equation}
Let $d_{n_j}=C_{15,j}\beta^{-1/3}_{n_j}$. Note that $\beta_{n_j}^{-1/3}\ll d_{n_j}\ll\rho^2_{n_j}/\beta_{n_j}$. Then (\ref{contra}) contradicts (\ref{leftclaim}), and (\ref{lupper}) follows.

We next prove equation (\ref{llower}) under the additional assumption that for all $n$, the birth rate function $b(x)$ is non-decreasing and the death rate function $d(x)$ is non-increasing. Define
\begin{align*}
S_1&=\big\{i\in\mathcal{N}_t:X_i(t_1)< \bar{L}\big\},\\
S_2&=\big\{i\in\mathcal{N}_t\setminus S_1: X_i(v)\geq \bar{L} \;\text{for all}\; v\in[t_1,t_2]\big\},\\
S_3&=\mathcal{N}_t\setminus (S_1\cup S_2).
\end{align*}
For $j=1,2,3$, we further denote 
\[
m_t^{S_i}=\min\{X_{i}(t),i\in S_i\}.
\]
We have
\begin{equation}\label{minS12}
m_t=\min\big\{m_t^{S_1},m_t^{S_2},m_t^{S_3}\big\}.
\end{equation}

For $S_1$, we will show that particles below $\bar{L}$ at time $t_1$ will not have descendants survive until time $t$. For $x<\bar{L}$, consider one process starting from a single particle at $x$ at time $t_1$, and another process starting from a single particle at $L^*$ at time $t_1$. Because of the monotonicity of the birth rate and the death rate, the probability that the first process will survive until time $t$ is dominated by the probability that the second process will survive until  time $t$, which is at most $\rho^2/\alpha$ by \ref{RS7}. Thus by  (\ref{leftnumt1}), for any $\eps>0$ and all positive constant $C_4$, if $n$ is sufficiently large,
\begin{align}\label{minS1}
P\bigg(m_t^{S_1}\leq\bar{L}-\frac{C_4}{\rho}\bigg)
\leq P(S_1\neq \emptyset)
\leq \frac{\rho^2}{\alpha}\frac{1}{\rho^2}\frac{\beta^{3/4}}{\rho^{9/4}}+\eps<2\eps.
\end{align}

For $S_2$, it is obvious that for all $C_4$,
\begin{equation}\label{minS2}
P\bigg(m_t^{S_2}\leq\bar{L}-\frac{C_4}{\rho}\bigg)=0.
\end{equation}

We next deal with $S_3$. Consider the process in which particles are killed when they hit $\bar{L}$ between times $t_1$ and $t$. Let $R$ be the number of particles that are killed at $\bar{L}$ between times $t_1$ and $t$. For $i\in 1,2,...,R$, let $r_i\in [t_1,t]$ be the time that this particle is killed at $\bar{L}$. For the $i$th particle that hits $\bar{L}$ between $t_1$ and $t$, consider a process without killing starts from this particle and let $K_i(v)$ be the number of descendants of this particle to the left of $\bar{L}$ at time $r_i+v$. Define
\[
K_i=\int_{0}^{2C_{10}/\rho^2}K_i(v)dv.
\]
Then
\begin{equation}\label{leftSisum}
\sum_{i=1}^{R}K_i\leq \int_{t_1}^{t_2}N_r\big((-\infty, \bar{L})\big)dr.
\end{equation}
For all $i=1,2,...,R$, by Tonelli's theorem and (\ref{density2.8}), and interchanging the roles of $y$ and $\bar{L}-y$ in the last step, we have for $n$ large enough,
\begin{align}\label{Silbint}
E[K_i]&=\int_{0}^{2C_{10}/\rho^2}E[K_i(v)]dv\nonumber\\
&=\int_{0}^{2C_{10}/\rho^2}\int_{-\infty}^{\bar{L}}p_v(\bar{L},y)dydv\nonumber\\
&=\int_{0}^{2C_{10}/\rho^2}\int_{-\infty}^{\bar{L}}\frac{1}{\sqrt{2\pi v}}\exp\bigg(\rho \bar{L}-\rho y-\frac{(\bar{L}-y)^2}{2v}-\frac{\rho^2v}{2}+\frac{\beta(\bar{L}+y)v}{2}+\frac{\beta^2v^3}{24}\bigg)dydv\nonumber\\
&\geq \int_{0}^{2C_{10}/\rho^2}\int_{0}^{\infty}\frac{1}{\sqrt{2\pi v}}\exp\bigg(\rho y-\frac{y^2}{2v}-C_{10}+\frac{\beta(2\bar{L}-y)v}{2}\bigg)dydv.
\end{align}
Since $\rho\gg C_{10}\beta/\rho^2$, for $n$ sufficiently large, we have for all $y>0$ and $v\leq2C_{10}/\rho^2$,
\begin{equation}\label{Silby}
\rho y+\frac{\beta(2\bar{L}-y)v}{2}\geq -\frac{5}{4}C_{10}+o(1).
\end{equation}
Thus by (\ref{Silbint}) and (\ref{Silby}), we get for $n$ sufficiently large,
\begin{equation}\label{Silb}
E[K_i]\geq \int_{0}^{2C_{10}/\rho^2}\frac{1}{2}e^{-9C_{10}/4+o(1)}dv\geq \frac{C_{10}}{\rho^2}e^{-5C_{10}/2}.
\end{equation}
Note that random variables $\{K_i\}_{i=1}^{R}$ are independently and identically distributed. Moreover, conditioned on $\mathcal{F}_{t_1}$, the random variables $R$ and $K_i$ are independent. Let $B\in\mathcal{F}_{t_1}$ be the event defined in Lemma \ref{Leftnumber}. By (\ref{leftSisum}), we have
\begin{equation}\label{leftSicond}
E[R|\mathcal{F}_{t_1}]E[K_i|\mathcal{F}_{t_1}]1_{B}=E\bigg[\sum_{i=1}^{R}K_i\bigg|\mathcal{F}_{t_1}\bigg]1_{B}\leq E\bigg[\int_{t_1}^{t_2}N_r\big((-\infty, \bar{L})\big)dr\cdot 1_{B}\bigg|\mathcal{F}_{t_1}\bigg].
\end{equation}
Note that $E[K_i|\mathcal{F}_{t_1}]=E[K_i]$. By (\ref{leftnumoccup}), (\ref{Silb}) and (\ref{leftSicond}), we have for $n$ sufficiently large,
\begin{equation}\label{r(t1t2)ub}
E[R1_{B}]\leq \frac{1}{C_{10}}e^{5C_{10}/2}\frac{1}{\rho^2}\frac{\beta^{3/4}}{\rho^{9/4}}.
\end{equation}

For every $i=1,...,R$, we consider three branching Brownian motions. All three processes start from a single particle at $\bar{L}$ at time $r_i$. The first process has inhomogeneous birth rate $b(x)$ and death rate $d(x)$. Each particle moves as Brownian motion with drift $-\rho$. The second process is constructed based on the first process with the extra restriction that particles are killed upon hitting $0$. In the third process, the birth rate is the constant $b(0)$ and the death rate is the constant $d(0)$. Each particle moves as standard Brownian motion without drift. We denote by $m_{t-r_i}$ and $m^{0}_{t-r_i}$ the minimal displacement at time $t$ in the first and second processes respectively. We let $\bar{M}^{0}_{t-r_i}$ be the maximal position that is ever reached by a particle before time $t$ in the second process. Furthermore, for the third process, we denote by $m$ the all-time minimum and $M$ the all-time maximum. Because of the monotonicity of $b(x)$ and $d(x)$, we can couple the second and third processes such that $M$ stochastically dominates $\bar{M}^{0}_{t-r_i}$. Taking the drift into consideration, we also have that $m^{0}_{t-r_i}+\rho(t-r_i)$ stochastically dominates $m$. Note that in the third process, since $b(0)=d(0)$, the branching is critical and the process dies out eventually. According to equation (1.7) of Sawyer and Fleischman \cite{FS79}, we have for $x$ large enough
\begin{equation}\label{fsmM}
P(m<\bar{L}-x)\leq \frac{6}{x^2},\qquad P(M>\bar{L}+x)\leq \frac{6}{x^2}.
\end{equation}
Thus, by the construction of the first and the second processes, we have
\begin{align*}
P\bigg(m_{t-r_i}<\bar{L}-\frac{C_4}{\rho}\bigg)
&\leq P\bigg(\bigg\{m_{t-r_i}<\bar{L}-\frac{C_4}{\rho}\bigg\}\cap\bigg\{\bar{M}^0_{t-r_i}<\bar{L}+\frac{C_4}{\rho}\bigg\}\bigg)+P\bigg(\bar{M}^0_{t-r_i}\geq\bar{L}+\frac{C_4}{\rho}\bigg)\nonumber\\
&= P\bigg(m^0_{t-r_i}<\bar{L}-\frac{C_4}{\rho}\bigg)+P\bigg(\bar{M}^0_{t-r_i}\geq \bar{L}+\frac{C_4}{\rho}\bigg).
\end{align*}
Note that $t-r_i\leq 2C_{10}/\rho^2$ for all $r_i$. By the stochastic dominance relations and equation (\ref{fsmM}), for $C_4>2C_{10}$, we have for $n$ large enough,
\begin{equation}\label{fs79apply}
P\bigg(m_{t-r_i}<\bar{L}-\frac{C_4}{\rho}\bigg)
\leq P\bigg(m<\bar{L}-\frac{C_4-2C_{10}}{\rho}\bigg)+P\bigg(M\geq \bar{L}+\frac{C_4}{\rho}\bigg)
\leq \frac{12}{(C_4-2C_{10})^2}\rho^2.
\end{equation}
From (\ref{Aprob}), (\ref{r(t1t2)ub}) and (\ref{fs79apply}), we can choose constant $C_{4}$ large enough such that for $n$ large enough,
\begin{equation}\label{minS3}
P\bigg(m_t^{S_3}\leq \bar{L}-\frac{C_4}{\rho}\bigg)\leq E\bigg[\sum_{i=1}^{R}P\bigg(m_{t-r_i}<\bar{L}-\frac{C_4}{\rho}\bigg)1_{B}\bigg]+\eps\leq \frac{48}{C_{10}(C_{4}-2C_{10})^2}e^{5C_{10}/2}\frac{\beta^{3/4}}{\rho^{9/4}}+\eps<2\eps.
\end{equation}
For any $\kappa>0$, by choosing $\eps$ appropriately, equation (\ref{llower}) follows from (\ref{minS12}), (\ref{minS1}), (\ref{minS2}) and ~(\ref{minS3}).

Next we consider the case when $\beta (t-t(\bar{L}))/\rho\rightarrow\tau\in (0,\infty)$ as $n\rightarrow\infty$. Choose time $v<t$ for which 
\[
\frac{\rho^{2/3}}{\beta^{8/9}}\ll v-t(\bar{L})\ll \frac{\rho}{\beta}.
\]
Let $r=t-v$. By Remark \ref{RS1}, the previous argument still holds with $Z(r)$ in place of $Z(0)$ and $Y(r)$ in place of $Y(0)$. As a result, equations (\ref{lupper}) and (\ref{llower}) follow in this case. 

Finally, when $t$ satisfies (\ref{corassumpl}), we prove equations (\ref{lupper}) and (\ref{llower}) by contradiction. Since the proofs of equations (\ref{lupper}) and (\ref{llower}) are essentially the same, we only prove equation (\ref{lupper}). Suppose equation (\ref{lupper}) does not hold true. Then there exists $\kappa>0$ such that for all positive constants $C_3$, we have for infinitely many $n$,
\[
P\bigg(m_{t_n,n}\leq \bar{L}_n+\frac{C_3}{\beta_n^{1/3}}\bigg)\leq 1-\kappa.
\]
We can therefore choose a sequence of positive integers $(n_j)_{j=1}^{\infty}$ and another sequence of positive constants $(C_{3,j})_{j=1}^{\infty}$, both of which tend to infinity as $j\rightarrow\infty$, such that 
\begin{equation}\label{contrasub}
P\bigg(m_{t_{n_j},n_j}\leq \bar{L}_{n_j}+\frac{C_{3,j}}{\beta_{n_j}^{1/3}}\bigg)\leq 1-\kappa.
\end{equation}
For every subsequence $(n_j)_{j=1}^{\infty}$,  there exists a sub-subsequence $(n_{j_k})_{k=1}^{\infty}$ such that 
\begin{equation}\label{contratn}
\lim_{k\rightarrow\infty}\frac{\beta_{n_{j_k}}\big(t_{n_{j_k}}-t_{n_{j_k}}\big(\bar{L}_{n_{j_k}}\big)\big)}{\rho_{n_{j_k}}}=\tau\in [0,\infty).
\end{equation}
Fix any positive constant $C_3$. Then $C_3<C_{3,j_k}$ for $k$ sufficiently large. From (\ref{contrasub}), we have for $k$ large enough,
\begin{equation}\label{contralupper}
P\bigg(m_{t_{n_{j_k}},n_{j_k}}\leq \bar{L}_{n_{j_k}}+\frac{C_{3}}{\beta_{n_{j_k}}^{1/3}}\bigg)\leq 1-\kappa.
\end{equation}
On the other hand, note that $(t_{n_{j_k}})_{k=1}^{\infty}$ satisfies the assumptions of one of the previous two cases by (\ref{contratn}). Therefore, equation (\ref{lupper}) holds with $n_{j_k}$ in place of $n$, which contradicts (\ref{contralupper}). As a result, equation (\ref{lupper}) holds true. Equation (\ref{corleft}) follows from (\ref{lupper}) and (\ref{llower}).
\qedwhite

\section{Proofs of lemmas}\label{Lemmasec}

In this section, we will prove all the lemmas except Lemma \ref{Secondmomentg}, whose proof is deferred until Section \ref{2momsec}.
\\

\noindent\textit{Proof of Lemma \ref{YrestrictionZ}.}
First, let us prove equation (\ref{L^*-y}). When $z\leq 0$, equation (\ref{L^*-y}) holds automatically. It remains to consider the case $z>0$. If $z\ll \sqrt{\rho/\beta}$, then
\[
l\lesssim \sqrt{\frac{\rho}{\beta}}\ll \frac{\rho^2}{\beta}\asymp L^*-z.
\] 
If $z\gtrsim \sqrt{\rho/\beta}$, then according to (\ref{tL-z}), (\ref{clb}) and (\ref{cc0})
\[
l\lesssim \frac{1}{c_0\rho}\ll \frac{c^2\rho^2}{\beta}\asymp L^*-z.
\]
For $z$ satisfying (\ref{zn}) and $y\in[z-l,z+l]$, we have $L^*-y\geq L^*-z-l\gg \beta^{-1/3}$, which proves equation (\ref{L^*-y}).

We next prove equation (\ref{y-Ldagger}). When $z\geq 0$, equation (\ref{y-Ldagger}) is obvious. It remains to consider the case $z<0$. If $-z\ll \sqrt{\rho/\beta}$, then
\[
l\lesssim  \sqrt{\frac{\rho}{\beta}}\ll \frac{\rho^2}{\beta}\asymp z-L^\dagger.
\]
If $\sqrt{\rho/\beta}\lesssim-z\ll \rho^2/\beta$, then according to (\ref{c0lb}),
\[
l\lesssim \frac{1}{|c_0|\rho} \lesssim \sqrt{\frac{\rho}{\beta}}\ll z-L^\dagger.
\]
If $-z\asymp \rho^2/\beta$, then
\[
l\lesssim\frac{1}{|c_0|\rho}\asymp\frac{1}{\rho}\ll \frac{1}{\beta^{1/3}}\ll z-L^\dagger.
\]
For $z$ satisfying (\ref{zn}) and $y\in[z-l,z+l]$, we have $y-L^\dagger\geq z-L^\dagger-l\gg\beta^{-1/3}$, which proves equation (\ref{y-Ldagger}).
\qedwhite
\\

\noindent\textit{Proof of Lemma \ref{Sy}.} Following similar calculations to those in (\ref{cortdiff}), with the help of the Taylor expansion $\sqrt{1-x}= 1-x/2+O(x^2)$, we have for all $y\in [z-l,z+l]$,
\begin{align*}
|t(y)-t(z)|
&=\bigg|\sqrt{\frac{2}{\beta}(L^*-z)}\bigg(1-\frac{y-z}{2(L^*-z)}+O\bigg(\Big(\frac{y-z}{L^*-z}\Big)^2\bigg)-1\bigg)\bigg|\\
&\leq \frac{l}{\sqrt{2\beta(L^*-z)}}+O\bigg(\frac{l^2}{\sqrt{\beta}(L^*-z)^{3/2}}\bigg).
\end{align*}
Expressing the above formula in terms of $c$ according to (\ref{defcc0}), we obtain for all $y\in [z,z+l]$,
\[
|t(y)-t(z)|\leq \frac{l}{c\rho}+O\bigg(\frac{\beta l^2}{c^3\rho^3}\bigg).
\]
If $|z|\gtrsim \sqrt{\rho/\beta}$, then $l\lesssim 1/(|c_0|\rho)$. By formulas (\ref{clb}) and (\ref{cc0}), we get
\[
\frac{l}{c\rho}\lesssim \frac{1}{c|c_0|\rho^2}\ll \beta^{-2/3},\qquad \frac{\beta l^2}{c^3\rho^3}\lesssim \frac{\beta}{c_0^2c^3\rho^5}\ll \beta^{-2/3}.
\]
If $|z|\ll \sqrt{\rho/\beta}$, then $l\lesssim \sqrt{\rho/\beta}$ and $c\asymp 1$. We get
\[
\frac{l}{c\rho}\lesssim \frac{1}{c\sqrt{\rho\beta}}\ll \beta^{-2/3}, \qquad \frac{\beta l^2}{c^3\rho^3}\lesssim \frac{1}{c^3\rho^2}\ll \beta^{-2/3}.
\]
Combining the above three formulas, equation (\ref{sy}) follows. Moreover, by (\ref{L^*-y}), we have $t(y)\gg \beta^{-2/3}$ for all $y\in [z-l,z+l]$ and equation (\ref{sy}) implies (\ref{tyasymtz}).
\qedwhite
\\

\noindent\textit{Proof of Lemma \ref{Gy}.}
Note that for $0<x<1$, we have $(1-x)^{3/2}=1-3x/2+3x^2/8+O(x^3)$. For all $y$,
\begin{align*}
g(y)-g(z)&=\rho (z-y)+\frac{2\sqrt{2\beta}}{3}(L^*-z)^{3/2}\bigg(1-\Big(1-\frac{y-z}{L^*-z}\Big)^{3/2}\bigg)\\
&= \rho (z-y)+\frac{2\sqrt{2\beta}}{3}(L^*-z)^{3/2}\bigg(1-1+\frac{3(y-z)}{2(L^*-z)}-\frac{3(y-z)^2}{8(L^*-z)^{2}}+O\Big(\frac{|y-z|^3}{(L^*-z)^3}\Big)\bigg)\\
&=\rho (z-y)-\sqrt{2\beta(L^*-z)}(z-y)-\sqrt{\frac{\beta}{8(L^*-z)}}(z-y)^2+O\bigg(\frac{\beta^{1/2}|z-y|^3}{(L^*-z)^{3/2}}\bigg).
\end{align*}
Because $L^*-z=c^2\rho^2/(2\beta)$, the above equation can be expressed in terms of $c$ as
\begin{equation}\label{diffg}
g(y)-g(z)=\rho(1-c)(z-y)-\frac{\beta}{2c\rho}(z-y)^2+O\bigg(\frac{\beta^2|z-y|^3}{c^3\rho^3}\bigg).
\end{equation}
 For all $y\in [z-l,z+l]$, we have
\begin{equation}\label{fingyz}
|g(y)-g(z)|\leq \rho l|1-c|+\frac{\beta l^2}{2c\rho}+O\bigg(\frac{\beta^2l^3}{c^3\rho^3}\bigg).
\end{equation}
If $|z|\gtrsim\sqrt{\rho/\beta}$, then according to (\ref{1-c=c_0}), (\ref{c0lb}) and (\ref{cc0}), we see that
\begin{equation}\label{fingyz2}
\frac{\beta l^2}{2c\rho}\lesssim\rho l|1-c|\lesssim \frac{|1-c|}{|c_0|}= \frac{1}{1+c} \asymp1,\qquad \frac{\beta^2 l^3}{c^3\rho^3}\lesssim\frac{\beta^2}{|cc_0|^3\rho^6}\lesssim \frac{\beta^{1/2}}{\rho^{3/2}} \ll1.
\end{equation}
If $|z|\ll \sqrt{\rho/\beta}$, then according to (\ref{1-c=c_0}) and (\ref{cc02}), we see that
\begin{equation}\label{fingyz3}
\rho l|1-c| \lesssim\frac{\rho^{3/2}}{\beta^{1/2}}|c_0|\ll1,\qquad\frac{\beta^2l^3}{c^3\rho^3}\ll\frac{\beta l^2}{2c\rho}\lesssim \frac{1}{c}\asymp 1.
\end{equation}
The lemma follows from (\ref{fingyz}), (\ref{fingyz2}), and (\ref{fingyz3}).
\qedwhite
\\

\noindent\textit{Proof of Lemma \ref{Lemma2.4equivg}.} We are going to express $p_t(x,y)$ in terms of $s_y$ and $w$. Writing $t=t(y)-s_y$ and using (\ref{density2.8}), we have
\begin{align}\label{(2.8)sy}
p_t(x,y)&=\frac{1}{\sqrt{2\pi t}}\exp\bigg(\rho (L^*-w)-\rho y-\frac{(L^*-y-w)^2}{2t(y)}\sum_{k=0}^{\infty}\bigg(\frac{s_y}{t(y)}\bigg)^k-\frac{\rho^2}{2}(t(y)-s_y)\nonumber\\
&\hspace{0.2in}+\frac{\beta(y+L^*-w)(t(y)-s_y)}{2}+\frac{\beta^2(t(y)-s_y)^3}{24}\bigg)\nonumber\\
&=\frac{1}{\sqrt{2\pi t}}\exp\bigg(\rho L^*-\rho y-\frac{(L^*-y)^2}{2t(y)}-\frac{\rho^2 t(y)}{2}+\frac{\beta(y+L^*)t(y)}{2}+\frac{\beta^2t(y)^3}{24}\bigg)\nonumber\\
&\hspace{0.2in}\times\exp\bigg(-\rho w-\frac{\beta wt(y)}{2}+\frac{\beta(L^*-y)s_y}{2}+\frac{\beta ws_y}{2}-\frac{\beta^2s_y^3}{24}+\frac{\beta^2t(y)s_y^2}{8}-\frac{\beta^2t(y)^2s_y}{8}\nonumber\\
&\hspace{0.2in}-\frac{(L^*-y)^2}{2t(y)}\bigg(\frac{s_y}{t(y)}+\frac{s_y^2}{t(y)^2}+\frac{s_y^3}{t(y)^3}\bigg)+\frac{w(L^*-y)}{t(y)}\bigg(1+\frac{s_y}{t(y)}\bigg)\nonumber\\
&\hspace{0.2in}-\frac{(L^*-y)^2}{2t(y)}\sum_{k=4}^{\infty}\bigg(\frac{s_y}{t(y)}\bigg)^k-\frac{w^2}{2t(y)}\sum_{k=0}^{\infty}\bigg(\frac{s_y}{t(y)}\bigg)^k+\frac{w(L^*-y)}{t(y)}\sum_{k=2}^{\infty}\bigg(\frac{s_y}{t(y)}\bigg)^k\bigg).
\end{align}
Using that $t(y)=\sqrt{2(L^*-y)/\beta}$ and $y+L^*=\rho^2/\beta-(L^*-y)$, we get
\begin{align}\label{(2.8)g(z)}
&\rho L^*-\rho y-\frac{(L^*-y)^2}{2t(y)}-\frac{\rho^2t(y)}{2}+\frac{\beta(y+L^*)t(y)}{2}+\frac{\beta^2t(y)^3}{24}\nonumber\\
&\hspace{0.1in}=\rho(L^*-y)-\frac{\beta^{1/2}(L^*-y)^{3/2}}{2\sqrt{2}}-\frac{\rho^2t(y)}{2}+\bigg(\frac{\rho^2t(y)}{2}-\frac{\beta^{1/2}(L^*-y)^{3/2}}{\sqrt{2}}\bigg)+\frac{\beta^{1/2}(L^*-y)^{3/2}}{6\sqrt{2}}\nonumber\\
&\hspace{0.1in}=\rho (L^*-y)-\frac{2\sqrt{2\beta}}{3}(L^*-y)^{3/2}\nonumber\\
&\hspace{0.1in}=g(y).
\end{align}
Also notice that
\begin{align}\label{threeinfsums}
&-\frac{(L^*-y)^2}{2t(y)}\sum_{k=4}^{\infty}\bigg(\frac{s_y}{t(y)}\bigg)^k-\frac{w^2}{2t(y)}\sum_{k=0}^{\infty}\bigg(\frac{s_y}{t(y)}\bigg)^k+\frac{w(L^*-y)}{t(y)}\sum_{k=2}^{\infty}\bigg(\frac{s_y}{t(y)}\bigg)^k\nonumber\\
&\hspace{0.1in}=-\frac{1}{2(t(y)-s_y)}\bigg((L^*-y)\Big(\frac{s_y}{t(y)}\Big)^2-w\bigg)^2\nonumber\\
&\hspace{0.1in}=-\frac{(\beta s_y^2-2w)^2}{8t}.
\end{align}
For all $w$, according to (\ref{(2.8)sy}), (\ref{(2.8)g(z)}) and (\ref{threeinfsums}), replacing $t(y)$ with $\sqrt{2(L^*-y)/\beta}$, we have
\begin{align}\label{(2.8)exact}
p_t(x,y)&=\frac{1}{\sqrt{2\pi t}}\exp\bigg(g(y)-\rho w-\frac{\beta wt(y)}{2}+\frac{\beta(L^*-y)s_y}{2}+\frac{\beta ws_y}{2}-\frac{\beta^2s_y^3}{24}\nonumber\\
&\hspace{0.2in}+\frac{\beta^2t(y)s_y^2}{8}-\frac{\beta^2t(y)^2s_y}{8}-\frac{(L^*-y)^2}{2t(y)}\bigg(\frac{s_y}{t(y)}+\frac{s_y^2}{t(y)^2}+\frac{s_y^3}{t(y)^3}\bigg)\nonumber\\
&\hspace{0.2in}+\frac{w(L^*-y)}{t(y)}\bigg(1+\frac{s_y}{t(y)}\bigg)-\frac{(\beta s_y^2-2w)^2}{8t}\bigg)\nonumber\\
&=\frac{1}{\sqrt{2\pi t}}\exp\bigg(g(y)-\rho w+\beta ws_y-\frac{\beta^2}{6}s_y^3-\frac{(\beta s_y^2-2w)^2}{8t}\bigg).
\end{align}
For all $w\in\mathbb{R}$, $s< t(z)$ and $y\in\mathbb{R}$, since $-(\beta s_y^2-2w)^2/8t\leq 0$, formula (\ref{lem2.4egallwsmally}) follows. Furthermore, if $s\asymp\beta^{-2/3}$, then by Lemma \ref{Sy}, we have $s_y\asymp\beta^{-2/3}$ for all $y\in [z-l,z+l]$.  Then for $|w|\lesssim \beta^{-1/3}$, we get
\[
\frac{(\beta s_y^2-2w)^2}{8t}=o(1).
\]
and therefore (\ref{(2.8)exact}) implies (\ref{lem2.4egsmallw}).
\qedwhite
\\

\noindent\textit{Proof of Lemma \ref{Unifgamma}.}
Equation (\ref{unifgammaz}) follows from Lemma 6.2 in \cite{RS2020} directly. We use a similar strategy as in the proof of Lemma 6.2 in \cite{RS2020} to prove (\ref{unifgamma}). For every $y\in [z-l,z+l]$, define 
\begin{equation*}
f_y(x) = \begin{cases}
              e^{2^{-1/3}\beta^{2/3} s_yx} & 0<x<2^{1/3}C_6 \\
             0  & \text{otherwise}
       \end{cases} \quad \text{and} \quad 
f(x) = \begin{cases}
             e^{2^{-1/3}C_5x}  & 0<x<2^{1/3}C_6 \\
             0  & \text{otherwise}.
       \end{cases}
\end{equation*}
We can express $\Gamma_y$ in terms of the function $f_y$ by writing
\begin{align}\label{Gammarep}
\Gamma_y&=\exp\bigg(2^{-1/3}\beta^{2/3}s_y\gamma_1-\frac{\beta^2s_y^3}{6}\bigg)\sum_{j\in\mathcal{N}_u}e^{\rho X_j(u)}e^{\beta s_y(L-X_j(u))}1_{\{L-C_6\beta^{-1/3}<X_j(u)<L\}}\nonumber\\
&=\exp\bigg(2^{-1/3}\beta^{2/3}s_y\gamma_1-\frac{\beta^2s_y^3}{6}\bigg)\sum_{j\in\mathcal{N}_u}e^{\rho X_j(u)}f_y\Big((2\beta)^{1/3}\big(L-X_j(u)\big)\Big).
\end{align}
According to Lemma \ref{Sy}, we see that $\beta^{2/3}s_y\rightarrow C_5$ uniformly for $y\in[z-l,z+l]$. Thus
\begin{equation}\label{s-syunif}
\exp\bigg(2^{-1/3}\beta^{2/3}s_y\gamma_1-\frac{\beta^2s_y^3}{6}\bigg)\rightarrow \exp\bigg(2^{-1/3}\gamma_1C_5-\frac{C_5^3}{6}\bigg), \quad \text{as}\;n\rightarrow\infty.
\end{equation}
Also, for every $\eta>0$, if $n$ is sufficiently large, then for all $x$,
\[
\sup_{y\in[z,z+l]}|f_y(x)-f(x)|<\eta f(x).
\]
Therefore, for every $\eta>0$, if $n$ is large enough, then for all $y\in [z-l,z+l]$,
\begin{align}\label{Phify-Phif}
\bigg|\sum_{j\in\mathcal{N}_u}e^{\rho X_j(u)}&f\Big((2\beta)^{1/3}\big(L-X_j(u)\big)\Big)- \sum_{j\in\mathcal{N}_u}e^{\rho X_j(u)}f_y\Big((2\beta)^{1/3}\big(L-X_j(u)\big)\Big)\bigg|\nonumber\\
&\leq \eta\sum_{j\in\mathcal{N}_u}e^{\rho X_j(u)}f\Big((2\beta)^{1/3}\big(L-X_j(u)\big)\Big).
\end{align}
Furthermore, since $u$ satisfies (\ref{assump5.8}), equation (\ref{5.8}) implies that
\begin{equation}\label{Xi4conv}
\frac{1}{Z(0)}\sum_{j\in\mathcal{N}_u}e^{\rho X_j(u)}f\Big((2\beta)^{1/3}\big(L-X_j(u)\big)\Big)\rightarrow_{p}\frac{1}{Ai'(\gamma_1)^2}\int_{0}^{2^{1/3}C_6}e^{2^{-1/3}C_5y}Ai(\gamma_1+y)dy.
\end{equation}
As a result, equation (\ref{unifgamma}) follows from (\ref{C1C2}), (\ref{Gammarep})--(\ref{Xi4conv}).
\qedwhite
\\

\noindent\textit{Proof of Lemma \ref{Integralofg}.} First, let us summarize properties of the function $g(y)$ that will be useful throughout the proof. For $y\in(-\infty,L^*)$, we have
\[
g'(y)=-\rho+\sqrt{2\beta(L^*-y)},\quad g''(y)=-\sqrt{\frac{\beta}{2(L^*-y)}}<0.
\]
The function $g(y)$ is increasing in the interval $(-\infty, 0)$ and decreasing in the interval $[0,L^*)$. The derivative of $g(y)$ is decreasing and $g(y)$ is a concave function. Thus $g(y)$ is bounded above by its first order Taylor approximation. We have for all $x_1,x_2\in (-\infty, L^*]$,
\begin{equation}\label{concave}
g(x_2)\leq g(x_1)+g'(x_1)(x_2-x_1).
\end{equation}

First consider the case $z\geq 0$. By (\ref{concave}) and the fact that the derivative of $g$ is decreasing, we have for all $y\in [z,z+d]$,
\[
g(y)\geq g(z)+g'(y)(y-z)\geq g(z)+g'(z+d)(y-z).
\]
Also noticing that $t(y)$ is a decreasing function of $y$, we have
\begin{align*}
\int_{z}^{z+d}\frac{1}{\sqrt{2\pi t(y)}}e^{g(y)}dy
&\geq \frac{1}{\sqrt{2\pi t(z)}}e^{g(z)}\int_{z}^{z+d} e^{g'(z+d)(y-z)}dy\\
&=\frac{1}{\sqrt{2\pi t(z)}}e^{g(z)}\frac{1}{|g'(z+d)|}(1-e^{dg'(z+d)}).
\end{align*}
According to the definitions of $c_0$ in (\ref{defc0}) and $c$ in (\ref{defcc0}), we get
\[
dg'(z+d)\leq dg'(z)=-\frac{C_{11}}{1+c}\leq -\frac{C_{11}}{2}.
\]
Therefore,
\begin{equation}\label{intofgrhs}
\int_{z}^{z+d}\frac{1}{\sqrt{2\pi t(y)}}e^{g(y)}dy\geq \frac{1}{\sqrt{2\pi t(z)}}e^{g(z)}\frac{1}{|g'(z+d)|}(1-e^{-C_{11}/2}).
\end{equation}
Moreover, since $t(z)=c\rho/\beta$ and $|g'(z+d)|\leq \rho$, we see that
\begin{equation}\label{intofgrhs2}
\int_{z}^{z+d}\frac{1}{\sqrt{2\pi t(y)}}e^{g(y)}dy\gtrsim \frac{\beta^{1/2}}{c^{1/2}\rho^{3/2}}e^{g(z)}.
\end{equation}

On the other hand, we will separate the integral on the left hand side of (\ref{lemintegralofg1}) into two parts and upper bound each of them. Define 
\begin{equation}\label{hdef}
h=\beta^{-1/6}\sqrt{L^*-z}=\frac{c\rho}{\sqrt{2}\beta^{2/3}}.
\end{equation} We claim that $d\ll h$ and $h\ll L^*-z$. Indeed, since $z\gtrsim\sqrt{\rho/\beta}$, equation (\ref{cc0}) gives
\[
c_0c\gtrsim \frac{\beta^{1/2}}{\rho^{3/2}}\gg \frac{\beta^{2/3}}{\rho^2},
\]
which implies 
\[
d\asymp \frac{1}{c_0\rho}\ll \frac{c\rho}{\beta^{2/3}}\asymp h.
\]
Furthermore, because $c\gg \beta^{1/3}/\rho$, we have 
\begin{equation}\label{hsmall}
h\asymp \frac{c\rho}{\beta^{2/3}}\ll \frac{c^2\rho^2}{\beta}\asymp L^*-z.
\end{equation}
We denote
\begin{equation}\label{intofgsmalllarge}
\int_{z+d}^{L^*}\frac{1}{\sqrt{2\pi t(y)}}e^{g(y)}dy=\int_{z+d}^{z+h}\frac{1}{\sqrt{2\pi t(y)}}e^{g(y)}dy+\int_{z+h}^{L^*}\frac{1}{\sqrt{2\pi t(y)}}e^{g(y)}dy=:K_1+K_2.
\end{equation}
We first consider $K_1$. By (\ref{concave}), we have $g(y)\leq g(z+d)+g'(z+d)(y-z-d)$ for $y\in [z+d,z+h]$. Hence,
\begin{align}\label{intofgsmally1}
K_1\leq \frac{1}{\sqrt{2\pi t(z+h)}}e^{g(z+d)}\int_{z+d}^{z+h}e^{g'(z+d)(y-z-d)}dy\leq \frac{1}{\sqrt{2\pi t(z+h)}}e^{g(z+d)} \frac{1}{|g'(z+d)|}.
\end{align}
Since $z\gtrsim \sqrt{\rho/\beta}$, by (\ref{1-c=c_0}) and (\ref{cc0}), we have 
\begin{equation}\label{rhod(1-c)}
\rho d(1-c)\gg \frac{\beta^2d^3}{c^3\rho^3}.
\end{equation}
According to equations (\ref{1-c=c_0}) and (\ref{diffg}), we have for $n$ sufficiently large,
\begin{equation}\label{g(z+d)}
g(z+d)= g(z)-\rho d(1-c)-\frac{\beta d^2}{2c\rho}+O\bigg(\frac{\beta^2d^3}{c^3\rho^3}\bigg)\leq g(z)-\frac{\rho d(1-c)}{2}=g(z)-\frac{C_{11}}{2(1+c)}\leq g(z)-\frac{C_{11}}{4}.
\end{equation}
Also by (\ref{hsmall}), we have for $n$ large
\[
t(z+h)=\sqrt{1-\frac{h}{L^*-z}}t(z)\geq \frac{1}{2}t(z).
\]
Combining the above two observations with (\ref{intofgsmally1}), we have for $n$ large enough,
\begin{equation}\label{intofgsmally}
K_1\leq e^{-C_{11}/4}\frac{1}{\sqrt{\pi t(z)}}e^{g(z)}\frac{1}{|g'(z+d)|}.
\end{equation}
We next consider $K_2$. Recalling that the function $g(y)$ is decreasing when $y\in [0,L^*]$ and $t(y)=\sqrt{2/\beta}\sqrt{L^*-y}$, we get
\begin{equation}\label{intofglargey1}
K_2\leq e^{g(z+h)}\int_{z}^{L^*}\frac{\beta^{1/4}}{2^{3/4}\sqrt{\pi}}(L^*-y)^{-1/4}dy=\frac{2^{5/4}\beta^{1/4}}{3\sqrt{\pi}}(L^*-z)^{3/4}e^{g(z+h)}.
\end{equation}
We are going to apply the same argument that led to (\ref{g(z+d)}). Because $z\gtrsim \sqrt{\rho/\beta}$, we have $\rho c_0\gg\beta^{2/3}/(c\rho)$ by (\ref{cc0}). Thus by (\ref{1-c=c_0}) and (\ref{hdef}),  we get
\[
\rho h(1-c)\asymp \rho hc_0\gg \frac{\beta^{2/3}}{c\rho}h\asymp \frac{\beta^2h^3}{c^3\rho^3}.
\]
According to (\ref{diffg}) and (\ref{hdef}), since $cc_0\gg\beta^{7/12}/\rho^{7/4}$ by (\ref{cc0}), we have for $n$ sufficiently large,
\[
g(z+h)\leq g(z)-\frac{\rho h(1-c)}{2}=g(z)-\frac{cc_0\rho^2}{2\sqrt{2}\beta^{2/3}(1+c)}\leq g(z)- \frac{\rho^{1/4}}{\beta^{1/12}}.
\]
Combining this result with (\ref{intofglargey1}), we get for $n$ large,
\begin{equation}\label{intofglargey2}
K_2 \leq\frac{2^{5/4}\beta^{1/4}}{3\sqrt{\pi}}(L^*-z)^{3/4}\exp\bigg(g(z)-\frac{\rho^{1/4}}{\beta^{1/12}}\bigg)\asymp \frac{c^{3/2}\rho^{3/2}}{\beta^{1/2}}\exp\bigg(g(z)-\frac{\rho^{1/4}}{\beta^{1/12}}\bigg).
\end{equation}
Furthermore, since $c\leq 1$ and $\rho^3\gg\beta$, we notice that
\begin{equation}\label{intofglargey3}
\frac{c^{3/2}\rho^{3/2}}{\beta^{1/2}}\exp\bigg(g(z)-\frac{\rho^{1/4}}{\beta^{1/12}}\bigg)\ll\frac{\beta^{1/2}}{c^{1/2}\rho^{3/2}}e^{g(z)},
\end{equation}
As a result, equations (\ref{intofgrhs2}), (\ref{intofglargey2}) and (\ref{intofglargey3}) imply
\begin{equation}\label{intofglargey}
K_2\ll \int_{z}^{z+d}\frac{1}{\sqrt{2\pi t(y)}}e^{g(y)}dy.
\end{equation}
For any $\eta>0$, choosing $C_{11}$ large enough such that $\eta(1-e^{-C_{11}/2})/\sqrt{2}>e^{-C_{11}/4}$, equation (\ref{lemintegralofg1}) follows from (\ref{intofgrhs}), (\ref{intofgsmalllarge}), (\ref{intofgsmally}) and (\ref{intofglargey}).

Next consider the case $z\leq 0$. The proof is similar to the previous case. By (\ref{concave}) and the fact that the derivative of $g$ is decreasing, we have for all $y\in [z-d,z]$,
\[
g(y)\geq g(z)-g'(y)(z-y)\geq g(z)-g'(z-d)(z-y).
\]
Also, note that $t(y)$ is a decreasing function of $y$. Thus,
\begin{align*}
\int_{z-d}^{z}\frac{1}{\sqrt{2\pi t(y)}}e^{g(y)}dy
&\geq \frac{1}{\sqrt{2\pi t(z-d)}}e^{g(z)}\int_{z-d}^{z} e^{-g'(z+d)(z-y)}dy\\
&=\frac{1}{\sqrt{2\pi t(z-d)}}e^{g(z)}\frac{1}{g'(z-d)}(1-e^{-dg'(z-d)}).
\end{align*}
According to the definitions of $c_0$ in (\ref{defc0}) and $c$ in (\ref{defcc0}), since $c\in [1,3/2)$, we get
\[
-dg'(z-d)\leq -dg'(z)=-\frac{C_{11}}{1+c}\leq -\frac{2C_{11}}{5}.
\]
Therefore,
\begin{equation}\label{intofgrhs2}
\int_{z-d}^{z}\frac{1}{\sqrt{2\pi t(y)}}e^{g(y)}dy\geq \frac{1}{\sqrt{2\pi t(z-d)}}e^{g(z)}\frac{1}{g'(z-d)}(1-e^{-2C_{11}/5}).
\end{equation}

On the other hand, since $g(y)\leq g(z-d)+g'(z-d)(y-z+d)$ by (\ref{concave}) and $t(y)$ is decreasing, we get
\begin{align}\label{intofgsmall2}
\int_{-\infty}^{z-d}\frac{1}{\sqrt{2\pi t(y)}}e^{g(y)}dy
&\leq \frac{1}{\sqrt{2\pi t(z-d)}}e^{g(z-d)}\int_{-\infty}^{z-d}e^{g'(z-d)(y-z+d)}dy\nonumber\\
&\leq \frac{1}{\sqrt{2\pi t(z-d)}}e^{g(z-d)}\frac{1}{g'(z-d)}.
\end{align}
We will apply the argument that led to (\ref{g(z+d)}) again. Note that under the current scenario $c_0<0$ and $1\leq c<3/2$. From (\ref{1-c=c_0}), (\ref{diffg}) and (\ref{rhod(1-c)}), we get for $n$ sufficiently large,
\[
g(z-d)\leq g(z)+\frac{\rho d(1-c)}{2}=g(z)-\frac{C_{11}}{2(1+c)}\leq g(z)-\frac{C_{11}}{5}.
\]
Combining the above formula with (\ref{intofgsmall2}), we get for $n$ large
\begin{equation}\label{intofglhs}
\int_{-\infty}^{z-d}\frac{1}{\sqrt{2\pi t(y)}}e^{g(y)}dy\leq \frac{1}{\sqrt{2\pi t(z-d)}}e^{g(z)}\frac{1}{g'(z-d)}e^{-C_{11}/5}.
\end{equation}
For any $\eta>0$, choosing $C_{11}$ large enough such that $\eta(1-e^{-2C_{11}/5})>e^{-C_{11}/5}$, equation (\ref{lemintegralofg2}) follows from (\ref{intofgrhs2}) and (\ref{intofglhs}). Therefore, for any $\eta>0$, we can choose $C_{11}$ large enough such that (\ref{C6eta})--(\ref{lemintegralofg2}) all hold and the lemma follows. 
\qedwhite
\\

\noindent\textit{Proof of Lemma \ref{Lemma2.4equivd}.}
First consider the case $z > 0$. For $x\leq L$ and $y\in [\zeta,\infty)$, from (\ref{density2.8}),
\begin{equation}\label{lem2.4d2.81}
\frac{p_t(x,y)}{p_t(x,\zeta)}=\exp\bigg(-(y-\zeta)\bigg(\rho-\frac{\beta t}{2}-\frac{2x-\zeta-y}{2t}\bigg)\bigg).
\end{equation}
Using that $x\leq L$ and $y\geq \zeta$, we have
\begin{align*}
\rho-\frac{\beta t}{2}-\frac{2x-\zeta-y}{2t}
&\geq \rho-\frac{\beta t(z)}{2}+\frac{\beta s}{2}-\frac{2(L-z-2d)}{2t}\\
&\geq \rho-\frac{\beta t(z)}{2}+\frac{\beta s}{2}-\frac{L^*-z}{t}+\frac{2^{-1/3}\beta^{-1/3}\gamma_1}{t}.
\end{align*}
Note that we can expand $1/t$ as a geometric sum and thus
\[
\frac{L^*-z}{t}=\frac{L^*-z}{t(z)}+\frac{(L^*-z)s}{t(z)^2}+\frac{(L^*-z)s^2}{t(z)^3}\sum_{k=0}^{\infty}\bigg(\frac{s}{t(z)}\bigg)^k=\frac{L^*-z}{t(z)}+\frac{(L^*-z)s}{t(z)^2}+\frac{(L^*-z)s^2}{t(z)^2t}.
\]
Recall from (\ref{tL-z}) that $t(z)=c\rho/\beta$ and $L^*-z=c^2\rho^2/(2\beta)$. Therefore, from the above two formulas, we get
\begin{align*}
\rho-\frac{\beta t}{2}-\frac{2x-\zeta-y}{2t}&\geq \rho(1-c)-\frac{\beta s^2}{2t}+\frac{2^{-1/3}\beta^{-1/3}\gamma_1}{t}.
\end{align*}
Since $z\gtrsim \sqrt{\rho/\beta}$, we have $\rho (1-c)\gg \beta^{2/3}/(c\rho)$ by (\ref{1-c=c_0}) and (\ref{cc0}). Thus,
\[
\bigg|-\frac{\beta s^2}{2t}+\frac{2^{-1/3}\beta^{-1/3}\gamma_1}{t}\bigg|\asymp \frac{\beta^{2/3}}{c\rho}\ll \rho(1-c).
\]
Therefore, for $n$ sufficiently large, we have for all $x\leq L$ and $y\in [\zeta,\infty)$, 
\begin{equation}\label{lem2.4egallwlargeylb}
\rho-\frac{\beta t}{2}-\frac{2x-\zeta-y}{2t}\geq \frac{\rho}{2}(1-c).
\end{equation}
Equation (\ref{lem2.4d1}) follows from (\ref{lem2.4d2.81}) and (\ref{lem2.4egallwlargeylb}).

Next consider the case $z < 0$. We are going to apply an argument that is similar to the proof of (\ref{lem2.4egallwsmally}). Writing $r = z - y$ and expressing $p_t(x,y)$ in terms of $s$, $w$, $r$ and $z$, we have
\begin{align}\label{(2.8)sr}
p_t(x,y)&=\frac{1}{\sqrt{2\pi t}}\exp\bigg(\rho (L^*-w)-\rho(z-r) -\frac{(L^*-z+r-w)^2}{2t(z)}\sum_{k=0}^{\infty}\bigg(\frac{s}{t(z)}\bigg)^k-\frac{\rho^2}{2}(t(z)-s)\nonumber\\
&\hspace{0.2in}+\frac{\beta(z-r+L^*-w)(t(z)-s)}{2}+\frac{\beta^2(t(z)-s)^3}{24}\bigg)\nonumber\\
&=\frac{1}{\sqrt{2\pi t}}\exp\bigg(\rho L^*-\rho z-\frac{(L^*-z)^2}{2t(z)}-\frac{\rho^2 t(z)}{2}+\frac{\beta(z+L^*)t(z)}{2}+\frac{\beta^2t(z)^3}{24}\bigg)\nonumber\\
&\hspace{0.2in}\times\exp\bigg(-\rho w+\rho r-\frac{\beta (r+w)t(z)}{2}+\frac{\beta(L^*-z)s}{2}+\frac{\beta (r+w)s}{2}-\frac{\beta^2s^3}{24}+\frac{\beta^2t(z)s^2}{8}\nonumber\\
&\hspace{0.2in}-\frac{\beta^2t(z)^2s}{8}-\frac{(L^*-z)^2}{2t(z)}\bigg(\frac{s}{t(z)}+\frac{s^2}{t(z)^2}+\frac{s^3}{t(z)^3}\bigg)-\frac{(r-w)(L^*-z)}{t(z)}\bigg(1+\frac{s}{t(z)}\bigg)\nonumber\\
&\hspace{0.2in}-\frac{(L^*-z)^2}{2t(z)}\sum_{k=4}^{\infty}\bigg(\frac{s}{t(z)}\bigg)^k-\frac{(r-w)^2}{2t(z)}\sum_{k=0}^{\infty}\bigg(\frac{s}{t(z)}\bigg)^k-\frac{(r-w)(L^*-z)}{t(z)}\sum_{k=2}^{\infty}\bigg(\frac{s}{t(z)}\bigg)^k\bigg).
\end{align}
By a computation similar to the one leading to (\ref{threeinfsums}), we get
\begin{align*}
&-\frac{(L^*-z)^2}{2t(z)}\sum_{k=4}^{\infty}\bigg(\frac{s}{t(z)}\bigg)^k-\frac{(r-w)^2}{2t(z)}\sum_{k=0}^{\infty}\bigg(\frac{s}{t(z)}\bigg)^k-\frac{(r-w)(L^*-z)}{t(z)}\sum_{k=2}^{\infty}\bigg(\frac{s}{t(z)}\bigg)^k\\
&\hspace{0.2in}=-\frac{(\beta s^2+2r-2w)^2}{8t}\\
&\hspace{0.2in}\leq 0.
\end{align*}
Combining the previous two formulas with (\ref{tL-z}) and (\ref{(2.8)g(z)}), we have
\begin{align*}
p_t(x,y)&\leq \frac{1}{\sqrt{2\pi t}}\exp\bigg(g(z)-\rho w+\rho r-\frac{\beta (r+w)t(z)}{2}+\frac{\beta(L^*-z)s}{2}+\frac{\beta (r+w)s}{2}-\frac{\beta^2s^3}{24}\\
&\hspace{0.2in}+\frac{\beta^2t(z)s^2}{8}-\frac{\beta^2t(z)^2s}{8}-\frac{(L^*-z)^2}{2t(z)}\bigg(\frac{s}{t(z)}+\frac{s^2}{t(z)^2}+\frac{s^3}{t(z)^3}\bigg)\\
&\hspace{0.2in}-\frac{(r-w)(L^*-z)}{t(z)}\bigg(1+\frac{s}{t(z)}\bigg)\bigg)\\
&= \frac{1}{\sqrt{2\pi t}}\exp\bigg(g(z)-(c-1)\rho r-\rho w +\beta sw-\frac{\beta^2s^3}{6}\bigg),
\end{align*}
which is (\ref{lem2.4d2}). Furthermore, let 
\[
c_{\zeta}=\sqrt{1-\frac{\zeta}{L^*}}.
\]
Note that that $t=t(\zeta)-s_{\zeta}$ and $s_{\zeta}\asymp\beta^{-2/3}$ by (\ref{sy}). Thus for $y\leq \zeta$, equation (\ref{lem2.4d2}) holds with $\zeta$ in place of $z$, $c_{\zeta}$ in place of $c$ and $s_{\zeta}$ in place of $s$, so we have
\begin{equation}\label{lem2.4d2var1}
p_{t}(x,y)\leq \frac{1}{\sqrt{2\pi t}}\exp\bigg(g(\zeta)-(c_{\zeta}-1)\rho (\zeta-y)-\rho w +\beta s_{\zeta}w-\frac{\beta^2s_{\zeta}^3}{6}\bigg).
\end{equation}
Since $\zeta\leq z<0$, we have $c_{\zeta}>c>1$. Therefore, equation (\ref{lem2.4d2var}) follows from (\ref{lem2.4d2var1}).
\qedwhite

\section{Second moment estimate}\label{2momsec}

\subsection{Proof of Lemma \ref{Secondmomentg}}
A key step in the proof of Lemma \ref{Secondmomentg} is the following second moment estimate, which will be proved in Section 5.2. Recall that $p_t^L(x,y)$ is the density of the process in which particles are killed at $L$.

\begin{Lemma}\label{Secondmoment}
For every $z$ satisfying (\ref{zn}), (\ref{zn1}) and (\ref{zn2}), let 
\[
0\leq s\lesssim \beta^{-2/3},\qquad t=t(z)-s.
\]
Suppose $0\leq L-x\lesssim \beta^{-1/3}$. Then
\begin{align}\label{secondmoment}
\int_{0}^{t}\int_{-\infty}^{L}p_u^L(x,r)\Big(p_{t-u}^{L}(r,z)\Big)^2drdu\lesssim \frac{\beta^{2/3}}{t(z)\rho^4}\exp\bigg(\rho x-2\rho z+\rho L-\frac{4\sqrt{2\beta}}{3}(L^*-z)^{3/2}\bigg).
\end{align}
\end{Lemma}
Note that equation (\ref{secondmoment}) means that the ratio between the left hand side and the right hand side is bounded above by a positive constant uniformly for all $n$ and all $z$ satisfying (\ref{zn}), (\ref{zn1}) and (\ref{zn2}).
\\

\noindent\textit{Proof of Lemma \ref{Secondmomentg}.} According to the standard second moment formula (see e.g. Theorem 2.2 in \cite{Sawyer1976}), we have
\begin{equation}\label{2ndmomentformulag}
E[N_t^L(\mathcal{I})^2]
\lesssim \int_{\mathcal{I}}p_t^L(x,y)dy+2\int_{0}^{t}\int_{-\infty}^{L}p_u^L(x,r)\bigg(\int_{\mathcal{I}}p_{t-u}^L(r,y)dy\bigg)^2drdu.
\end{equation}

Regarding the first term in (\ref{2ndmomentformulag}), we upper bound $p_t^L(x,y)$ by $p_t(x,y)$ and then apply (\ref{lem2.4egallwsmally}) to get
\[
\int_{\mathcal{I}}p_t^L(x,y)dy\leq \int_{\mathcal{I}}\frac{1}{\sqrt{2\pi t}}\exp\bigg(g(y)-\rho(L^*-x)+\beta(L^*-x)(t(y)-t)-\frac{\beta^2}{6}(t(y)-t)^3\bigg)dy.
\]
On account of (\ref{sy}), we observe that $0\leq t(y)-t \asymp \beta^{-2/3}$  for $y\in \mathcal{I}$. Also notice that $|L^*-x|\lesssim \beta^{-1/3}$. Therefore, we get
\begin{align}\label{2ndfirst1}
\int_{\mathcal{I}}p_t^L(x,y)dy
&\lesssim \int_{\mathcal{I}}\frac{1}{\sqrt{2\pi t(y)}}\exp\bigg(g(y)-\rho L^*+\rho x\bigg)dy\nonumber\\
&=\frac{\beta^{2/3}}{\rho^4}e^{\rho x+\rho L-2\rho L^*}\bigg(\int_{\mathcal{I}}\frac{1}{\sqrt{2\pi t(y)}}e^{g(y)}dy\bigg)^{2} \frac{\rho^4}{\beta^{2/3}}e^{-\rho L+\rho L^*}\bigg(\int_{\mathcal{I}}\frac{1}{\sqrt{2\pi t(y)}}e^{g(y)}dy\bigg)^{-1}.
\end{align}
For $z$ satisfying (\ref{zn}), we see that $g(z)\geq 0$ and $t(y)\leq 2t(z)$ for all $y\in \mathcal{I}$ when $n$ is large enough. Also note that $c< 3/2$, $\rho/\beta^{1/3}\gg 1$ and $\gamma_1<0$. By (\ref{tL-z}), we get for $n$ large,
\begin{align}\label{2ndfirst2}
\frac{\rho^4}{\beta^{2/3}}e^{-\rho L+\rho L^*}\bigg(\int_{\mathcal{I}}\frac{1}{\sqrt{2\pi t(y)}}e^{g(y)}dy\bigg)^{-1}
&\leq \frac{\rho^4}{\beta^{2/3}}e^{-\rho L+\rho L^*}\bigg(\frac{l}{2\sqrt{\pi t(z)}}\bigg)^{-1}\nonumber\\
&=\frac{2\sqrt{\pi}c^{1/2}\rho}{l}\bigg(\frac{\rho}{\beta^{1/3}}\bigg)^{7/2}\exp\bigg(\frac{\gamma_1\rho}{2^{1/3}\beta^{1/3}}\bigg)\nonumber\\
&\ll1.
\end{align}
By (\ref{2ndfirst1}) and (\ref{2ndfirst2}), we have
\begin{equation}\label{2ndgfirst}
\int_{\mathcal{I}}p_t^L(x,y)dy\ll \frac{\beta^{2/3}}{\rho^4}e^{\rho x+\rho L-2\rho L^*}\bigg(\int_{\mathcal{I}}\frac{1}{\sqrt{2\pi t_y}}e^{g(y)}dy\bigg)^{2}.
\end{equation}

Regarding the second part of (\ref{2ndmomentformulag}), by the Cauchy-Schwarz inequality and Tonelli's theorem, we have
\begin{equation}\label{Cauchy-Schwarzg}
\int_{0}^{t}\int_{-\infty}^{L}p_u^L(x,r)\bigg(\int_{\mathcal{I}}p_{t-u}^L(r,y)dy\bigg)^2drdu
\leq l\int_{\mathcal{I}}\int_{0}^{t}\int_{-\infty}^{L}p_u^L(x,r)\big(p_{t-u}^{L}(r,y)\big)^2drdudy.
\end{equation}
We want to apply Lemma \ref{Secondmoment} to upper bound the above expression. First, by Lemma \ref{YrestrictionZ}, we know that (\ref{zn}) holds with $y$ in place of $z$. Also, by Lemma \ref{Sy}, for all $y\in \mathcal{I}$,
\[
t(y)=t(z)\pm o(\beta^{-2/3}).
\]
As a result, the assumptions in Lemma \ref{Secondmoment} are satisfied, and we can apply Lemma \ref{Secondmoment} to get
\begin{align}\label{2ndmomentpartIIglemmaplugin}
&l\int_{\mathcal{I}}\int_{0}^{t}\int_{-\infty}^{L}p_u^L(x,r)\big(p_{t-u}^{L}(r,y)\big)^2drdudy\nonumber\\
&\hspace{0.2in}\lesssim\frac{l\beta^{2/3}}{\rho^4}\int_{\mathcal{I}}\frac{1}{2\pi t(y)}\exp\bigg(\rho x-2\rho y+\rho L-\frac{4\sqrt{2\beta}}{3}(L^*-y)^{3/2}\bigg)dy\nonumber\\
&\hspace{0.2in}=\frac{\beta^{2/3}}{\rho^4}e^{\rho x+\rho L-2\rho L^*}l\int_{\mathcal{I}}\frac{1}{2\pi t(y)}e^{2g(y)}dy.
\end{align}
According to Lemma \ref{Gy}, we have for all $y\in \mathcal{I}$,
\[
e^{g(z)}\asymp e^{g(y)}.
\]
From the previous equation and (\ref{tyasymtz}), we get 
\begin{equation}\label{2ndmomentpartIIginvCS}
l\int_{\mathcal{I}}\frac{1}{2\pi t(y)}e^{2g(y)}dy\asymp\bigg(\int_{\mathcal{I}}\frac{1}{\sqrt{2\pi t(y)}}e^{g(y)}dy\bigg)^2.
\end{equation}
By equations (\ref{Cauchy-Schwarzg}), (\ref{2ndmomentpartIIglemmaplugin}) and (\ref{2ndmomentpartIIginvCS}), we obtain
\begin{equation}\label{2ndmomentpartIIg}
\int_{0}^{t}\int_{-\infty}^{L}p_u^L(x,r)\bigg(\int_{\mathcal{I}}p_{t-u}^L(r,y)dy\bigg)^2drdu\lesssim \frac{\beta^{2/3}}{\rho^4}e^{\rho x+\rho L-2\rho L^*}\bigg(\int_{\mathcal{I}}\frac{1}{\sqrt{2\pi t(y)}}e^{g(y)}dy\bigg)^2.
\end{equation}
The lemma follows from (\ref{2ndmomentformulag}), (\ref{2ndgfirst}) and (\ref{2ndmomentpartIIg}).
\qedwhite

\subsection{Proof of Lemma \ref{Secondmoment}}
The proof of Lemma \ref{Secondmoment} will be divided into the following three lemmas.
\begin{Lemma}\label{Small2ndmoment}
For every $z$ satisfying (\ref{zn}), (\ref{zn1}) and (\ref{zn2}), let 
\[
0\leq s\lesssim \beta^{-2/3},\qquad t=t(z)-s, \qquad u_1=\beta^{-7/12}(L^*-z)^{1/4}.
\]
Suppose $0\leq L-x\lesssim \beta^{-1/3}$. Then
\begin{align}\label{secondmomentI}
I_1:&=\int_{0}^{u_1}\int_{-\infty}^{L}p_u^L(x,r)\Big(p_{t-u}^{L}(r,z)\Big)^2drdu\lesssim \frac{\beta^{2/3}}{t(z)\rho^4}\exp\bigg(\rho x-2\rho z+\rho L-\frac{4\sqrt{2\beta}}{3}(L^*-z)^{3/2}\bigg).
\end{align}
\end{Lemma}
\begin{Lemma}\label{Large2ndmoment}
For every $z$ satisfying 
\begin{equation}\label{znvar1}
L^*-z\gg \beta^{-1/3}\log^{4/3} (\rho/\beta^{1/3}), \qquad z-L^\dagger\gg \beta^{-1/3},
\end{equation}
let 
\[
0\leq s\lesssim \beta^{-2/3},\qquad t=t(z)-s,  \qquad u_1=\beta^{-7/12}(L^*-z)^{1/4}.
\]
Suppose $0\leq L-x\lesssim \beta^{-1/3}$. Then
\begin{align}\label{secondmomentII}
I_2:&=\int_{u_1}^{t}\int_{-\infty}^{L}p_u^L(x,r)\Big(p_{t-u}^{L}(r,z)\Big)^2drdu\lesssim \frac{\beta^{2/3}}{t(z)\rho^4}\exp\bigg(\rho x-2\rho z+\rho L-\frac{4\sqrt{2\beta}}{3}(L^*-z)^{3/2}\bigg).
\end{align}
\end{Lemma}
\begin{Lemma}\label{2ndmomentgap}
For every positive $z$ satisfying 
\begin{equation}\label{znvar2}
\beta^{-1/3}\ll L^*-z\lesssim \beta^{-1/3}\log^{4/3} (\rho/\beta^{1/3}),
\end{equation}
let 
\[
0\leq s\lesssim \beta^{-2/3},\qquad t=t(z)-s,  \qquad u_1=\beta^{-7/12}(L^*-z)^{1/4}.
\]
Suppose $0\leq L-x\lesssim \beta^{-1/3}$. Then (\ref{secondmomentII}) holds.
\end{Lemma}

\noindent\textit{Proof of Lemma \ref{Secondmoment}.} 
It suffices to show that for every subsequence $(n_j)_{j=1}^{\infty}$, there exists a sub-subsequence $(n_{j_k})_{k=1}^{\infty}$, such that
\begin{align}\label{2ndmomsub}
&\int_{0}^{t_{n_{j_k}}}\int_{-\infty}^{L_{n_{j_k}}}p_u^{L_{n_{j_k}}}(x_{n_{j_k}},r)\Big(p_{t_{n_{j_k}}-u_{n_{j_k}}}^{L_{n_{j_k}}}(r,z_{n_{j_k}})\Big)^2drdu\nonumber\\
&\hspace{0.2in}\lesssim \frac{\beta^{2/3}_{n_{j_k}}}{t_{n_{j_k}}(z_{n_{j_k}})\rho_{n_{j_k}}^4}\exp\bigg(\rho_{n_{j_k}} x_{n_{j_k}}-2\rho_{n_{j_k}} z_{n_{j_k}}+\rho_{n_{j_k}} L_{n_{j_k}}-\frac{4\sqrt{2\beta_{n_{j_k}}}}{3}(L_{n_{j_k}}^*-z_{n_{j_k}})^{3/2}\bigg).
\end{align}
Given a subsequence $(n_j)_{j=1}^{\infty}$, there exists a further subsequence $(n_{j_k})_{k=1}^{\infty}$ such that one of the following holds:
\begin{enumerate}
\item $L_{n_{j_k}}^*-z_{n_{j_k}}\gg \beta^{-1/3}_{n_{j_k}}\log^{4/3} (\rho_{n_{j_k}}/\beta_{n_{j_k}}^{1/3})$ and $z_{n_{j_k}}-L_{n_{j_k}}^\dagger\gg \beta^{-1/3}_{n_{j_k}}$.
\item $\beta^{-1/3}_{n_{j_k}}\ll L_{n_{j_k}}^*-z_{n_{j_k}}\lesssim \beta_{n_{j_k}}^{-1/3}\log^{4/3} (\rho_{n_{j_k}}/\beta_{n_{j_k}}^{1/3})$.
\end{enumerate}
In case 1, equation (\ref{2ndmomsub}) follows from Lemmas \ref{Small2ndmoment} and \ref{Large2ndmoment}. In case 2, equation (\ref{2ndmomsub}) follows from Lemmas \ref{Small2ndmoment} and \ref{2ndmomentgap}.
\qedwhite
\\

The second moment estimate relies on delicate estimates of the density. Different approximations to the density $p^L_t(x,y)$ were obtained in \cite{RS2020}. The following results come from Lemmas 2.6, 2.7 and 2.8 in \cite{RS2020}.
\begin{Lemma}\label{RS20202.678}
For all $t\geq 0$ and $x,y<L$, we have
\begin{equation}\label{lemma2.6and2.7}
p_t^L(x,y)\lesssim \min\bigg\{\frac{1}{t^{1/2}},\frac{(L-x)(L-y)}{t^{3/2}}\bigg\}\exp\bigg(\rho x-\rho y-\frac{(y-x)^2}{2t}-\frac{\rho^2 t}{2}+\beta Lt\bigg).
\end{equation}
Moreover, when $t\geq 2\beta^{-2/3}$, $0\leq L-x\lesssim\beta^{-1/3}$ and $y<L$, we have
\begin{align}\label{lemma2.8}
p_t^L(x,y)\lesssim& \frac{\beta^{1/3}(L-x)}{\sqrt{t}}\max\bigg\{1,\frac{1}{\beta^{1/3} t}\Big(L-y-\frac{\beta t^2}{2}\Big)\bigg\}\nonumber\\
&\times \exp\bigg(\rho x-\rho y-\frac{(y-x)^2}{2t}-\frac{\rho^2 t}{2}+\frac{\beta(x+y)t}{2}+\frac{\beta^2t^3}{24}+\frac{1}{2\beta^{1/3}t}\Big(L-y-\frac{\beta t^2}{2}\Big)\bigg).
\end{align}
\end{Lemma}
It remains now to prove Lemmas \ref{Small2ndmoment}, \ref{Large2ndmoment} and \ref{2ndmomentgap}. 
\\

\noindent\textit{Proof of Lemma \ref{Small2ndmoment}.} 
For all $u\in [0,u_1]$, we see that $t-u\geq t/2\gg \beta^{-2/3}$ for $n$ sufficiently large. We will bound $p_u^L(x,r)$ by equation (\ref{lemma2.6and2.7}) and $p^L_{t-u}(r,z)$ by equations (\ref{density2.8}) and (\ref{lemma2.8}). We have
\begin{align*}
I_1&\lesssim\int_{0}^{u_1}\int_{-\infty}^{L}\min\bigg\{\frac{1}{u^{1/2}}, \frac{(L-x)(L-r)}{u^{3/2}}\bigg\}\exp\bigg(\rho x-\rho r-\frac{(x-r)^2}{2u}-\frac{\rho^2u}{2}+\beta Lu\bigg)\\
&\hspace{0.2in}\times \frac{1}{t-u}\Bigg(1_{\{L-r>\beta^{-1/3}\}}+\beta^{2/3}(L-r)^2\bigg(\max\bigg\{1,\frac{1}{\beta^{1/3}(t-u)}\bigg(L-z-\frac{\beta(t-u)^2}{2}\bigg)\bigg\}\bigg)^2\\
&\hspace{0.4in}\times\exp\bigg(\frac{1}{\beta^{1/3}(t-u)}\Big(L-z-\frac{\beta(t-u)^2}{2}\Big)\bigg)1_{\{0\leq L-r\leq \beta^{-1/3}\}}\Bigg)\\
&\hspace{0.2in}\times\exp\bigg(2\rho r-2\rho z-\frac{(r-z)^2}{t-u}-\rho^2(t-u)+\beta(r+z)(t-u)+\frac{\beta^2(t-u)^3}{12}\bigg)drdu.
\end{align*}
Denote
\[
M(u,r,x)=\min\bigg\{\frac{1}{u^{1/2}}, \frac{(L-x)r}{u^{3/2}}\bigg\},
\]
and
\begin{align*}
N(u,r,z)&=1_{\{r>\beta^{-1/3}\}}+\beta^{2/3}r^2\bigg(\max\bigg\{1,\frac{1}{\beta^{1/3}(t-u)}\bigg(L-z-\frac{\beta(t-u)^2}{2}\bigg)\bigg\}\bigg)^2\\
&\hspace{0.6in}\times\exp\bigg(\frac{1}{\beta^{1/3}(t-u)}\bigg(L-z-\frac{\beta(t-u)^2}{2}\bigg)\bigg)1_{\{0\leq r\leq \beta^{-1/3}\}}.
\end{align*}
Interchanging the roles of $r$ and $L-r$, we have
\begin{align*}
I_1&\lesssim\int_{0}^{u_1}\int_{0}^{\infty}\frac{1}{t-u}M(u,r,x)N(u,r,z)\exp\bigg(\rho x-2\rho z+\rho L-\rho r-\frac{(x-L+r)^2}{2u}-\frac{(L-r-z)^2}{t-u}\nonumber\\
&\hspace{0.2in}+\frac{\rho^2u}{2}-\rho^2t+\beta Lu+\beta(L-r)(t-u)+\beta z(t-u)+\frac{\beta^2(t-u)^3}{12}\bigg)drdu.
\end{align*}
Since $t=t(z)-s$, $L=L^*-(2\beta)^{-1/3}\gamma_1$ and
\begin{equation}\label{1/(t-u)}
\frac{1}{t-u}=\frac{1}{t(z)}\sum_{k=0}^{\infty}\bigg(\frac{u+s}{t(z)}\bigg)^k,
\end{equation}
we can express $-(L-r-z)^2/(t-u)$ as
\begin{align}\label{geometric}
&-\frac{(L^*-z)^2}{t(z)}\sum_{k=0}^{\infty}\bigg(\frac{u+s}{t(z)}\bigg)^k-\frac{((2\beta)^{-1/3}\gamma_1+r)^2}{t-u}+\frac{2(L^*-z)r}{t-u}+\frac{2(L^*-z)(2\beta)^{-1/3}\gamma_1}{t-u}\nonumber\\
&\hspace{0.3in}\leq -\frac{(L^*-z)^2}{t(z)}-\frac{(L^*-z)^2}{t(z)^2}(u+s)-\frac{(L^*-z)^2}{t(z)^3}(u+s)^2-\frac{(L^*-z)^2}{t(z)^4}(u+s)^3+\frac{2(L^*-z)r}{t-u}\nonumber\\
&\hspace{0.5in}+\frac{2(L^*-z)(2\beta)^{-1/3}\gamma_1}{t-u}.
\end{align}
Rearranging all the terms, $I_1$ can be further bounded as follows:
\begin{align}\label{I1simplify}
I_1&\lesssim \exp\bigg(\rho x-2\rho z+\rho L-\rho^2t(z)+\beta L^* t(z) +\beta zt(z)+\frac{\beta^2t(z)^3}{12}-\frac{(L^*-z)^2}{t(z)}\bigg)\int_{0}^{u_1}\int_{0}^{\infty}\frac{1}{t-u}\nonumber\\
&\hspace{0.15in}\times M(u,r,x)N(u,r,z)\exp\bigg(-\rho r-\frac{((L-x)-r)^2}{2u}-\frac{(L^*-z)^2}{t(z)^2}(u+s)-\frac{(L^*-z)^2}{t(z)^3}(u+s)^2\nonumber\\
&\hspace{0.15in}-\frac{(L^*-z)^2}{t(z)^4}(u+s)^3+\frac{2(L^*-z)r}{t-u}+\frac{2(L^*-z)(2\beta)^{-1/3}\gamma_1}{t-u}+\frac{\rho^2u}{2}+\rho^2s-\beta(2\beta)^{-1/3}\gamma_1t\nonumber\\
&\hspace{0.15in}-\beta(t-u)r-\beta z(u+s)-\beta L^*s-\frac{\beta^2(u+s)^3}{12}-\frac{\beta^2t(z)^2(u+s)}{4}+\frac{\beta^2t(z)(u+s)^2}{4}\bigg)drdu.
\end{align}
Note that 
\begin{equation}\label{2c^3/3}
-\rho^2t(z)+\beta L^* t(z) +\beta zt(z)+\frac{\beta^2t(z)^3}{12}-\frac{(L^*-z)^2}{t(z)}=-\frac{4\sqrt{2\beta}}{3}(L^*-z)^{3/2}.
\end{equation}
Also,
\begin{align*}
\frac{2(L^*-z)(2\beta)^{-1/3}\gamma_1}{t-u}-\beta(2\beta)^{-1/3}\gamma_1t
&< \frac{2(L^*-z)(2\beta)^{-1/3}\gamma_1}{t(z)}-\beta(2\beta)^{-1/3}\gamma_1t(z)+\beta(2\beta)^{-1/3}\gamma_1s\\
&=\beta(2\beta)^{-1/3}\gamma_1 s\\
&<0,
\end{align*}
and
\[
\frac{\rho^2 u}{2}+\rho^2 s-\beta z(u+s)-\beta L^* s=\beta (L^*-z)(u+s).
\]
Combining the above four formulas, we get
\begin{align*}
I_1&\lesssim \frac{1}{t}\exp\bigg(\rho x-2\rho z+\rho L-\frac{4\sqrt{2\beta}}{3}(L^*-z)^{3/2}\bigg)\\
&\hspace{0.2in}\int_{0}^{u_1}\int_{0}^{\infty}M(u,r,x)N(u,r,z)\exp\bigg(-r\bigg(\rho+\beta(t-u)-\frac{2(L^*-z)}{t-u}\bigg)-\frac{((L-x)-r)^2}{2u}\\
&\hspace{0.2in}-\frac{(L^*-z)^2}{t(z)^2}(u+s)-\frac{(L^*-z)^2}{t(z)^3}(u+s)^2-\frac{(L^*-z)^2}{t(z)^4}(u+s)^3+\beta(L^*-z)(u+s)\\
&\hspace{0.2in}-\frac{\beta^2(u+s)^3}{12}-\frac{\beta^2t(z)^2(u+s)}{4}+\frac{\beta^2t(z)(u+s)^2}{4}\bigg)drdu.
\end{align*}
Observe that
\[
-\frac{(L^*-z)^2}{t(z)^2}(u+s)+\beta(L^*-z)(u+s)-\frac{\beta^2t(z)^2(u+s)}{4}=0,
\]
\[
-\frac{(L^*-z)^2}{t(z)^3}(u+s)^2+\frac{\beta^2t(z)(u+s)^2}{4}=0,
\]
\[
-\frac{(L^*-z)^2}{t(z)^4}(u+s)^3-\frac{\beta^2(u+s)^3}{12}\leq -\frac{\beta^2u^3}{3}.
\] 
Therefore,
\begin{align}\label{I1J1J2}
I_1&\lesssim \frac{1}{t}\exp\bigg(\rho x-2\rho z+\rho L-\frac{4\sqrt{2\beta}}{3}(L^*-z)^{3/2}\bigg)\int_{0}^{u_1}\int_{0}^{\infty}M(u,r,x)N(u,r,z)\nonumber\\
&\hspace{0.2in}\exp\bigg(-r\bigg(\rho+\beta(t-u)-\frac{2(L^*-z)}{t-u}\bigg)-\frac{((L-x)-r)^2}{2u}-\frac{\beta^2u^3}{3}\bigg)drdu.
\end{align}
Denote the double integral as $J$. Write $J=J_1+J_2$ where $J_1$ is the portion of the double integral for which $r>\beta^{-1/3}$ and $J_2$ is the portion of the double integral for which $0\leq r\leq \beta^{-1/3}$. We first estimate $J_1$. Since $u\leq u_1=\beta^{-7/12}(L^*-z)^{1/4}$ and $s\ll u_1$, we see that $\beta(u+s)\ll \rho$. Thus, according to (\ref{1/(t-u)}), for $n$ sufficiently large we get
\begin{align}\label{I1J1intermediate}
\rho+\beta(t-u)-\frac{2(L^*-z)}{t-u}
&= \rho+t(z)\beta-\beta(u+s)-\frac{2(L^*-z)}{t(z)}\sum_{k=0}^{\infty}\bigg(\frac{u+s}{t(z)}\bigg)^{k}\nonumber\\
&=\rho-\beta(u+s)-\frac{2(L^*-z)}{t(z)}\sum_{k=1}^{\infty}\bigg(\frac{u+s}{t(z)}\bigg)^{k}\nonumber\\
&=\rho-\beta(u+s)-\frac{\beta t(z)(u+s)}{t-u}\nonumber\\
&\geq \frac{\rho}{2}.
\end{align}
Therefore,
\begin{align*}
J_1&\lesssim \int_{0}^{u_1}\int_{\beta^{-1/3}}^{\infty}M(u,r,x)N(u,r,z)e^{-\rho r/2}dr du\\
&\leq \int_{0}^{\rho^{-2}}\frac{1}{u^{1/2}}\int_{\beta^{-1/3}}^{\infty}e^{-\rho r/2}drdu+\int_{\rho^{-2}}^{u_1}\int_{\beta^{-1/3}}^{\infty}\frac{(L-x)r}{u^{3/2}}e^{-\rho r/2}drdu\\
&\leq e^{-\rho/2\beta^{1/3}}\bigg(\frac{4}{\rho^2}+\frac{4(L-x)}{\beta^{1/3}}+\frac{8(L-x)}{\rho}\bigg).
\end{align*}
Since $L-x\lesssim \beta^{-1/3}$, for $n$ sufficiently large,
\[
e^{-\rho/2\beta^{1/3}}\bigg(\frac{4}{\rho^2}+\frac{4(L-x)}{\beta^{1/3}}+\frac{8(L-x)}{\rho}\bigg)\lesssim e^{-\rho/2\beta^{1/3}}\beta^{-2/3}=\frac{\beta^{2/3}}{\rho^4}\Big(\frac{\rho}{\beta^{1/3}}\Big)^4e^{-\rho/2\beta^{1/3}}\ll \frac{\beta^{2/3}}{\rho^4}.
\]
Combining the above two equations, we have
\begin{equation}\label{I1J1}
J_1\ll \frac{\beta^{2/3}}{\rho^4}.
\end{equation}
Next, we estimate $J_2$. Note that
\begin{align}\label{1/beta^{1/3}(t-u)}
\frac{1}{\beta^{1/3}(t-u)}\bigg(L-z-\frac{\beta(t-u)^2}{2}\bigg)&=\frac{1}{\beta^{1/3}(t-u)}\bigg(-(2\beta)^{-1/3}\gamma_1-\frac{\beta(u+s)^2}{2}+\beta t(z)(u+s)\bigg) \nonumber \\
&= \frac{1}{\beta^{1/3}(t-u)}\bigg(-\frac{\beta(u+s)^2}{2}+\beta t(z)(u+s)\bigg) + o(1).
\end{align}
We will expand $1/(t-u)$ as a geometric sum. Using that $t=t(z)-s$, $u\leq u_1=\beta^{-7/12}(L^*-z)^{1/4}$ and $s\lesssim \beta^{-2/3}$, we see that $u + s \ll t(z)$.  Using also (\ref{tL-z}) and (\ref{1/(t-u)}), we get
\begin{align*}
\frac{1}{\beta^{1/3}(t-u)}\bigg(L-z-\frac{\beta(t-u)^2}{2}\bigg)&\leq\beta^{2/3}(u+s)\sum_{k=0}^{\infty}\bigg(\frac{u+s}{t(z)}\bigg)^k+o(1)\\
&= \beta^{2/3}(u+s) \bigg(1 + O \Big( \frac{u+s}{t(z)} \Big) \bigg) + o(1) \\
&=\beta^{2/3}u+O(1).
\end{align*}
Equation (\ref{I1J1intermediate}) and the previous formula imply that
\begin{align*}
J_2&\lesssim\int_{0}^{u_1}\int_{0}^{\beta^{-1/3}}M(u,r,x)\beta^{2/3}r^2\Big(\max\big\{1,\beta^{2/3}u\big\}\Big)^2\\
&\hspace{0.2in}\times\exp\bigg(-\frac{\rho r}{2}-\frac{(L-x-r)^2}{2u}-\frac{\beta^2u^3}{3}+\beta^{2/3}u\bigg)drdu\\
&\lesssim \int_{0}^{\rho^{-2}}\frac{1}{u^{1/2}}\int_{0}^{\beta^{-1/3}}\beta^{2/3}r^2e^{-\rho r/2}drdu\\
&\hspace{0.2in}+\int_{\rho^{-2}}^{3\beta^{-2/3}}\frac{1}{u^{3/2}}\int_{0}^{\beta^{-1/3}}(L-x)\beta^{2/3}r^3\exp\bigg(-\frac{\rho r}{2}-\frac{(L-x-r)^2}{2u}\bigg)drdu\\
&\hspace{0.2in}+\int_{3\beta^{-2/3}}^{u_1}\frac{1}{u^{3/2}}\int_{0}^{\beta^{-1/3}}(L-x)\beta^{2/3}r^3\exp\bigg(-\frac{\rho r}{2}-\frac{(L-x-r)^2}{2u}\bigg)\\
&\hspace{0.4in}\times\exp\bigg(-\frac{(\beta^{2/3}u)^3}{3}+\beta^{2/3}u+2\log (\beta^{2/3} u)\bigg)drdu.
\end{align*}
When $x\geq 3$, we have $-x^3/3+x+2\log x\leq 0$. Therefore,
\begin{align}\label{J2}
J_2&\lesssim \int_{0}^{\rho^{-2}}\frac{1}{u^{1/2}}\int_{0}^{\beta^{-1/3}}\beta^{2/3}r^2e^{-\rho r/2}drdu\nonumber\\
&\hspace{0.2in}+\int_{\rho^{-2}}^{u_1}\frac{1}{u^{3/2}}\int_{0}^{\beta^{-1/3}}(L-x)\beta^{2/3}r^3\exp\bigg(-\frac{\rho r}{2}-\frac{(L-x-r)^2}{2u}\bigg)drdu\nonumber\\
&\lesssim \frac{\beta^{2/3}}{\rho^4}+(L-x)\beta^{2/3}\int_{\rho^{-2}}^{u_1}\frac{1}{u^{3/2}}\int_{0}^{\beta^{-1/3}}r^3\exp\bigg(-\frac{\rho r}{2}-\frac{(L-x-r)^2}{2u}\bigg)drdu.
\end{align}
Note that $u_1\lesssim \rho^{1/2}\beta^{-5/6}$. By applying exactly the same calculation as in equations (8.54), (8.55) and (8.56) of \cite{RS2020}, we see that
\begin{equation}\label{I1J2}
(L-x)\beta^{2/3}\int_{\rho^{-2}}^{u_1}\frac{1}{u^{3/2}}\int_{0}^{\beta^{-1/3}}r^3\exp\bigg(-\frac{\rho r}{2}-\frac{(L-x-r)^2}{2u}\bigg)drdu\lesssim \frac{\beta^{2/3}}{\rho^4}.
\end{equation}
Combining equations (\ref{I1J1J2}), (\ref{I1J1}), (\ref{J2}) and (\ref{I1J2}), we have
\begin{equation*}
I_1\lesssim\frac{\beta^{2/3}}{t\rho^4}\exp\bigg(\rho x-2\rho z+\rho L-\frac{4\sqrt{2\beta}}{3}(L^*-z)^{3/2}\bigg),
\end{equation*}
which implies (\ref{secondmomentI}).\qedwhite
\\

\noindent\textit{Proof of Lemma \ref{Large2ndmoment}.}
The restriction (\ref{znvar1}) is equivalent to
\begin{equation}\label{znvar1c}
c\gg \frac{\beta^{1/3}}{\rho}\log^{2/3}\bigg(\frac{\rho}{\beta^{1/3}}\bigg),\qquad \frac{3}{2}-c\gg \frac{\beta^{2/3}}{\rho^2}.
\end{equation}
In particular, $c\gg \beta^{1/3}\rho^{-1}$ and $0<c<3/2$.
For $u_1\leq u\leq t$, we will bound both $p_u^L(x,r)$ and $p_{t-u}^{L}(r,z)$ by equation (\ref{density2.8}). Similar to the calculation for $I_1$, interchanging the roles of $r$ and $L-r$, we have
\begin{align}\label{I2firstformula}
I_2&\lesssim \exp\bigg(\rho x-2\rho z+\rho L-\rho^2 t+\beta Lt+\beta zt\bigg)\nonumber\\
&\hspace{0.2in}\times\int_{u_1}^{t}u^{-1/2}(t-u)^{-1}\exp\bigg(\frac{\rho^2 u}{2}+\frac{\beta xu}{2}-\beta zu-\frac{\beta Lu}{2}+\frac{\beta^2 u^3}{24}+\frac{\beta^2(t-u)^3}{12}\bigg)\nonumber\\
&\hspace{0.5in}\int_{0}^{\infty}\exp\bigg(-r\bigg(\rho+\beta t-\frac{\beta u}{2}\bigg)-\frac{(L-x-r)^2}{2u}-\frac{(L-z-r)^2}{t-u}\bigg)drdu.
\end{align}
Using (\ref{tL-z}), to prove (\ref{secondmomentII}), it is equivalent to show that
\begin{equation}\label{secondmomentpartIIc}
I_2\lesssim \frac{\beta^{5/3}}{c\rho^5}\exp\bigg(\rho x-2\rho z+\rho L-\frac{2c^3}{3}\frac{\rho^3}{\beta}\bigg).
\end{equation}

To simplify $I_2$, we are going to estimate $K$ in two ways depending on the value of $u$. This cutoff value $u_2$ will be defined based on $y_c$ introduced below. Denote
\begin{equation}\label{Deltadef}
\Delta=\sqrt{\frac{9c^2}{4}+c+1},\qquad y_c=\frac{3}{2}+\frac{1}{c}-\frac{1}{c}\Delta.
\end{equation}
Note that $y_c<1$ since for $0<c\leq 3/2$,
\begin{align}\label{1-yc}
1-y_c
&=-\frac{1}{2}+\frac{(\Delta-1)(\Delta+1)}{c(\Delta+1)}\nonumber\\
&=\frac{9c}{4(\Delta+1)}-\bigg(\frac{1}{2}-\frac{1}{\Delta+1}\bigg)\nonumber\\
&=\frac{9c}{4(\Delta+1)}-\frac{(\Delta-1)(\Delta+1)}{2(\Delta+1)^2}\nonumber\\
&=\frac{9c}{4(\Delta+1)}-\frac{c}{2(\Delta+1)^2}-\frac{9c^2}{8(\Delta+1)^2}\nonumber\\
&\geq\frac{1}{(\Delta+1)^2}\bigg(\frac{9c}{2}-\frac{c}{2}-\frac{9c^2}{8}\bigg)\nonumber\\
&>0.
\end{align}
Choose a constant $C_{16}>0$ small enough such that the following hold: 
\begin{equation}\label{C81}
\frac{1}{3}-\frac{3C_{16}}{2}>C_{16},
\end{equation}
\begin{equation}\label{C82}
-\frac{4}{3(3+\sqrt{17})(1+\sqrt{17}/2)}+2C_{16}<-C_{16},
\end{equation}
\begin{equation}\label{C83}
-\frac{7}{2(3+\sqrt{17})(1+\sqrt{17}/2)}+\frac{3}{2}C_{16}+6C_{16}<0,
\end{equation}
\begin{equation}\label{C84}
\frac{-8}{9(3\cdot \frac{3}{2}+2\cdot\frac{11}{4})(\frac{11}{4}+1)}+9C_{16}<-C_{16}.
\end{equation}
Then define
\begin{equation}\label{u2}
u_2=\bigg(\frac{3}{2}+\frac{1}{c}-\frac{1}{c}\sqrt{\frac{9c^2}{4}+c+1}-C_{16} c\bigg)t.
\end{equation}
Since $y_c<1$ by (\ref{1-yc}), we see that $u_2<t$. Also, when $0<c\leq 3/2$, one can show that $3/2+1/c-\Delta/c-C_{16}c$ is a decreasing function of $c$. Thus by (\ref{C81}),
\begin{equation}\label{u2C8}
u_2 \geq \bigg(\frac{3}{2}+\frac{2}{3}-\frac{2}{3}\sqrt{\frac{9}{4}\cdot\frac{9}{4}+\frac{3}{2}+1}-\frac{3C_{16}}{2}\bigg)t=\bigg(\frac{1}{3}-\frac{3C_{16}}{2}\bigg)t>C_{16}t
\end{equation}
and $u_2>u_1$ for sufficiently large $n$. 

Denote the inner integral in (\ref{I2firstformula}) as $K$.  When $u_1<u\leq u_2$, letting $a=(\rho+\beta t-\beta u/2)(t-u)$ and $b=L-z$, $K$ can be written as
\begin{equation}\label{Kfirst}
K=\int_{0}^{\infty}\exp\bigg(-\frac{(r-(b-a/2))^2+ab-a^2/4}{t-u}-\frac{(L-x-r)^2}{2u}\bigg)dr.
\end{equation}
For $u_1\leq u\leq u_2$, we claim that $b-a/2\leq 0$. Because $b - a/2$ is an increasing function of $u$, it is sufficient to show that for $u=u_2$, we have $b-a/2\leq 0$. Recall from (\ref{tL-z}) that $t=c\rho/\beta-s$ and $L^*-z=c^2\rho^2/2\beta$. Writing 
\begin{equation}\label{y2}
y_2=\frac{u_2}{t}=y_c-C_{16}c,
\end{equation}
we have that
\begin{align}\label{b-a/2}
b-\frac{a}{2}
&=(1-c_0)\frac{\rho^2}{2\beta}-(2\beta)^{-1/3}\gamma_1-\frac{1}{2}\bigg((1+c)\rho-\beta s-\frac{cy_2\rho}{2}+\frac{\beta y_2s}{2}\bigg)(1-y_2)\bigg(c\frac{\rho}{\beta}-s\bigg)\nonumber\\
&=\frac{\rho^2}{\beta}\bigg(-\frac{c^2y_2^2}{4}+y_2\Big(\frac{3c^2}{4}+\frac{c}{2}\Big)-\frac{c}{2}\bigg)+\rho s(1-y_2)\bigg(\frac{1}{2}+c-\frac{cy_2}{2}\bigg)+O(\beta^{-1/3}).
\end{align}
We observe that $y_c$ is one root of
\begin{equation}\label{root}
-\frac{c^2y^2}{4}+y\Big(\frac{3c^2}{4}+\frac{c}{2}\Big)-\frac{c}{2}=0.
\end{equation}
Therefore, according to (\ref{y2}) and (\ref{root}), equation (\ref{b-a/2}) implies that
\[
b-\frac{a}{2}
=-\frac{c^2\rho^2}{\beta}\bigg(\frac{C_{16}^2c^2}{4}+\frac{C_{16}\Delta}{2}\bigg)+\rho s\bigg(1-y_c+C_{16}c\bigg)\bigg(\frac{1}{2}+c-\frac{cy_2}{2}\bigg)+O(\beta^{-1/3}).
\]
Because $1-y_c\leq 9c/4(\Delta+1)$ by the second line of (\ref{1-yc}) and $1/2+c-cy_2/2\leq 2$ for $0<c\leq 3/2$, it follows that
\[
b-\frac{a}{2}\leq -\frac{C_{16}\Delta c^2\rho^2}{2\beta}+2\rho sc\bigg(\frac{9}{4(\Delta+1)}+C_{16}\bigg)+O(\beta^{-1/3}).
\]
Since $c\gg \beta^{1/3}\rho^{-1}$ by (\ref{znvar1}) and $s\lesssim \beta^{-2/3}$, we see that $c^2\rho^2/\beta\gg \rho s c$ and $c^2\rho^2/\beta\gg \beta^{-1/3}$. As a result, for $u=u_2$, and thus for all $u_1<u\leq u_2$, for $n$ sufficiently large,
\begin{equation}\label{b-a/2<0}
b-\frac{a}{2}\leq 0.
\end{equation}
From (\ref{Kfirst}) and (\ref{b-a/2<0}), we obtain that for $u_1<u\leq u_2$,
\begin{equation}\label{Kfirstupperbound}
K\leq \int_{0}^{\infty}\exp\bigg(-\frac{(L-z)^2}{t-u}-\frac{(L-x-r)^2}{2u}\bigg)dr\leq \sqrt{2\pi u}\exp\bigg(-\frac{(L-z)^2}{t-u}\bigg).
\end{equation}
When $u_2\leq u\leq t$, we bound $K$ by the formula for the moment generating function of the normal distribution to get
\begin{align}\label{Ksecondupperbound}
K&\leq \int_{-\infty}^{\infty}\exp\bigg(-r\Big(\rho+\beta t-\frac{\beta u}{2}-\frac{L-x}{u}-\frac{2(L-z)}{t-u}\Big)-\frac{r^2}{2}\frac{t+u}{u(t-u)}-\frac{(L-x)^2}{2u}-\frac{(L-z)^2}{t-u}\bigg)dr\nonumber\\
&=\sqrt{\frac{2\pi u(t-u)}{t+u}}\exp\bigg(\frac{u(t-u)}{2(t+u)}\Big(\rho+\beta t-\frac{\beta u}{2}-\frac{L-x}{u}-\frac{2(L-z)}{t-u}\Big)^2-\frac{(L-x)^2}{2u}-\frac{(L-z)^2}{t-u}\bigg).
\end{align}
According to (\ref{I2firstformula}), (\ref{Kfirstupperbound}) and (\ref{Ksecondupperbound}), we have
\begin{align}\label{defR1R2}
I_2&\lesssim \exp\bigg(\rho x-2\rho z+\rho L-\rho^2 t+\beta Lt+\beta zt\bigg)\nonumber\\
&\hspace{0.4in}\times\int_{u_1}^{u_2}(t-u)^{-1}\exp\bigg(\frac{\rho^2 u}{2}+\frac{\beta xu}{2}-\beta zu-\frac{\beta Lu}{2}+\frac{\beta^2 u^3}{24}+\frac{\beta^2(t-u)^3}{12}-\frac{(L-z)^2}{t-u}\bigg)du\nonumber\\
&\hspace{0.2in}+\exp\bigg(\rho x-2\rho z+\rho L-\rho^2 t+\beta Lt+\beta zt\bigg)\nonumber\\
&\hspace{0.4in}\times\int_{u_2}^{t}\sqrt{\frac{1}{(t-u)(t+u)}}\exp\bigg(\frac{\rho^2 u}{2}+\frac{\beta xu}{2}-\beta zu-\frac{\beta Lu}{2}+\frac{\beta^2 u^3}{24}+\frac{\beta^2(t-u)^3}{12}\bigg)\nonumber\\
&\hspace{0.6in}\times\exp\bigg(\frac{u(t-u)}{2(t+u)}\Big(\rho+\beta t-\frac{\beta u}{2}-\frac{L-x}{u}-\frac{2(L-z)}{t-u}\Big)^2-\frac{(L-x)^2}{2u}-\frac{(L-z)^2}{t-u}\bigg)du\nonumber\\
&=:R_1+R_2.
\end{align}

We first estimate $R_1$. Let $y_1=u_1/t$ and $y_2=u_2/t$. After making the change of variables $u=yt$ and writing $z=c_0 L^*$ and $t=c\rho/\beta-s$, we obtain
\begin{align*}
R_1&=\exp\bigg(\rho x-2\rho z+\rho L-\rho^2\Big(c\frac{\rho}{\beta}-s\Big)+\beta \Big(\frac{\rho^2}{2\beta}-(2\beta)^{-1/3}\gamma_1\Big)\Big(c\frac{\rho}{\beta}-s\Big)+\beta c_0\frac{\rho^2}{2\beta}\Big(c\frac{\rho}{\beta}-s\Big)\bigg)\\
&\hspace{0.2in}\int_{y_1}^{y_2}\frac{1}{1-y}\exp\bigg(\frac{\rho^2y}{2}\Big(c\frac{\rho}{\beta}-s\Big)-\frac{\beta(L-x)y}{2}\Big(c\frac{\rho}{\beta}-s\Big)-\beta c_0\frac{\rho^2}{2\beta}\Big(c\frac{\rho}{\beta}-s\Big)y+\frac{\beta^2y^3}{24}\Big(c\frac{\rho}{\beta}-s\Big)^3\\
&\hspace{0.2in}+\frac{\beta^2(1-y)^3}{12}\Big(c\frac{\rho}{\beta}-s\Big)^3-\frac{((1-c_0)\rho^2/2\beta-(2\beta)^{-1/3}\gamma_1)^2}{(c\rho/\beta-s)(1-y)}\bigg)dy.
\end{align*}
Observing that
\begin{align*}
-\frac{((1-c_0)\rho^2/2\beta-(2\beta)^{-1/3}\gamma_1)^2}{(c\rho/\beta-s)(1-y)}
&\leq  -\frac{(1-c_0)^2\rho^3}{4c\beta(1-y)}\sum_{k=0}^{\infty}\Big(\frac{s\beta}{c\rho}\Big)^k\\
&\leq -\frac{(1-c_0)^2\rho^3}{4c\beta(1-y)}-\frac{(1-c_0)^2\rho^2s}{4c^2(1-y)}-\frac{(1-c_0)^2\rho\beta s^2}{4c^3(1-y)},
\end{align*}
we get
\begin{align}\label{R1intermediate}
R_1&\lesssim \exp\bigg(\rho x-2\rho z+\rho L-\frac{c^3}{2}\frac{\rho^3}{\beta}+\frac{c^2\rho^2s}{2}-2^{-1/3}\gamma_1c\frac{\rho}{\beta^{1/3}}\bigg)\nonumber\\
&\hspace{0.2in}\int_{y_1}^{y_2}\frac{1}{1-y}\exp\bigg(\frac{\rho^3}{\beta}\Big(\frac{cy}{2}-\frac{cc_0y}{2}+\frac{c^3y^3}{24}+\frac{c^3(1-y)^3}{12}-\frac{(1-c_0)^2}{4c(1-y)}\Big)-\frac{\beta(L-x)y}{2}\Big(c\frac{\rho}{\beta}-s\Big)\nonumber\\
&\hspace{0.2in}+\rho^2s\Big(-\frac{y}{2}+\frac{c_0y}{2}-\frac{c^2y^3}{8}-\frac{c^2(1-y)^3}{4}-\frac{(1-c_0)^2}{4c^2(1-y)}\Big)+c\rho\beta s^2\Big(\frac{y^3}{8}+\frac{(1-y)^3}{4}\Big)\nonumber\\
&\hspace{0.2in}-\beta^2 s^3\Big(\frac{y^3}{24}+\frac{(1-y)^3}{12}\Big)-\frac{(1-c_0)^2\rho\beta s^2}{4c^3(1-y)}\bigg)dy.
\end{align}
Note that for $s\lesssim \beta^{-2/3}$ and all $y\in [y_1,y_2]$, 
\[
0<-2^{-1/3}\gamma_1c\frac{\rho}{\beta^{1/3}}+c\rho\beta s^2\Big(\frac{y^3}{8}+\frac{(1-y)^3}{4}\Big)=O\Big(\frac{c\rho}{\beta^{1/3}}\Big).
\]
Since $\sqrt{1-c_0}=c$, the upper bound of $R_1$ in (\ref{R1intermediate}) can further be expressed as
\begin{align}\label{R1M2}
R_1&\lesssim \exp\bigg(\rho x-2\rho z+\rho L-\frac{2c^3}{3}\frac{\rho^3}{\beta}+O\Big(\frac{c\rho}{\beta^{1/3}}\Big)\bigg)\int_{y_1}^{y_2}\frac{1}{1-y}\exp\bigg(-\frac{c\rho\beta s^2}{4(1-y)}\bigg)\nonumber\\
&\hspace{0.2in}\exp\bigg(\frac{c^3\rho^3}{4\beta}\Big(-\frac{y^3}{6}+y^2+y+1-\frac{1}{1-y}\Big)+c^2\rho^2s\Big(\frac{1}{2}-\frac{y}{2}-\frac{y^3}{8}-\frac{(1-y)^3}{4}-\frac{1}{4(1-y)}\Big)\bigg)dy.
\end{align}
Let 
\[
h(y)=\frac{1}{2}-\frac{y}{2}-\frac{y^3}{8}-\frac{(1-y)^3}{4}-\frac{1}{4(1-y)}
\]
Since for $y\in[0,1]$,
\[
h'(y)=-\frac{1}{2}-\frac{3y^2}{8}+\frac{3(1-y)^2}{4}-\frac{1}{4(1-y)^2}=\frac{1}{4}\bigg(1-\frac{1}{(1-y)^2}\bigg) +\frac{3y}{2}\bigg(\frac{y}{4}-1\bigg)\leq 0,
\]
we get $h(y)\leq h(0)=0$.
Also for all $y\in [y_1,y_2]$,
\[
\frac{c^3}{4}\Big(-\frac{y^3}{6}+y^2+y+1-\frac{1}{1-y}\Big)=-\frac{c^3y^3}{4}\Big(\frac{1}{6}+\frac{1}{1-y}\Big)\leq -\frac{c^3y_1^3}{24}.
\]
Thus,
\begin{equation}\label{R1estimate}
R_1\lesssim \exp\bigg(\rho x-2\rho z+\rho L-\frac{2c^3}{3}\frac{\rho^3}{\beta}+O\Big(\frac{c\rho}{\beta^{1/3}}\Big)-\frac{c^3y_1^3\rho^3}{24\beta}\bigg)\int_{y_1}^{y_2}\frac{1}{1-y}\exp\bigg(-\frac{c\rho\beta s^2}{4(1-y)}\bigg)dy.
\end{equation}
By (\ref{znvar1c}), we have that
\[
\frac{c^3y_1^3\rho^3}{24\beta}=\frac{c^3\rho^3}{24\beta}\Big(\frac{\beta^{-7/12}(L^*-z)^{1/4}}{t}\Big)^3\sim \frac{c^3\rho^3}{24\beta}\Big(\frac{\beta^{-7/12}(c^2\rho^2/2\beta)^{1/4}}{c\rho/\beta}\Big)^3=\frac{c^{3/2}\rho^{3/2}}{24\cdot 2^{3/4}\cdot\beta^{1/2}}\gg\frac{c\rho}{\beta^{1/3}}.
\]
Also, because $c\rho\beta s^2\gg 1$, after changing variables twice, we get
\begin{align}\label{intR1wrty}
\int_{y_1}^{y_2}\frac{1}{1-y}\exp\bigg(-\frac{c\rho\beta s^2}{4(1-y)}\bigg)dy
&\leq\int_{0}^{1}\frac{1}{y}\exp\bigg(-\frac{c\rho\beta s^2}{4y}\bigg)dy\nonumber\\
&=\int_{1}^{\infty}\frac{1}{y}\exp\bigg(-\frac{c\rho\beta s^2y}{4}\bigg)dy\nonumber\\
&\ll 1.
\end{align}
Therefore, in equation (\ref{R1estimate}), the term $O(c\rho/\beta^{1/3})$ can be absorbed into $-c^3y_1^3\rho^3/24\beta$ in the exponent and the integral can be neglected. By (\ref{znvar1c}), we conclude that
\begin{align}\label{R1}
R_1&\lesssim  \frac{\beta^{5/3}}{c\rho^5}\exp\bigg(\rho x-2\rho z+\rho L-\frac{2c^3}{3}\frac{\rho^3}{\beta}\bigg)\frac{c\rho^5}{\beta^{5/3}}\exp\bigg(-\frac{c^{3/2}\rho^{3/2}}{48\beta^{1/2}}\bigg)\nonumber\\
&\ll  \frac{\beta^{5/3}}{c\rho^5}\exp\bigg(\rho x-2\rho z+\rho L-\frac{2c^3}{3}\frac{\rho^3}{\beta}\bigg).
\end{align}
Here we want to point out that this is the only place in the proof of Lemma \ref{Large2ndmoment} where we need to use the assumption $c\gg \beta^{1/3}\rho^{-1}\log^{2/3}(\rho\beta^{-1/3})$ instead of the weaker one $c\gg \beta^{1/3}\rho^{-1}$.

Next, we estimate $R_2$. Letting $u=yt$, by similar computations as for $R_1$, we have
\begin{align}\label{R2beginning}
R_2&\lesssim \exp\bigg(\rho x-2\rho z+\rho L-\frac{c^3}{2}\frac{\rho^3}{\beta}+\frac{c^2\rho^2s}{2}-2^{-1/3}\gamma_1c\frac{\rho}{\beta^{1/3}}\bigg)\nonumber\\
&\hspace{0.2in}\int_{y_2}^{1}\frac{1}{\sqrt{(1+y)(1-y)}}\exp\bigg(\frac{c^3\rho^3}{\beta}\Big(\frac{y}{2}+\frac{y^3}{24}+\frac{(1-y)^3}{12}\Big)-c^2\rho^2s\Big(\frac{y}{2}+\frac{y^3}{8}+\frac{(1-y)^3}{4}\Big)\nonumber\\
&\hspace{0.2in}+c\rho\beta s^2\Big(\frac{y^3}{8}+\frac{(1-y)^3}{4}\Big)-\beta^2s^3\Big(\frac{y^3}{24}+\frac{(1-y)^3}{12}\Big)\bigg)\nonumber\\
&\hspace{0.2in}\times\exp\bigg(\frac{ty(1-y)}{2(1+y)}\bigg(\rho+\beta t-\frac{\beta yt}{2}-\frac{L-x}{yt}-\frac{2(L-z)}{t(1-y)}\bigg)^2-\frac{(L-x)^2}{2yt}-\frac{(L-z)^2}{t(1-y)}\bigg)dy.
\end{align}
Denote the exponent on the last line of (\ref{R2beginning}) as $A$. For $x<L$ and $y\in [y_2,1]$, we get
\begin{align*}
A&=\frac{ty(1-y)}{2(1+y)}\bigg[\bigg(\rho+\beta t-\frac{\beta yt}{2}\bigg)^2+\frac{(L-x)^2}{y^2t^2}+\frac{4(L-z)^2}{t^2(1-y)^2}-\frac{2(L-x)}{yt}\bigg(\rho+\beta t-\frac{\beta yt}{2}\bigg)\\
&\hspace{0.2in}-\frac{4(L-z)}{t(1-y)}\bigg(\rho+\beta t-\frac{\beta yt}{2}\bigg)+\frac{4(L-x)(L-z)}{t^2y(1-y)}\bigg]-\frac{(L-x)^2}{2yt}-\frac{(L-z)^2}{t(1-y)}\\
&=\frac{t}{2}\frac{y(1-y)}{1+y}\bigg(\rho+\beta t-\frac{\beta yt}{2}\bigg)^2-\frac{(L-x)^2}{t(1+y)}-\frac{(L-z)^2}{t(1+y)}-\frac{1-y}{1+y}(L-x)\bigg(\rho+\beta t-\frac{\beta yt}{2}\bigg)\\
&\hspace{0.2in}-\frac{2y}{1+y}(L-z)\bigg(\rho+\beta t-\frac{\beta yt}{2}\bigg)+\frac{2(L-x)(L-z)}{t(1+y)}\\
&\leq \frac{t}{2}\frac{y(1-y)}{1+y}\bigg(\rho+\beta t-\frac{\beta yt}{2}\bigg)^2-\frac{(L-z)^2}{t(1+y)}-\frac{2y}{1+y}(L-z)\bigg(\rho+\beta t-\frac{\beta yt}{2}\bigg)+\frac{2(L-x)(L-z)}{t(1+y)}.
\end{align*}
Recalling that $t=c\rho/\beta-s$ and $L^*-z=c^2\rho^2/2\beta$, we have
\begin{align}\label{Aint}
A&\leq \frac{y(1-y)}{2(1+y)}\Big(\frac{c\rho}{\beta}-s\Big)\bigg(\rho+\beta\Big(\frac{c\rho}{\beta}-s\Big)-\frac{\beta y}{2}\Big(\frac{c\rho}{\beta}-s\Big)\bigg)^2-\frac{(c^2\rho^2/2\beta-(2\beta)^{-1/3}\gamma_1)^2}{(1+y)(c\rho/\beta-s)}\nonumber\\
&\hspace{0.15in}-\frac{2y}{1+y}\frac{c^2\rho^2}{2\beta}\bigg(\rho+\beta\Big(\frac{c\rho}{\beta}-s\Big)-\frac{\beta y}{2}\Big(\frac{c\rho}{\beta}-s\Big)\bigg)+\frac{2(L-x)(c^2\rho^2/2\beta-(2\beta)^{-1/3}\gamma_1)}{(c\rho/\beta-s)(1+y)}.
\end{align}
Because $L-x\lesssim\beta^{-1/3}$ and $s\lesssim \beta^{-2/3}$, we see that for $y\in [y_2,1]$,
\[
\frac{2(L-x)(c^2\rho^2/2\beta-(2\beta)^{-1/3}\gamma_1)}{(c\rho/\beta-s)(1+y)}\lesssim \frac{c\rho}{\beta^{1/3}}.
\]
We also observe that $\rho\beta s^2\lesssim\rho/\beta^{1/3}$, $\beta^2s^3\ll \rho/\beta^{1/3}$ and
\[
-\frac{(c^2\rho^2/2\beta-(2\beta)^{-1/3}\gamma_1)^2}{(1+y)(c\rho/\beta-s)}\leq -\frac{c^4\rho^4/4\beta^2}{(1+y)c\rho/\beta}\sum_{k=0}^{\infty}\bigg(\frac{s}{c\rho/\beta}\bigg)^k
\leq-\frac{c^3\rho^3}{4\beta(1+y)}-\frac{c^2\rho^2s}{4(1+y)}.
\]
Therefore, equation (\ref{Aint}) implies that
\begin{align}\label{A}
A
\leq& \frac{\rho^3}{\beta}\bigg[\frac{y(1-y)}{2(1+y)}c\bigg(1+c-\frac{cy}{2}\bigg)^2-\frac{c^3}{4(1+y)}-\frac{y}{1+y}c^2\bigg(1+c-\frac{cy}{2}\bigg)\bigg]+\rho^2s\bigg[-\frac{c^2}{4(1+y)}\nonumber\\
&-\Big(1-\frac{y}{2}\Big)\Big(1+c-\frac{cy}{2}\Big)\frac{cy(1-y)}{1+y}-\frac{y(1-y)}{2(1+y)}\Big(1+c-\frac{cy}{2}\Big)+\frac{y}{1+y}\Big(1-\frac{y}{2}\Big)c^2\bigg]+O\bigg(\frac{\rho}{\beta^{1/3}}\bigg).
\end{align}
By (\ref{R2beginning}) and (\ref{A}), we have
\begin{align}\label{R2fh}
R_2&\lesssim\exp\bigg(\rho x-2\rho z+\rho L-\frac{2c^3}{3}\frac{\rho^3}{\beta}+O\Big(\frac{\rho}{\beta^{1/3}}\Big)\bigg)\int_{y_2}^{1}\frac{1}{\sqrt{(1+y)(1-y)}}\nonumber\\
&\hspace{0.2in}\times\exp\Bigg(\frac{\rho^3}{\beta}\bigg(\frac{c^3y}{4}+\frac{c^3y^2}{4}-\frac{c^3y^3}{24}+\frac{c^3}{4}+\frac{y(1-y)}{2(1+y)}c\Big(1+c-\frac{cy}{2}\Big)^2-\frac{c^3}{4(1+y)}\nonumber\\
&\hspace{0.2in}-\frac{y}{1+y}c^2\Big(1+c-\frac{cy}{2}\Big)\bigg)+\rho^2s\bigg(\frac{c^2}{2}-\frac{c^2y}{2}-\frac{c^2y^3}{8}-\frac{c^2(1-y)^3}{4}-\frac{c^2}{4(1+y)}\nonumber\\
&\hspace{0.2in}-\Big(1-\frac{y}{2}\Big)\Big(1+c-\frac{cy}{2}\Big)\frac{cy(1-y)}{1+y}-\frac{y(1-y)}{2(1+y)}\Big(1+c-\frac{cy}{2}\Big)^2+\frac{y}{1+y}\Big(1-\frac{y}{2}\Big)c^2\bigg)\Bigg)dy.
\end{align}
Define
\begin{equation}\label{phi}
\phi(y)=-\frac{2y^3}{3}+\Big(\frac{10}{3}+\frac{2}{c}\Big)y^2-\Big(\frac{2}{c^2}+\frac{6}{c}\Big)y+\frac{2}{c^2},
\end{equation}
and
\begin{equation}\label{psi}
\psi(y)=\frac{c^2y^3}{2}-\Big(\frac{5c^2}{2}+c\Big)y^2+\Big(2c^2+3c+\frac{1}{2}\Big)y-2c-\frac{1}{2}.
\end{equation}
After algebraic calculation, equation (\ref{R2fh}) is equivalent to
\begin{align}\label{R2fhshort}
R_2\lesssim& \exp\bigg(\rho x-2\rho z+\rho L-\frac{2c^3}{3}\frac{\rho^3}{\beta}+O\Big(\frac{\rho}{\beta^{1/3}}\Big)\bigg)\int_{y_2}^{1}\frac{1}{\sqrt{(1+y)(1-y)}}\nonumber\\
&\times\exp\bigg(\frac{c^3\rho^3}{\beta}\frac{y}{4(1+y)}\phi(y)+\rho^2s\frac{y}{1+y}\psi(y)\bigg)dy.
\end{align}
Below we will obtain the upper bounds for $\phi(y)$ and $\psi(y)$ in the cases $z\geq 0$ and $z<0$. 

Let us first study $\phi (y)$. Note that for every $c$, we have $\phi''(y)>0$ for all $y\in [y_2,1]$. Therefore, for $y\in[y_2,1]$, the function $\phi(y)$ reaches its maximum either at $1$ or at $y_2$. When $y=1$, we have for all $c\in (0,3/2)$,
\begin{equation}\label{phi1}
\phi(1)=\frac{8}{3}-\frac{4}{c}=-\frac{8(3/2-c)}{3c}<0.
\end{equation}
For $c\in(0,3/2)$, since $0<y_2= y_c-C_{16}c\leq 1$, after rearranging terms, we have
\begin{align*}
\phi(y_2)&=\frac{2}{3}y_2^2(1-y_c+C_{16}c)+\Big(\frac{8}{3}+\frac{2}{c}\Big)y_2^2-\Big(\frac{2}{c^2}+\frac{6}{c}\Big)(y_c-C_{16}c)+\frac{2}{c^2}\\
&\leq \frac{2}{3}y_c^2(1-y_c)+\Big(\frac{8}{3}+\frac{2}{c}\Big)y_c^2-\Big(\frac{2}{c^2}+\frac{6}{c}\Big)y_c+\frac{2}{c^2}+\frac{2C_{16}c}{3}+\frac{2C_{16}}{c}+6C_{16}\\
&=\phi(y_c)+\frac{2C_{16}c}{3}+\frac{2C_{16}}{c}+6C_{16}.
\end{align*}
Since $y_c$ satisfies (\ref{root}), we have $y_c^2=y_c(3+2/c)-2/c$. Replacing $y_c^2$ with $y_c(3+2/c)-2/c$ in the first step and replacing $y_c$ with $3/2+1/c-\Delta/c$ in the second step, we get
\begin{align*}
\phi(y_c)
&= 4y_c-\frac{2}{3c^2}y_c-\frac{8}{3c}+\frac{2}{3c^2}\\
&=6-\frac{4\Delta}{c}+\frac{4}{3c}-\frac{1}{3c^2}-\frac{2}{3c^3}+\frac{2\Delta}{3c^3}\\
&=\frac{2(9c^2-4\Delta^2)}{c(3c+2\Delta)}+\frac{4}{3c}-\frac{1}{3c^2}+\frac{2(\Delta^2-1)}{3c^3(\Delta+1)}\\
&=\frac{-8c-8}{(3c+2\Delta)c}+\frac{4}{3c}+\frac{3}{2c(\Delta+1)}-\frac{1}{3c^2}+\frac{2}{3c^2(\Delta+1)}.
\end{align*}
For $c\in (0,3/2)$, we have $\Delta>1$. According to the above two formulas, we have
\begin{align}\label{phi24}
\phi(y_2)&\leq \frac{-8c-8}{(3c+2\Delta)c}+\frac{4}{3c}+\frac{3}{2c(\Delta+1)}+\frac{2C_{16}c}{3}+\frac{2C_{16}}{c}+6C_{16}\nonumber\\
&= \frac{-24\Delta c+3c-14\Delta-48+16\Delta^2}{6c(3c+2\Delta)(\Delta+1)}+\frac{2C_{16}c}{3}+\frac{2C_{16}}{c}+6C_{16}.
\end{align}
If $z\geq 0$, then $c\in (0,1]$ and $1<\Delta\leq \sqrt{17}/2$. We have
\[
-24\Delta c+3c<-21c,
\]
and
\[-14\Delta-48+16\Delta^2=14\Delta(\Delta-1)-48+2\Delta^2\leq 14\cdot \frac{\sqrt{17}}{2}\bigg(\frac{\sqrt{17}}{2}-1\bigg)-48+2\bigg(\frac{\sqrt{17}}{2}\bigg)^2<-8.
\]
Therefore, combining the above two observations with equations (\ref{C82}), (\ref{C83}) and (\ref{phi24}), we obtain
\begin{align}\label{phiy2}
\phi(y_2)\leq \bigg(-\frac{4}{3(3+\sqrt{17})(1+\sqrt{17}/2)}+2C_{16}\bigg)\frac{1}{c}-\frac{7}{2(3+\sqrt{17})(1+\sqrt{17}/2)}+\frac{2}{3}C_{16}+6C_{16}<-\frac{C_{16}}{c}.
\end{align}
According to equations (\ref{phi1}) and (\ref{phiy2}), we get when $z\geq0$, or equivalently $c\in (0,1]$,
\begin{equation}\label{phimax1}
\max_{y\in [y_2,1]}\phi(y)\leq\max\bigg\{-\frac{8(3/2-c)}{3c},-\frac{C_{16}}{c}\bigg\}=-\frac{C_{16}}{c}.
\end{equation}
If $z<0$, then $c\in (1,3/2)$ and $\sqrt{17}/2<\Delta<11/4$. We have
\begin{align*}
-24\Delta c+3c-14\Delta-48+16\Delta^2
&<-21\Delta c+14\Delta(\Delta-1)+2\Delta^2-48\\
&<-21\cdot\frac{\sqrt{17}}{2}\cdot 1+14\cdot\frac{11}{4}\Big(\frac{11}{4}-1\Big)+2\cdot\frac{11}{4}\cdot\frac{11}{4}-48\\
&<-8.
\end{align*}
Therefore, combining the above two observations with equations (\ref{C84}) and (\ref{phi24}), since $c\in (0,3/2)$, we obtain
\[
\phi(y_2)\leq \frac{-8}{9(3\cdot \frac{3}{2}+2\cdot\frac{11}{4})(\frac{11}{4}+1)}+9C_{16}<-C_{16}.
\]
According to equations (\ref{phi1}) and (\ref{phiy2}), we get when $z\leq0$, or equivalently $c\in (1,3/2)$,
\begin{equation}\label{phimax2}
\max_{y\in [y_2,1]}\phi(y)\leq\max\bigg\{-\frac{8(3/2-c)}{3c},-C_{16}\bigg\}.
\end{equation}

We next study $\psi(y)$. Let us first consider the case $z\geq 0$, or equivalently, $c\in (0,1]$. For all $y\in [y_2,1]$, we have
\begin{equation}\label{psi1}
\psi(y)\leq \frac{c^2y^2}{2}-\Big(\frac{5c^2}{2}+c\Big)y^2+(2c^2+c)y=c(2c+1)y(1-y)\leq 3c.
\end{equation}
If $z<0$, we claim that for all $y\in [y_2,1]$ and $c\in (0,3/2)$, we have
\begin{equation}\label{psiclaim}
\psi(y)\leq -3c^2\phi(y)/4.
\end{equation}
Indeed, for every $y$, we can view $-3c^2\phi(y)/4-\psi(y)$ as a quadratic function of $c$:
\[
-\frac{3c^2\phi(y)}{4}-\psi(y)=-2yc^2+c\bigg(-\frac{1}{2}y^2+\frac{3}{2}y+2\bigg)+y-1=: \varphi(c).
\]
Note that for every $y\in (0,1]$, the quadratic function $\varphi(c)$ for $c\in [1,3/2]$ reaches its minimum at either $c=1$ or $c=3/2$. Since for all $y\in (0,1]$, we have
\[
\varphi(1)=-\frac{1}{2}(y^2-y-2)>0,\qquad \varphi\Big(\frac{3}{2}\Big)=-\frac{1}{4}(3y^2+5y-8)\geq 0,
\]
the claim follows.

Now it remains to apply the upper bounds of $\phi(y)$ and $\psi(y)$ in (\ref{R2fhshort}). By (\ref{u2C8}), for all $y\in [y_2,1]$, we have 
\begin{equation} \label{y/(1+y)}
C_{16}/2\leq y/(1+y)\leq 1.
\end{equation}
When $z\geq 0$, combining (\ref{R2fhshort}) with (\ref{phimax1}), (\ref{psi1}) and (\ref{y/(1+y)}), we obtain
\begin{equation}\label{R2f}
R_2\lesssim \exp\bigg(\rho x-2\rho z+\rho L-\frac{2c^3}{3}\frac{\rho^3}{\beta}-\frac{C_{17}^2}{8}\frac{c^2\rho^3}{\beta}+3c\rho^2s+O\Big(\frac{\rho}{\beta^{1/3}}\Big)\bigg).
\end{equation}
Recalling that $c\gg \beta^{1/3}/\rho$ and $s\lesssim \beta^{-2/3}$, we have
\[
\frac{c^2\rho^3}{\beta}\gg c\rho^2s,\quad \frac{c^2\rho^3}{\beta}\gg \frac{\rho}{\beta^{1/3}}, \quad\exp\bigg(-\frac{C_{17}^2}{16}\frac{c^2\rho^3}{\beta}\bigg)\ll\frac{\beta^{5/3}}{c\rho^5}.
\]
Consequently, 
\begin{equation}\label{R2}
R_2\lesssim \frac{\beta^{5/3}}{c\rho^5}\exp\bigg(\rho x-2\rho z+\rho L-\frac{2c^3}{3}\frac{\rho^3}{\beta}\bigg).
\end{equation}
If $z<0$, since $c^3\rho^3/\beta\gg c^2\rho^2 s$, according to (\ref{R2fhshort}), (\ref{phimax2}), (\ref{psiclaim}) and (\ref{y/(1+y)}), we have for $n$ sufficiently large
\begin{align*}
R_2&\lesssim \exp\bigg(\rho x-2\rho z+\rho L-\frac{2c^3}{3}\frac{\rho^3}{\beta}+O\Big(\frac{\rho}{\beta^{1/3}}\Big)\bigg)\int_{y_2}^{1}\frac{1}{\sqrt{(1+y)(1-y)}}\exp\bigg(\frac{c^3\rho^3}{\beta}\frac{y\phi(y)}{8(1+y)}\bigg)dy\\
&\lesssim  \exp\bigg(\rho x-2\rho z+\rho L-\frac{2c^3}{3}\frac{\rho^3}{\beta}-\frac{C_{16}}{16}\frac{c^3\rho^3}{\beta}\min\bigg\{\frac{8(3/2-c)}{3c},C_{16}\bigg\}+O\Big(\frac{\rho}{\beta^{1/3}}\Big)\bigg).
\end{align*}
By (\ref{znvar1c}), since $c>1$, we have
\[
\frac{c^3\rho^3}{\beta}\min\bigg\{\frac{8(3/2-c)}{3c},C_{16}\bigg\}\gg\frac{\rho}{\beta^{1/3}}, \quad\exp\bigg(-\frac{c^3\rho^3}{\beta}\min\bigg\{\frac{8(3/2-c)}{3c},C_{16}\bigg\}\bigg)\ll\frac{\beta^{5/3}}{c\rho^5}.
\]
Therefore, equation (\ref{R2}) also holds when $z<0$. 

Finally, combining (\ref{R1}) and (\ref{R2}), equation (\ref{secondmomentpartIIc}) is proved and the lemma follows.
\qedwhite
\\

\noindent\textit{Proof of Lemma \ref{2ndmomentgap}.} Recall that in the proof of Lemma \ref{Large2ndmoment}, the only place where we used the assumption $L^*-z\gg  \beta^{-1/3}\log^{4/3} (\rho/\beta^{1/3})$ is equation (\ref{R1}). Thus to prove Lemma \ref{2ndmomentgap}, it is sufficient to prove that for $z$ satisfying $\beta^{-1/3}\ll L^*-z\lesssim \beta^{-1/3}\log^{4/3} (\rho/\beta^{1/3})$, or equivalently, $\beta^{1/3}\rho^{-1}\ll c\lesssim \beta^{1/3}\rho^{-1}\log^{2/3}(\rho/\beta^{1/3})$, we have
\begin{equation}\label{2ndsmallgapgoal}
I_2':=\int_{u_1}^{u_2}\int_{-\infty}^{L}p_u^L(x,r)\Big(p_{t-u}^{L}(r,z)\Big)^2drdu\lesssim \frac{\beta^{5/3}}{c\rho^5}\exp\bigg(\rho x-2\rho z+\rho L-\frac{2c^3}{3}\frac{\rho^3}{\beta}\bigg).
\end{equation}
The portion of the double integral in $I_2$ for which $u_2\leq u\leq t$ has been dealt with in Lemma \ref{Large2ndmoment}.
By (\ref{gamma1}), we can choose constant $C_{17}>0$ such that
\begin{equation}\label{C9}
2^{-1/3}\gamma_1+C_{17}+1<-\frac{1}{2}
\end{equation}
According to equations (\ref{tL-z}) and (\ref{y2}), and the fourth equality in (\ref{1-yc}), we get
\begin{equation}\label{t-u2}
t-u_2=t(1-y_2)=t(1-y_c+C_{16}c)\leq t(z)\bigg(\frac{9c}{4(\Delta+1)}+C_{16}c\bigg)\asymp \frac{c^2\rho}{\beta}\ll\beta^{-2/3}.
\end{equation}
Thus $t-u_2< 2\beta^{-2/3}$ for $n$ large enough. Therefore, for $n$ large enough, we can write
\[
I_2'=P_1+P_2+P_3,
\]
where $P_1$ is the part of the double integral for which $L-C_{17}\beta^{-1/3}\leq r\leq L$ and $u_1\leq u\leq t-2\beta^{-2/3}$, $P_2$ is the part of the double integral for which $L-C_{17}\beta^{-1/3}\leq r\leq L$ and $t-2\beta^{-2/3}<u\leq u_2$, and $P_3$ is the part of the double integral for which $r<L-C_{17}\beta^{-1/3}$ and $u_1\leq u\leq u_2$.

To bound $P_1$, we are going to bound $p_u^L(x,r)$ by (\ref{lemma2.6and2.7}) and $p_{t-u}^L(r,z)$ by (\ref{lemma2.8}). We get
\begin{align*}
P_1
&\lesssim \int_{u_1}^{t-2\beta^{-2/3}}\int_{L-C_{17}\beta^{-1/3}}^{L}\frac{(L-x)(L-r)}{u^{3/2}}\exp\bigg(\rho x-\rho r-\frac{(x-r)^2}{2u}-\frac{\rho^2u}{2}+\beta Lu\bigg)\\
&\qquad \times\frac{\beta^{2/3}(L-r)^2}{t-u} \bigg(\max\bigg\{1,\frac{1}{\beta^{1/3}(t-u)}\bigg(L-z-\frac{\beta(t-u)^2}{2}\bigg)\bigg\}\bigg)^2 \\
&\qquad \times \exp\bigg(2\rho r-2\rho z-\frac{(r-z)^2}{t-u}-\rho^2(t-u)+\beta(r+z)(t-u)+\frac{\beta^2(t-u)^3}{12}\\
&\qquad \qquad+\frac{1}{\beta^{1/3}(t-u)}\bigg(L-z-\frac{\beta(t-u)^2}{2}\bigg)\bigg)drdu.
\end{align*}
Note that $L-x\lesssim \beta^{-1/3}$ and $t=t(z)-s$. Interchanging the roles of $r$ and $L-r$, we have
\begin{align}\label{M11}
P_1&\lesssim\int_{u_1}^{t-2\beta^{-2/3}} \hspace{-0.1in}\frac{\beta^{1/3}}{u^{3/2}(t-u)} \bigg(\max\bigg\{1,\frac{1}{\beta^{1/3}(t-u)}\bigg(L-z-\frac{\beta(t-u)^2}{2}\bigg)\bigg\}\bigg)^2\exp\bigg(\rho x-2\rho z+\rho L\nonumber\\
&\hspace{0.15in}-\frac{\rho^2u}{2}+\beta Lu-\rho^2\big(t(z)-s-u\big)+\frac{\beta^2(t(z)-s-u)^3}{12}+\frac{1}{\beta^{1/3}(t-u)}\Big(L-z-\frac{\beta(t-u)^2}{2}\Big)\bigg)\nonumber\\
&\hspace{0.15in}\int_{0}^{C_{17}\beta^{-1/3}}r^3\exp\bigg(-\rho r-\frac{(L-x-r)^2}{2u}-\frac{(L-r-z)^2}{t-u}+\beta(L-r+z)(t(z)-s-u)\bigg)drdu.
\end{align}
We first estimate the term $(L-z-\beta(t-u)^2/2)/(\beta^{1/3}(t-u))$. For $u\geq u_1$, we see that
\[
\frac{\beta(u+s)^2}{2}\geq \frac{\beta u_1^2}{2}=\frac{\beta}{2}\big(\beta^{-7/12}(L^*-z)^{1/4}\big)^2\asymp \frac{c\rho}{\beta^{2/3}}\gg \beta^{-1/3}\asymp (2\beta)^{-1/3}|\gamma_1|.
\]
Thus by (\ref{1/beta^{1/3}(t-u)}), we have for $n$ large enough
\begin{equation}\label{lem2.8bound}
\frac{1}{\beta^{1/3}(t-u)}\bigg(L-z-\frac{\beta(t-u)^2}{2}\bigg)\leq \frac{\beta t(z)(u+s)}{\beta^{1/3}(t-u)}=\frac{\beta^{2/3}t(z)(u+s)}{t-u}.
\end{equation}
Furthermore, we note that 
\[
\frac{\beta^{2/3}t(z)(u+s)}{t-u}\geq \frac{\beta^{2/3}t(z)u_1}{t(z)}=\beta^{2/3}u_1\asymp\frac{c^{1/2}\rho^{1/2}}{\beta^{1/6}}\gg 1.
\]
Moreover, we will upper bound the term $-(L-r-z)^2/(t-u)$ by  (\ref{geometric}). By  (\ref{geometric}), (\ref{M11}) and (\ref{lem2.8bound}), after rearranging terms, we get for $n$ large
\begin{align*}
P_1&\lesssim \exp\bigg(\rho x-2\rho z+\rho L-\rho^2t(z)+\beta L^* t(z) +\beta zt(z)+\frac{\beta^2t(z)^3}{12}-\frac{(L^*-z)^2}{t(z)}\bigg)\\
&\hspace{0.15in}\int_{u_1}^{t-2\beta^{-2/3}}\frac{\beta^{5/3}t(z)^2(u+s)^2}{u^{3/2}(t-u)^3}\exp\bigg(-\frac{\rho^2u}{2}+\beta Lu+\rho^2(u+s)-\frac{\beta^2(u+s)^3}{12}+\frac{\beta^2t(z)(u+s)^2}{4}\\
&\hspace{0.15in}-\frac{\beta^2t(z)^2(u+s)}{4}+\frac{\beta^{2/3}t(z)(u+s)}{t-u}-\beta L(u+s)-2^{-1/3}\beta^{2/3}\gamma_1t(z)-\beta z(u+s)\\
&\hspace{0.15in}-\frac{(L^*-z)^2(u+s)}{t(z)^2}-\frac{(L^*-z)^2(u+s)^2}{t(z)^3}-\frac{(L^*-z)^2(u+s)^3}{t(z)^4}\bigg)\int_{0}^{C_{17}\beta^{-1/3}}r^3\exp\bigg(-\rho r\\
&\hspace{0.15in}-\beta r(t(z)-s-u)-\frac{(L-x-r)^2}{2u}+\frac{2(L^*-z)((2\beta)^{-1/3}\gamma_1+r)}{t-u}\bigg)drdu.
\end{align*}
Notice that since $t(z)=\sqrt{2/\beta}\sqrt{L^*-z}$, $L=L^*-(2\beta)^{-1/3}\gamma_1$ and $s\lesssim \beta^{-2/3}$, we observe that
\begin{align*}
&-\frac{\rho^2 u}{2}+\beta L u+\rho^2(u+s)-\frac{\beta^2t(z)^2(u+s)}{4}-\beta L(u+s)-\beta z(u+s)-\frac{(L^*-z)^2(u+s)}{t(z)^2}\\
&\hspace{0.2in}=\frac{\rho^2u}{2}+\rho^2s-\beta Ls-\beta z(u+s)-\beta (L^*-z)(u+s)\\
&\hspace{0.2in}=\Big(\frac{\rho^2}{2}-\beta L\Big)s\\
&\hspace{0.2in}=O(1).
\end{align*}
Also
\[
\frac{\beta^2t(z)(u+s)^2}{4}-\frac{(L^*-z)^2(u+s)^2}{t(z)^3}=0
\]
and
\[
 -\frac{\beta^2(u+s)^3}{12}-\frac{(L^*-z)^2(u+s)^3}{t(z)^4}=-\frac{\beta^2(u+s)^3}{3}.
\]
By the above four equations and (\ref{2c^3/3}), we can further bound $P_1$ as follows:
\begin{align*}
P_1&\lesssim \exp\bigg(\rho x-2\rho z+\rho L-\frac{4\sqrt{2\beta}}{3}(L^*-z)^{3/2}\bigg)\int_{u_1}^{t-2\beta^{-2/3}}\frac{\beta^{5/3}t(z)^2(u+s)^2}{u^{3/2}(t-u)^3}\\
&\hspace{0.4in}\times\exp\bigg(-\frac{\beta^2(u+s)^3}{3}+\frac{\beta^{2/3}t(z)(u+s)}{t-u}-2^{-1/3}\beta^{2/3}\gamma_1t(z)\\
&\hspace{0.6in}+\frac{2(L^*-z)(2^{-1/3}\gamma_1+C_{17})\beta^{-1/3}}{t-u}\bigg)\times\int_{0}^{C_{17}\beta^{-1/3}}r^3e^{-\rho r}drdu.
\end{align*}
Recall that $c\gg \beta^{1/3}/\rho$. By equation (\ref{tL-z}), for $n$ sufficiently large, for all $u_1\leq u\leq t-2\beta^{-2/3}$, we have
\[
-\frac{\beta^2(u+s)^3}{3}-2^{-1/3}\beta^{2/3}\gamma_1t(z)\leq -\frac{\beta^2u_1^3}{3}-\frac{\gamma_1c\rho}{2^{1/3}\beta^{1/3}}=-\frac{c^{3/2}\rho^{3/2}}{3\cdot 2^{3/4}\beta^{1/2}}-\frac{\gamma_1c\rho}{2^{1/3}\beta^{1/3}}\leq -\frac{c^{3/2}\rho^{3/2}}{6\beta^{1/2}}.
\]
Also, by equations (\ref{tL-z}) and (\ref{C9}), since $u+s\leq t(z)=c\rho/\beta$, we have for all $u_1\leq u\leq t-2\beta^{-2/3}$,
\begin{align*}
\frac{\beta^{2/3}t(z)(u+s)}{t-u}+\frac{2(L^*-z)(2^{-1/3}\gamma_1+C_{17})\beta^{-1/3}}{t-u}
&=\frac{c\rho}{\beta^{1/3}(t-u)}\bigg((u+s)+\frac{c\rho}{\beta}\big(2^{-1/3}\gamma_1+C_{17}\big)\bigg)\\
&\leq \frac{c^2\rho^2}{\beta^{4/3}(t-u)}\big(1+2^{-1/3}\gamma_1+C_{17}\big)\\
&\leq -\frac{c^2\rho^2}{2\beta^{4/3}(t-u)}.
\end{align*}
Combining the above three equations with (\ref{tL-z}), after some standard calculations, we get
\begin{align*}
P_1&\lesssim \frac{c^2}{\rho^2\beta^{1/3}}\exp\bigg(\rho x-2\rho z+\rho L-\frac{2c^3}{3}\frac{\rho^3}{\beta}-\frac{c^{3/2}\rho^{3/2}}{6\beta^{1/2}}\bigg)\\
&\hspace{0.2in}\int_{u_1}^{t-2\beta^{-2/3}}\frac{(u+s)^2}{u^{3/2}(t-u)^3}\exp\Big(-\frac{c^2\rho^2}{2\beta^{4/3}(t-u)}\Big)du.
\end{align*}
Note that for $n$ large, we have $u+s\leq 2u$ for all $u_1\leq u\leq t-2\beta^{-2/3}$. Let $v=t-u$. The integral in the previous equation can be upper bounded by
\begin{align}\label{1/v^3}
4\int_{u_1}^{t-2\beta^{-2/3}}\frac{u^{1/2}}{(t-u)^3}\exp\Big(-\frac{c^2\rho^2}{2\beta^{4/3}(t-u)}\Big)du
&\leq \frac{4c^{1/2}\rho^{1/2}}{\beta^{1/2}}\int_{2\beta^{-2/3}}^{t-u_1}v^{-3}\exp\Big(-\frac{c^2\rho^2}{2\beta^{4/3}v}\Big)dv\nonumber\\
&\leq  \frac{4c^{1/2}\rho^{1/2}}{\beta^{1/2}}\int_{0}^{\infty}v^{-3}\exp\Big(-\frac{c^2\rho^2}{2\beta^{4/3}v}\Big)dv\nonumber\\
&=\frac{4c^{1/2}\rho^{1/2}}{\beta^{1/2}}\cdot\frac{4\beta^{8/3}}{c^4\rho^4}.
\end{align} 
Combining the above two formulas, because $c\rho/\beta^{1/3}\gg 1$, we get
\begin{align}\label{M1}
P_1&\lesssim\frac{\beta^{11/6}}{c^{3/2}\rho^{11/2}}\exp\bigg(\rho x-2\rho z+\rho L-\frac{2c^3}{3}\frac{\rho^3}{\beta}-\frac{c^{3/2}\rho^{3/2}}{6\beta^{1/2}}\bigg)\nonumber\\
&\ll \frac{\beta^{5/3}}{c\rho^5}\exp\bigg(\rho x-2\rho z+\rho L-\frac{2c^3}{3}\frac{\rho^3}{\beta}\bigg).
\end{align}

We next estimate $P_2$. We are going to bound both $p_u^L(x,r)$ and $p_{t-u}^L(r,z)$ by (\ref{lemma2.6and2.7}). Interchanging the roles of $r$ and $L-r$, we get
\begin{align}\label{M21}
P_2&\lesssim (L-x)(L-z)^2\exp\bigg(\rho x-2\rho z+\rho L-\rho^2 t+2\beta Lt\bigg) \nonumber \\
&\qquad \times \int_{t-2\beta^{-2/3}}^{u_2}\frac{1}{u^{3/2}(t-u)^3}\exp\bigg(\frac{\rho^2u}{2}-\beta Lu\bigg)\nonumber\\
&\qquad \qquad  \int_{0}^{C_{17}\beta^{-1/3}}r^3\exp\bigg(-\rho r-\frac{(L-x-r)^2}{2u}-\frac{(L-r-z)^2}{t-u}\bigg)drdu.
\end{align}
By (\ref{C9}), we have for all $0\leq r\leq C_{17}\beta^{-1/3}$, 
\[
-\frac{(L-r-z)^2}{t-u}=-\frac{(L^*-z)^2}{t-u}-\frac{((2\beta)^{-1/3}\gamma_1+r)^2}{t-u}+\frac{2(L^*-z)((2\beta)^{-1/3}\gamma_1+r)}{t-u}\leq -\frac{(L^*-z)^2}{t-u}.
\]
Thus for $t-2\beta^{-2/3}\leq u\leq u_2$, the inner integral in (\ref{M21}) can be upper bounded by
\begin{align*}
\int_{0}^{\infty}r^3\exp\bigg(-\rho r-\frac{(L^*-z)^2}{t-u}\bigg)dr &\lesssim 
\frac{1}{\rho^4}\exp\Big(-\frac{(L^*-z)^2}{t-u}\Big) \\
&\leq \frac{1}{\rho^4}\exp\bigg(-\frac{(L^*-z)^2}{2(t-u)}-\frac{(L^*-z)^2}{2\cdot 2\beta^{-2/3}}\bigg).
\end{align*}
Then equation (\ref{M21}) becomes
\begin{align*}
P_2&\lesssim \frac{(L-x)(L-z)^2}{\rho^4}\exp\bigg(\rho x-2\rho z+\rho L-\rho^2t+2\beta Lt-\frac{\beta^{2/3}(L^*-z)^2}{4}\bigg)\\
&\hspace{0.2in}\int_{t-2\beta^{-2/3}}^{u_2}\frac{1}{u^{3/2}(t-u)^3}\exp\bigg(\frac{\rho^2u}{2}-\beta Lu-\frac{(L^*-z)^2}{2(t-u)}\bigg)du.
\end{align*}
Expressing $L-z$ and $t$ in terms of $c$, since $L-x\lesssim\beta^{-1/3}$, we get
\begin{align*}
P_2&\lesssim \frac{c^4}{\beta^{7/3}}\exp\bigg(\rho x-2\rho z+\rho L-2^{2/3}\gamma_1\frac{c\rho}{\beta^{1/3}}-\frac{c^4\rho^4}{16\beta^{4/3}}\bigg)\\
&\qquad \qquad \times \int_{t-2\beta^{-2/3}}^{u_2}\frac{1}{u^{3/2}(t-u)^3}\exp\bigg(-\frac{c^4\rho^4}{8\beta^2(t-u)}\bigg)du.
\end{align*}
By applying the same argument as in (\ref{1/v^3}), the integral in the previous equation can be upper bounded by 
\begin{align*}
\frac{1}{(t-2\beta^{-2/3})^{3/2}}\int_{t-2\beta^{-2/3}}^{u_2}\frac{1}{(t-u)^3}\exp\bigg(-\frac{c^4\rho^4}{8\beta^2(t-u)}\bigg)
&\lesssim \frac{\beta^{3/2}}{c^{3/2}\rho^{3/2}}\int_{0}^{\infty}\frac{1}{v^3}\exp\bigg(-\frac{c^4\rho^4}{8\beta^2v}\bigg)dv\\
&\asymp\frac{\beta^{11/2}}{c^{19/2}\rho^{19/2}}.
\end{align*}
Combining the above two equations, since $c\rho/\beta^{1/3}\gg 1$, we have
\begin{align}\label{M2}
P_2&\lesssim\frac{\beta^{19/6}}{c^{11/2}\rho^{19/2}}\exp\bigg(\rho x-2\rho z+\rho L-2^{2/3}\gamma_1\frac{c\rho}{\beta^{1/3}}-\frac{c^4\rho^4}{16\beta^{4/3}}\bigg)\nonumber\\
&\ll \frac{\beta^{5/3}}{c\rho^5}\exp\bigg(\rho x-2\rho z+\rho L-\frac{2c^3}{3}\frac{\rho^3}{\beta}\bigg).
\end{align}

It remains to estimate $P_3$. We are going to bound both $p_u^L(x,r)$ and $p_{t-u}^L(r,z)$ by (\ref{density2.8}). By a similar calculation as in (\ref{I2firstformula}), we get
\begin{align*}
P_3& \lesssim \exp\bigg(\rho x-2\rho z+\rho L-\rho^2 t+\beta Lt+\beta zt\bigg)\\
&\hspace{0.2in}\times\int_{u_1}^{u_2}u^{-1/2}(t-u)^{-1}\exp\bigg(\frac{\rho^2 u}{2}+\frac{\beta xu}{2}-\beta zu-\frac{\beta Lu}{2}+\frac{\beta^2 u^3}{24}+\frac{\beta^2(t-u)^3}{12}\bigg)\\
&\hspace{0.5in}\int_{C_{17}\beta^{-1/3}}^{\infty}\exp\bigg(-r\bigg(\rho+\beta t-\frac{\beta u}{2}\bigg)-\frac{(L-x-r)^2}{2u}-\frac{(L-z-r)^2}{t-u}\bigg)drdu\\
&\leq \exp\bigg(\rho x-2\rho z+\rho L-\rho^2 t+\beta Lt+\beta zt-\frac{C_{17}\rho}{\beta^{1/3}}\bigg)\\
&\hspace{0.2in}\times\int_{u_1}^{u_2}(t-u)^{-1}\exp\bigg(\frac{\rho^2 u}{2}+\frac{\beta xu}{2}-\beta zu-\frac{\beta Lu}{2}+\frac{\beta^2 u^3}{24}+\frac{\beta^2(t-u)^3}{12}\bigg)du.
\end{align*}
Note that the above upper bound is very similar to $R_1$ defined in (\ref{defR1R2}). Therefore, by carrying out the same calculation as for $R_1$ in (\ref{R1M2}), we obtain
\begin{align}\label{M31}
P_3&\lesssim \exp\bigg(\rho x-2\rho z+\rho L-\frac{2c^3}{3}\frac{\rho^3}{\beta}-\frac{C_{17}\rho}{\beta^{1/3}}+O\Big(\frac{c\rho}{\beta^{1/3}}\Big)\bigg)\nonumber\\
&\quad \times \int_{y_1}^{y_2}\frac{1}{1-y}\exp\bigg(\frac{c^3\rho^3}{4\beta}\Big(-\frac{y^3}{6}+y^2+y+1\Big)+c^2\rho^2s\Big(\frac{1}{2}-\frac{y}{2}-\frac{y^3}{8}-\frac{(1-y)^3}{4}\Big)\bigg)dy,
\end{align}
where $y_1=u_1/t$ and $y_2=u_2/t$. We see that for $c\lesssim \beta^{1/3}\rho^{-1}\log^{2/3}(\rho/\beta^{1/3})$ and $y\in [y_1,y_2]$,
\[
\frac{c^3\rho^3}{4\beta}\Big(-\frac{y^3}{6}+y^2+y+1\Big)\leq \frac{3c^3\rho^3}{4\beta}\ll \frac{\rho}{\beta^{1/3}},
\]
\[
c^2\rho^2s\Big(\frac{1}{2}-\frac{y}{2}-\frac{y^3}{8}-\frac{(1-y)^3}{4}\Big)\lesssim \frac{c^2\rho^2}{\beta^{2/3}}\ll \frac{\rho}{\beta^{1/3}}.
\]
Thus the integral in (\ref{M31}) can be bounded by
\begin{equation}\label{M3int1}
\exp\bigg(o\Big(\frac{\rho}{\beta^{1/3}}\Big)\bigg)\int_{y_1}^{y_2}\frac{1}{1-y}dy\leq \exp\bigg(o\Big(\frac{\rho}{\beta^{1/3}}\Big)\bigg)\int_{1-y_2}^{1}\frac{1}{v}dv=\exp\bigg(o\Big(\frac{\rho}{\beta^{1/3}}\Big)\bigg)\log\Big(\frac{1}{1-y_2}\Big).
\end{equation}
According to (\ref{Deltadef}) and (\ref{u2}), we see that $1-y_2=1-y_c+C_{16}c\geq C_{16}c$. Thus, since $c\gg\beta^{1/3}\rho^{-1}$, we have
\begin{equation}\label{M3int2}
\log\Big(\frac{1}{1-y_2}\Big)\leq \log\Big(\frac{1}{C_{16}c}\Big)\lesssim\log\Big(\frac{\rho}{\beta^{1/3}}\Big).
\end{equation}
Combining (\ref{M3int1}) and (\ref{M3int2}) with (\ref{M31}), since $c\lesssim\beta^{1/3}\rho^{-1}\log^{2/3}(\rho/\beta^{1/3})$, we get
\begin{align}\label{M3}
P_3&\lesssim \log\Big(\frac{\rho}{\beta^{1/3}}\Big)\exp\bigg(\rho x-2\rho z+\rho L-\frac{2c^3}{3}\frac{\rho^3}{\beta}-\frac{C_{17}\rho}{\beta^{1/3}}+O\Big(\frac{c\rho}{\beta^{1/3}}\Big)+o\Big(\frac{\rho}{\beta^{1/3}}\Big)\bigg)\nonumber\\
&\ll \frac{\beta^{5/3}}{c\rho^5}\exp\bigg(\rho x-2\rho z+\rho L-\frac{2c^3}{3}\frac{\rho^3}{\beta}\bigg).
\end{align}
Finally, equation (\ref{2ndsmallgapgoal}) follows from (\ref{M1}), (\ref{M2}) and (\ref{M3}) and the lemma follows.

\end{document}